\documentclass[]{article}

\usepackage[hmarginratio=1:1,top=32mm,columnsep=20pt]{geometry} % Document margins
\usepackage{geometry}
\usepackage[hang, small,labelfont=bf,up,textfont=it,up]{caption} % Custom captions under/above floats in tables or figures
\usepackage{booktabs} % Horizontal rules in tables
\usepackage{amssymb, epsfig,amssymb, latexsym}
\usepackage{bm}
\usepackage{amsfonts,psfrag,amsmath,bbm,color,url}
\usepackage{graphicx}
\usepackage{subfig}

\usepackage{amsfonts}
 \usepackage{amsmath}
 \usepackage{xcolor}
\usepackage{tikz}
\usepackage{pgfplots}
\usepackage{comment}
\usepackage{bm}
\usepackage{mathabx}
\usepackage{algorithm}
\usepackage{algpseudocode}
\DeclareMathAlphabet{\mathbbx}{U}{bbold}{m}{n}

\usepackage{cases}
\usepackage[output-decimal-marker={,}]{siunitx}
\usepackage{lineno,hyperref}
\definecolor{vargreen}{rgb}{0.0, 0.5, 0.0}
\definecolor{navyblue}{rgb}{0.0, 0.0, 0.5}

\newcommand{\blue}[1]{{\color{black} #1}}
\definecolor{mediumorchid}{rgb}{0.73, 0.33, 0.83}
\definecolor{crimson}{rgb}{0.86, 0.08, 0.24}
\definecolor{lightseagreen}{rgb}{0.13, 0.7, 0.67}
\definecolor{royalblue}{rgb}{0.25, 0.41, 0.88}
\definecolor{hotpink}{rgb}{1.0, 0.41, 0.71}
\definecolor{magenta}{rgb}{1.0, 0.0, 1.0}
\definecolor{goldenrod}{rgb}{0.85, 0.65, 0.13}

\definecolor{plum(traditional)}{rgb}{0.56, 0.27, 0.52}
\usepackage{appendix}

\usepackage[export]{adjustbox}
\usepackage{mdframed}
\usepackage{array}
\usepackage{tikz}
\usetikzlibrary{shapes}
\usetikzlibrary{plotmarks}

\usepackage{xpatch}
\usepackage{nomencl} % Nomenclature package
\makenomenclature
\RequirePackage{ifthen}
\let\oldref\ref
\renewcommand{\ref}[1]{(\oldref{#1})}

\DeclareMathAlphabet\mathbfcal{OMS}{cmsy}{b}{n}

\xpatchcmd{\thenomenclature}{%
  \section{\nomname}% Look for `\section*... etc.
}{}{\typeout{Success}}{\typeout{Failure}}

\renewcommand{\nomname}{B.$\quad$List of abbreviations and symbols}

\usepackage{ifthen}
\renewcommand{\nomgroup}[1]{%
  \ifthenelse{\equal{#1}{A}}{\item[\textbf{Abbreviations}]}{%
    \ifthenelse{\equal{#1}{G}}{\item[\textbf{Symbols}]}{%
      \ifthenelse{\equal{#1}{C}}{\item[\textbf{Abbreviations}]}{%
        \ifthenelse{\equal{#1}{S}}{\item[\textbf{Subscripts}]}{%
          \ifthenelse{\equal{#1}{Z}}{\item[\textbf{Mathematical Symbols}]}{}
        }% matches mathematical symbols
      }% matches Subscripts
    }% matches Abbreviations
  }% matches Greek Symbols
}% matches Roman Symbols

\usepackage{multirow}
\usepackage{comment}
\usepackage{pifont}% http://ctan.org/pkg/pifont

\newcommand{\highlight}[1]{\colorbox{gray!30}{\parbox{\dimexpr\linewidth-2\fboxsep\relax}{#1}}}

\title{\textbf{Data-driven Closure Strategies for Parametrized Reduced Order Models via Deep Operator Networks}}

\date{ }

\author{Anna Ivagnes  \\ \small SISSA, International School for Advanced Studies, \\ \small Mathematics Area, mathLab, Trieste, Italy. \\ \small  \href{mailto:aivagnes@sissa.it}{aivagnes@sissa.it} \normalsize \and Giovanni Stabile \\ \small Sant'Anna School of Advanced Studies
\\ \small  The Biorobotics Institute, %V.le R. Piaggio 34, 56025,
Pontedera, Pisa, Italy. \\ \small  \href{mailto:giovanni.stabile@santannapisa.it}{giovanni.stabile@santannapisa.it} \normalsize \and Gianluigi Rozza  \\ \small  SISSA, International School for Advanced Studies, \\\small  Mathematics Area, mathLab, Trieste, Italy. \\ \small \href{mailto:grozza@sissa.it}{grozza@sissa.it}}

\begin{document}

\maketitle

\nomenclature{$\text{FOM}$}{Full Order Model}
\nomenclature{$\text{ROM}$}{Reduced Order Model}
\nomenclature{$\text{NSE}$}{Navier--Stokes Equations}
\nomenclature{$\text{POD}$}{Proper Orthogonal Decomposition}
\nomenclature{$\text{PCA}$}{Principal Component Analysis}
\nomenclature{$\text{SVD}$}{Singular Value Decomposition}
\nomenclature{$\text{PPE-ROM}$}{ROM with pressure Poisson equation stabilisation}
\nomenclature{$\text{EV-ROM}$}{Eddy-Viscosity Reduced Order Model}
\nomenclature{$\text{DD-EV-ROM}$}{Data-driven Eddy-Viscosity Reduced Order Model (with standard networks' training)}
\nomenclature{$\text{DD-EV-ROM}^{\star}$}{Data-driven Eddy-Viscosity Reduced Order Model (with \emph{advanced} networks' training)}
\nomenclature{$\text{RANS}$}{Reynolds--Averaged Navier--Stokes}
\nomenclature{$\text{URANS}$}{Unsteady Reynolds--Averaged Navier--Stokes}
\nomenclature{$\text{MSE}$}{Mean Squared Error}
\nomenclature{$\text{NN}$}{Neural Network}
\nomenclature{$\text{LSTM}$}{Long--Short Term Memory neural network}
\nomenclature{$\text{SinNN}$}{Neural network with ad-hoc architecture (specified in Figure \ref{fig:sinNN})}
\nomenclature{$\text{LES}$}{Large Eddy Simulation}
\nomenclature{$\text{EFR}$}{Evolve--Filter--Relax}
\nomenclature{$\text{SST}$}{Shear Stress Transport}
\nomenclature{$\text{ODE}$}{Ordinary Differential Equation}

% Symbols
\nomenclature[G]{$\bm{u}$}{velocity field}
\nomenclature[G]{${p}$}{pressure field}
\nomenclature[G]{$\bm{\mu}$}{parameters considered for the ROM}
\nomenclature[G]{$\nu_t$}{eddy viscosity field}
\nomenclature[G]{$\bm{u}_r$}{reduced velocity field}
\nomenclature[G]{$p_r$}{reduced pressure field}
\nomenclature[G]{$\nu_{t_r}$}{reduced eddy viscosity field}

\nomenclature[G]{$\bm{a}$}{reduced vector of unknowns for velocity}
\nomenclature[G]{$\bm{b}$}{reduced vector of unknowns for pressure}
\nomenclature[G]{$\bm{g}$}{reduced vector of unknowns for eddy viscosity}
\nomenclature[G]{$\bm{a}^{proj}$}{projected reduced velocity coefficients (of dimension $N_u$)}
\nomenclature[G]{$\bm{b}^{proj}$}{projected reduced pressure coefficients (of dimension $N_p$)}

\nomenclature[G]{$\hat{\bm{a}}^{proj}$}{projected reduced velocity coefficients (of dimension $\hat{N}_u$)}
\nomenclature[G]{$\hat{\bm{b}}^{proj}$}{projected reduced pressure coefficients (of dimension $\hat{N}_p$)}

\nomenclature[G]{$Re$}{Reynolds Number}
\nomenclature[G]{$\mathcal{R}_{ij}$}{$(ij)_{th}$ component of the Reynolds stress tensor}
\nomenclature[G]{$\bm{E}_{ij}$}{$(ij)_{th}$ component of the averaged strain rate tensor}
\nomenclature[G]{${\nu}$}{kinematic viscosity}
\nomenclature[G]{${N_u^h}$}{number of unknowns for velocity at full-order level}
\nomenclature[G]{${N_p^h}$}{number of unknowns for pressure at full-order level}
\nomenclature[G]{$N_u$}{number of unknowns for velocity at reduced order level}
\nomenclature[G]{$N_p$}{number of unknowns for pressure at reduced order level}
\nomenclature[G]{$\hat{N}_u$}{number of velocity modes re-introduced in the ROM with correction terms}
\nomenclature[G]{$\hat{N}_p$}{number of pressure modes re-introduced in the ROM with correction terms}
\nomenclature[G]{$\bm{\tau}_u$}{correction term for the momentum equation}
\nomenclature[G]{$\bm{\tau}_p$}{correction term for the PPE}
%\nomenclature[G]{$\Toffline$}{training/offline set of time instances}
%\nomenclature[G]{$\Tonline$}{online set of time instances}
%\nomenclature[G]{$\nutrain$}{training/offline set of viscosity values}
%\nomenclature[G]{$\nutest$}{online set of viscosity values}
\nomenclature[G]{$N_{\text{networks}}$}{number of neural networks used to compute the confidence interval}
\nomenclature[G]{$\varepsilon_{\bm{u}}^{\nu^{\star}}(t)$}{relative error of velocity field with respect to the full order counterpart, for the parameter $\nu^{\star}$}
\nomenclature[G]{$\varepsilon_{p}^{\nu^{\star}}(t)$}{relative error of pressure field with respect to the full order counterpart, for the parameter $\nu^{\star}$}

\nomenclature[G]{$N_P$}{dimensionality of parameters}
%\nomenclature[G]{$N_T$}{number of time instances}
%\nomenclature[G]{$N_M$}{number of viscosity values}
\nomenclature[G]{$\mathbb{V}^u_{\text{POD}}$}{reduced basis space for velocity}
\nomenclature[G]{$\mathbb{V}^p_{\text{POD}}$}{reduced basis space for pressure}
\nomenclature[G]{$\mathbb{V}^{\nu_t}_{\text{POD}}$}{reduced basis space for pressure}
\nomenclature[G]{${\Omega}$}{bounded domain}
\nomenclature[G]{${\Gamma}$}{boundary of $\Omega$}
\nomenclature[G]{${\Gamma_{D_i}}$}{i-th boundary of $\Omega$, with non-homogeneous Dirichlet boundary condition}
\nomenclature[G]{$N_{\text{BC}}$}{number of boundaries of $\Omega$ with Dirichlet non-homogeneous boundary conditions}
\nomenclature[G]{$\bm{n}$}{outward normal vector}
\nomenclature[G]{$\bm{\varphi_i}$}{i-th POD basis function for velocity}
\nomenclature[G]{${\chi_i}$}{i-th POD basis function for pressure}
\nomenclature[G]{$\eta_{i}$}{i-th POD basis function for eddy viscosity}

\nomenclature[G]{$\mathcal{G}(\cdot)$}{mapping used to compute the reduced eddy viscosity field}
\nomenclature[G]{$\mathcal{M}(\cdot)$}{mapping used to compute the correction terms}

\nomenclature[G]{$\bm{{\mathcal{S}_u}}$}{snapshots matrix for the velocity field}
\nomenclature[G]{$\bm{{\mathcal{S}_p}}$}{snapshots matrix for the pressure field}
\nomenclature[G]{$\bm{{\mathcal{S}_{\nu_t}}}$}{snapshots matrix for the eddy viscosity field}
\nomenclature[G]{$\bm{f}(\cdot;\bm{\mu})$}{residual of the reduced momentum conservation equation}
\nomenclature[G]{$\bm{c}(\cdot; \bm{\mu})$}{residual of the reduced continuity equation}
\nomenclature[G]{$\bm{h}(\cdot; \bm{\mu})$}{residual of the reduced PPE}
\nomenclature[G]{$\bm{\tau}^{exact}$}{exact correction term}
\nomenclature[G]{$\bm{\tau}^{approx}$}{approximated correction term}
\nomenclature[G]{$\tau$}{weight used for boundary penalty term}
\nomenclature[G]{$N_{\text{seq}}$}{sequence length in the LSTM neural network}

\nomenclature[G]{$\otimes$}{tensor product}
\nomenclature[G]{$\bm{\nabla}\cdot$}{divergence operator}
\nomenclature[G]{$\bm{\nabla}\times$}{curl operator}
\nomenclature[G]{$\bm{\nabla}$}{gradient operator}
\nomenclature[G]{$\Delta$}{laplacian operator}
\nomenclature[G]{$\left\lVert \cdot\right\rVert_{L^2(\Omega)}$}{norm in $L^2(\Omega)$}
\nomenclature[G]{$( \cdot , \cdot )_{L^2(\Omega)}$}{inner product in $L^2(\Omega)$}

\begin{abstract}
\noindent In this paper, we propose an equation-based parametric Reduced Order Model (ROM), whose accuracy is improved with data-driven terms added into the reduced equations.
These additions have the aim of reintroducing contributions that in standard reduced-order approaches are not taken into account. In particular, in this work we focus on a Proper Orthogonal Decomposition (POD)-based formulation and our goal is to build a \emph{closure} or \emph{correction} model, aimed to re-introduce the contribution of the discarded modes.
The approach has been investigated in previous works such as~\cite{mohebujjaman2019physically, ivagnes2023pressure, ivagnes2023hybrid} and the goal of this manuscript is to extend the model to a \textbf{parametric} setting making use of machine learning procedures, and, in particular, of \emph{deep operator networks}.
More in detail, we model the closure terms through a deep operator network taking as input the reduced variables and the parameters of the problem.
We tested the methods on three test cases with different behaviors: the periodic turbulent flow past a circular cylinder, the unsteady turbulent flow in a channel-driven cavity, and the geometrically-parametrized backstep flow.
%In the first two test cases, the parameters are time and the Reynolds number, while the last test case is characterized by three geometrical parameters.
The performance of the machine learning-enhanced ROM is deeply studied in different modal regimes, and considerably improved the pressure and velocity accuracy with respect to the standard POD-Galerkin approach.
\end{abstract}

%----------------------------------------------------------------------------------------

%----------------------------------------------------------------------------------------
%	ARTICLE CONTENTS
%----------------------------------------------------------------------------------------

\section{Introduction}
\label{sec:intro}
Reduced Order Models (ROMs) \cite{degruyter1, degruyter2, degruyter3,rozza2008reduced, rozza2013reduced,aromabook} are a powerful tool used to reduce the computational effort of 
time-demanding simulations.
One of the fields where this class of techniques is widespread is Computational Fluid Dynamics,
where high-fidelity simulations may take days or weeks, even in the case of parallel computations
on many cores. For this reason, a simplified model is necessary to efficiently compute the solutions
for unseen configurations.

Most of the ROMs are built upon an \emph{offline-online} paradigm \cite{aromabook}.
The \emph{offline} stage consists of the computation of a large number of expensive high-fidelity simulations,
performed after setting the Full Order Model (FOM), which is usually a discretized version of complex PDEs, like the Navier--Stokes Equations (NSE).
The goal of this stage is the collection of the so-called \emph{snapshots}, namely the solutions of the simulations.
On the other hand, in the \emph{online} stage the full-order manifold is projected into a space with reduced dimensionality,
resulting in a reduced representation of the snapshots.
This reduction step may be assessed both with linear or nonlinear approaches. In particular, we employ the Proper Orthogonal Decomposition (POD) \cite{kosambi2016statistics, chatterjee2000introduction}, a linear technique which may be considered equivalent to the Principal Component Analysis (PCA) \cite{pearson1901liii, wold1987principal} and Singular Value Decomposition (SVD) \cite{lawson1995solving, golub1971singular, golub1965calculating}.

While the offline stage usually employs the numerical resolution of complex PDEs, in \emph{intrusive} ROMs the online stage allows for the resolution of ODEs, which consists of a reduced and simplified version of the FOM.

In particular, in this contribution, we focus on POD-Galerkin ROMs, that are based on a Galerkin projection \cite{noack1994low, bergmann2009enablers, kunisch2002galerkin, Azaiez2017}.
The main assumption of these models is that the solution may be approximated as a convex combination of a reduced number of global basis (the \emph{modes}),
whose coefficients are the reduced representations, namely the unknown variables in the ROM.
The above-mentioned linearity assumption, while considerably reducing the number of the system's unknowns, may make the model inaccurate, especially in the case of advection-dominated flows and/or where transport phenomena are dominant.

The goal of the ROM community is usually to retain a small number of modes, and to have efficient online simulations, namely a \emph{under-resolved} or \emph{marginally-resolved} regime. In this regime, the number of modes is enough to capture the dynamics of the system, but the POD-Galerkin ROM may lead to inaccurate results \cite{rowley2004model, iollo2000stability}.

Indeed, in advection-dominated cases, where the decay is slow, POD-Galerkin ROMs may require hundreds of modes to provide accurate results.

The inaccuracy may be caused by the formation of spurious oscillations in the approximated solution or by the ill-conditioning of the reduced system.
This motivates the use of data-driven techniques to stabilize and/or enhance the system solution. Some examples of works integrating ROMs with machine learning strategies are \cite{zancanaro2021hybrid, papapicco2022neural, pichi2023artificial, ivagnes2024enhancing, demo2023deeponet}.

Moreover, several researchers in the ROM community are recently putting their efforts into performing a system \emph{closure}, by adding to the system data-driven terms aimed to close the gap between the system solution and a reference optimal solution \cite{xie2020closure, ahmed2021closures, akhtar2012new}.
In this work, and in general in POD-based approaches, the data-driven \emph{correction} or closure terms may allow to reintegrate the contributions of the neglected modes into the system, as done in \cite{xie2018data, mohebujjaman2018physically, mohebujjaman2019physically}. 

In this framework, the key problem is to choose a proper and expressive closure model. Some proposed maps are either of \emph{quadratic} form, as in \cite{xie2018data, ivagnes2023pressure, barnett2022quadratic, geelen2023operator} or \emph{nonlinear} maps, for instance modeled through a neural network, as in \cite{barnett2023neural, san2018machine, san2018neural, xie2020closure}.

We propose here a nonlinear closure model through an efficient \emph{deep operator network} taking as input the reduced variables and the system's parameters.

In particular, the novelty of the present contribution consists in extending the previous work \cite{ivagnes2023hybrid} to a \textbf{parametrized} setting and to more challenging test cases.
The work \cite{ivagnes2023hybrid} proposed a quadratic closure model of the reduced variables, where the quadratic operator is constant in the parameters' space. We propose here to extend this model to a nonlinear learnable and parametric deep operator network.

We show numerical evidence of the method for three test cases in a turbulent regime:
\begin{itemize}
    \item[($\bm{a}$)] the periodic turbulent flow past a circular cylinder;
    \item[($\bm{b}$)] the unsteady channel-driven cavity flow;
    \item[($\bm{c}$)] the geometrically-parametrized flow over a backstep, considering three geometrical parameters.
\end{itemize}

The contribution is structured as follows.
\begin{itemize}
    \item The methodological framework in Section \ref{sec:methods}, including the FOM description, the foundations of POD-Galerkin ROMs, and the deep operator network proposed.
    \item The numerical results in Section \ref{sec:results}, focused on the above-mentioned turbulent test cases, with an additional discussion section \ref{subsec:discussion}.
    \item The conclusive part \ref{sec:conclusions}, which summarizes the key results obtained in the contribution and proposes future outlooks of the project.
\end{itemize}

\section{Numerical methods}
\label{sec:methods}

This Section is dedicated to the presentation of the numerical methods used in
this project. In particular, it is divided in the following parts:
\begin{itemize}
    \item the presentation of the Full Order Model (FOM) in Subsection
    \ref{subsec:fom}, which is used to compute the \emph{offline} solutions, named
    \emph{snapshots};
    \item a brief overview of the POD-Galerkin ROM approach for turbulent flows, in Subsection
    \ref{subsec:ev-roms}, which is the base model of our study;
    \item the closure approach used to enhance the standard
    POD-Galerkin formulation, in Subsection \ref{subsec:dd-ev-roms}. Here we briefly
    describe the procedure used to model the \emph{exact} closure/correction term.
    \item the deep operator networks' architecture used for the closure model and for the eddy viscosity modeling in Subsection \ref{subsec:networks}.
\end{itemize}

\subsection{Full Order Model (FOM)}
\label{subsec:fom}

We call the fluid domain $\Omega \in \mathbb{R}^d$ with $d=2$ (since we are
considering only 2-dimensional domains), $\Gamma$ its boundary. We indicate by
$\bm{\mu}$ the flow parameter vector. We consider $\bm{\mu} \in \mathbb{R}^P$, where $P$ is the dimensionality of the parameters. We consider time $t \in [0, T]$, $\bm{u}=\bm{u}(\bm{x}, \bm{\mu})$ is the flow velocity
vector field, $p=p(\bm{x}, \bm{\mu})$ is the pressure scalar field normalized by
the fluid density, and $\nu$ is the kinematic viscosity. In all the test cases we consider in the numerical results,
the FOM is based on the Reynolds--Averaged Navier--Stokes (RANS)
formulation for incompressible flows, which consists of a time-averaged version
of the Navier--Stokes equations.
More in detail, we consider a steady formulation for test case (c) (analyzed in Subsection \ref{subsec:test-case-c}), namely the geometrically-parametrized backstep, and an unsteady formulation (URANS) in test cases (a) and (b) (analyzed in Subsections \ref{subsec:test-case-a} and \ref{subsec:test-case-b}, respectively).

The main hypothesis that characterizes the RANS approach is the \emph{Reynolds
decomposition}~\cite{reynolds1895iv}. This theory is based on the assumption
that each flow field $\bm{s}$ can be expressed as the sum of its time
averaged value, indicated as $\overline{\bm{s}}$, and fluctuating parts,
indicated with $\bm{s}^{\prime}$. Such mean has different definitions depending
on the case of application. Classical choices are, for example, time averaging,
averaging along homogeneous directions or ensemble averaging.

We briefly recall here the standard URANS formulation for the incompressible
Navier-Stokes equations:
\begin{equation}
\begin{cases}
 \dfrac{\partial \overline{u}_{i}}{\partial x_i} = 0 ,\smallskip\\
    \dfrac{\partial \overline{u}_i}{\partial t} + \overline{u}_{j} \dfrac{\partial \overline{u}_{i} }{\partial x_j}=-\dfrac{\partial \overline{p}}{\partial x_i}+ \dfrac{\partial( 2 \nu\overline{\bm{E}}_{ij}- \mathcal{R}_{ij})}{\partial x_j},
    \label{RANS}
\end{cases}
\end{equation}

where the Einstein notation has been adopted, $\mathcal{R}_{ij}=\overline{u'_i
u'_j}$ is the Reynolds stress tensor, and
$\overline{\bm{E}}_{ij}=\dfrac{1}{2}\left(\dfrac{\partial
\overline{u}_{i}}{\partial x_j} + \dfrac{\partial \overline{u}_{j}}{\partial
x_i}\right)$ is the averaged strain rate tensor.

The URANS formulation in \eqref{RANS} needs to be coupled with a turbulence
model to close system \eqref{RANS}. In particular, we adopt the $\kappa-\omega$
\emph{Shear Stress Transport} (SST) model~\cite{menter1994two}.

This model belongs to the class of \emph{eddy viscosity models}, which are based
on the Boussinesq hypothesis, i.e. the turbulent stresses are related to the
mean velocity gradients as follows:
\[
 -\mathcal{R}_{ij}=2 \nu_t \overline{\bm{E}}_{ij} - \dfrac{2}{3} \kappa \delta_{ij},
    \label{bouss}
\]
where $\kappa=\frac{1}{2}\overline{u'_i u'_i}$ is the turbulent kinetic energy
and $\nu_t$ is the eddy viscosity. In general, in this case the RANS (or URANS) model is enriched with two additional transport equations, for $\kappa$ and $\omega$, respectively. For the complete model we refer the reader to
the original paper \cite{menter1994two}, and the extended RANS model including
the SST $\kappa-\omega$ equations can be found in \cite{hijazi2020data}.

The full order solutions of \ref{RANS} are computed by means of the open-source
software \emph{OpenFOAM}, which employs a finite volume discretization of the
RANS equations \cite{moukalled2016finite, jasak1996error}.

Following the finite volume method, developed and implemented in
\cite{moukalled2016finite, jasak1996error}, the computational domain is
discretized in polygonal control volumes and the partial differential equations
\ref{RANS} are integrated in each control volume and converted to algebraic
equations. In particular, the volume integrals are converted into surface
integrals by the divergence theorem, and discretized as sums of the fluxes at
the boundary faces of the control volumes~\cite{moukalled2016finite}. 

Typically, the velocity-pressure coupling is based on the so-called Pressure Poisson Equation approach (PPE), which consists in replacing the continuity equation with a dedicated pressure equation, which is obtained by taking the divergence of the momentum equation.
In such cases, the URANS system \eqref{RANS} at the full order level
can be written as follows:
\begin{equation}
    \begin{cases}
    \dfrac{\partial \overline{\bm{u}}}{\partial t}+ \nabla \cdot (\overline{\bm{u}} \otimes \overline{\bm{u}})=\nabla \cdot \left[-\overline{p} \bm{I} +(\nu+\nu_t) \left(\nabla \overline{\bm{u}} + (\nabla \overline{\bm{u}})^T \right) \right] & \text{ in } \Omega \times [0,T]\, ,\\
    \Delta \overline{p}=-\nabla \cdot (\nabla \cdot (\overline{\bm{u}} \otimes \overline{\bm{u}})) +\nabla \cdot \left[ \nabla \cdot \left( \nu_t \left(\nabla \overline{\bm{u}} +(\nabla \overline{\bm{u}})^T \right) \right) \right] & \text{ in }\Omega \, ,\\
    + \text{ Boundary Conditions }& \text{ on } \Gamma \times [0,T]\, ,\\
    + \text{ Initial Conditions } & \text{ in } (\Omega,0)\, .
    \end{cases}
    \label{RANS-PPE}
\end{equation}

\subsection{POD-Galerkin and Eddy-Viscosity ROMs}
\label{subsec:ev-roms}

In this section, we provide a brief overview of the standard POD-Galerkin ROM approach, and we also briefly recall the turbulence modeling at the reduced level. For more details we refer the reader to \cite{hijazi2020data}.
\medskip

We consider $N_{\mu}$ the number of training parameters considered to run the offline stage.
Once all the high-fidelity simulations are run, all the FOM snapshots, i.e., the FOM solutions for different parameter values  $\{\bm{\mu}_j\}_{j=1}^{N_{\mu}}$, can be collected. 
%in our case for different time instances $\{t_i\}_{i=1}^{N_T}$ and for different viscosity values $\{\nu_k\}_{k=1}^{N_M}$, where $N_{\mu}=N_T\times N_M$. In particular, each parameter is $[t_i, \nu_k], \, i=1, \dots, N_T, \, k=1, \dots, N_M$.
The POD is then applied on the full order velocity and pressure snapshots' matrices
\begin{equation*}
\mathcal{S}_u=\{\bm{u}(\bm{x},\bm{\mu}_1),...,\bm{u}(\bm{x},\bm{\mu}_{N_{\mu}})\}  \in \mathbb{R}^{N_u^h \times N_{\mu}}, \quad \mathcal{S}_p=\{p(\bm{x},\bm{\mu}_1),...,p(\bm{x},\bm{\mu}_{N_{\mu}})\}  \in \mathbb{R}^{N_p^h \times N_{\mu}},
\end{equation*}
and (in case of turbulent flows) on the eddy viscosity snapshots' matrix
\begin{equation*}
\mathcal{S}_{\nu_t}=\{\nu_t(\bm{x},\bm{\mu}_1),...,\nu_t(\bm{x},\bm{\mu}_{N_{\mu}})\}  \in \mathbb{R}^{N_{\nu_t}^h \times N_{\mu}},
\end{equation*}
If we call $N_c$ the number of the mesh cells considered, $N_u^h=d\, N_c$ and $N_p^h=N_{\nu_t}^h=N_c$ are the numbers of spatial degrees of freedom for the
velocity, pressure and eddy viscosity fields, respectively.

After the POD is applied to the snapshots' matrices, the following reduced POD
spaces are found:
\begin{equation}
\mathbb{V}^u_{\text{POD}}=\mbox{span}\{[\boldsymbol{\phi}_i]_{i=1}^{N_u}
\},\quad 
\mathbb{V}^p_{\text{POD}}=\mbox{span}\{[\chi_i]_{i=1}^{N_p}\}, \quad 
\mathbb{V}^{\nu_t}_{\text{POD}}=\mbox{span}\{[\eta_i]_{i=1}^{N_{\nu_t}}\},
\end{equation}
where $N_u \ll N_u^h$, $N_p \ll N_p^h$, $N_{\nu_t} \ll N_{\nu_t}^h$.
$[\boldsymbol{\phi}_i]_{i=1}^{N_u}$, $[\chi_i]_{i=1}^{N_p}$, and $[\eta_i]_{i=1}^{N_{\nu_t}}$ indicate the
velocity, pressure and eddy viscosity modes, respectively.

The combination of modes $(N_u, N_p, N_{\nu_t})$ has to be defined \emph{a priori}, based on the retained energy and on the knowledge of the system of interest.

The main hypothesis of the POD is that each field can be then approximated as a
convex combination of its modes, namely:
\begin{equation}
\begin{split}
\bm{u}(\bm{x},\bm{\mu}) \approx \bm{u}_r(\bm{x}, \bm{\mu})=&\sum_{i=1}^{N_u} a_i(\bm{\mu})\boldsymbol{\phi}_i(\bm{x}), \quad
p(\bm{x},\bm{\mu})\approx p_r(\bm{x},\bm{\mu})=\sum_{i=1}^{N_p} b_i(\bm{\mu})\chi_i(\bm{x}), \\
&\nu_t(\bm{x},\bm{\mu})\approx {\nu_t}_r(\bm{x},\bm{\mu})=\sum_{i=1}^{N_{\nu_t}} g_i(\bm{\mu})\eta_i(\bm{x}).
\end{split}
\label{eq:appfield}
\end{equation}

Supposing that the reduced fields provide an accurate approximation of the
original snapshots, the reduced order system can be written in a compact form
as:
\begin{equation}
    \begin{cases}
        \bm{f}(\bm{a}, \bm{b}, \bm{g}; \bm{\mu})=\bm{0},\\
        \bm{c}(\bm{a}, \bm{g}; \bm{\mu}) = \bm{0},
    \end{cases}
    \label{eq:compact-stand-rom}
\end{equation}
where $\bm{a}=(a_i)_{i=1}^{N_u}$,
$\bm{b}=(b_i)_{i=1}^{N_p}$, $\bm{g}=(g_i)_{i=1}^{N_{\nu_t}}$ are the vectors of coefficients for the velocity, pressure and eddy viscosity fields. $\bm{f}(\cdot)$ and $\bm{c}(\cdot)$ are the reduced momentum and continuity residual vectors, respectively. 

In \eqref{eq:compact-stand-rom} and in the following parts of the manuscript we
will omit the dependencies $\bm{a}(\bm{\mu})$ and $\bm{b}(\bm{\mu})$ for the
sake of brevity.

It may happen that the model \eqref{eq:compact-stand-rom} is not stable
in terms of velocity-pressure coupling and this may lead to spurious
oscillations, causing the residuals to increase as time evolves. For this
reason, there exist in literature different stabilization approaches for the
velocity-pressure coupling. Two of those approaches are the supremizer method
(SUP-ROM) and the PPE approach (PPE-ROM). The first
approach consists in enriching the velocity reduced space with additional modes
named \emph{supremizer modes}, which are directly computed either from the
pressure modes (\emph{exact} supremizer approach) or from the pressure snapshots
(\emph{approximated} supremizer approach) \cite{rozza2007stability,
ballarin2015supremizer, stabile2018finite, ali2020stabilized}. The second method
provides a reduced version of system \ref{RANS-PPE}, replacing the continuity equation with the PPE equation, also at the reduced level
\cite{akhtar2009stability,Stabile2017CAIM,stabile2018finite, noack2005need}.

In this project, we choose to adopt the PPE-ROM stabilization for two specific
reasons:
\begin{itemize}
    \item[(i)] it relies on a dedicated pressure equation and it allows for the
    introduction of specific \emph{pressure correction} terms only acting on the
    pressure accuracy, as will be specified in Section \ref{subsec:dd-ev-roms};
    \item[(ii)] it is consistent with velocity-pressure coupling at the FOM level, namely with the implementation used in OpenFOAM.
\end{itemize}

In the case of PPE-ROMs we report here the dynamical system corresponding to a reduced version of system \eqref{RANS-PPE}.
From now on, we will refer to this model as \emph{Eddy-Viscosity ROM}, namely \emph{EV-ROM}. This will be our baseline for the closure modeling.

\begin{equation}
    \begin{cases}
    - \bm{M} \dot{\bm{a}}+\nu(\bm{B}+\bm{B_T})\bm{a}-\bm{a}^T \bm{C} \bm{a}+ \bm{g}^T (\bm{C}_{\text{T1}} +\bm{C}_{\text{T2}}) \bm{a}-\bm{H}\bm{b}+\tau \left( \sum_{k=1}^{N_{\text{BC}}}(U_{\text{BC},k}\bm{D}^k-\bm{E}^k \bm{a})\right)=\bm{0} \, ,\\
    \bm{D}\bm{b}+ \bm{a}^T \bm{G} \bm{a} -\bm{g}^T(\bm{C}_{\text{T3
    }} +\bm{C}_{\text{T4
    }})\bm{a} - \nu \bm{N} \bm{a}- \bm{L}=\bm{0}\, ,
    \end{cases}
    \label{eq:ppe-rom-turb}
\end{equation}

In the above formulation, the POD operators read as follows:
\begin{equation}
\begin{split}
&(\bm{M})_{ij}=(\boldsymbol{\phi}_i,\boldsymbol{\phi}_j)_{L^2(\Omega)}, \quad (\bm{P})_{ij}=(\chi_i,\nabla \cdot \boldsymbol{\phi}_j)_{L^2(\Omega)}\, ,\quad (\bm{B})_{ij}=(\boldsymbol{\phi}_i,\nabla \cdot \nabla \boldsymbol{\phi}_j)_{L^2(\Omega)}, \\
&(\bm{B_T})_{ij}=(\boldsymbol{\phi}_i,\nabla \cdot (\nabla \boldsymbol{\phi}_j)^T)_{L^2(\Omega)},\quad (\bm{C})_{ijk}=(\boldsymbol{\phi}_i,\nabla \cdot (\boldsymbol{\phi}_j \otimes \boldsymbol{\phi}_k))_{L^2(\Omega)}, \quad (\bm{H})_{ij}=(\boldsymbol{\phi}_i,\nabla \chi_j)_{L^2(\Omega)}\, \\
&(\bm{D})_{ij}=(\nabla \chi_i,\nabla \chi_j)_{L^2(\Omega)}, \quad 
(\bm{G})_{ijk}=(\nabla \chi_i,\nabla \cdot (\boldsymbol{\phi}_j \otimes \boldsymbol{\phi}_k))_{L^2(\Omega)}, \\ &(\bm{N})_{ij}=(\bm{n} \times \nabla \chi_i,\nabla \boldsymbol{\phi}_j)_\blue{{L^2(\Gamma)}}, \quad (\bm{L})_{ij}=(\chi_i,\bm{n} \cdot \boldsymbol{R}_t)_\blue{{L^2(\Gamma)}}\,,
\end{split}
\label{eq:operators-rom}
\end{equation}
vector $\bm{n}$ is the normal unitary vector to the domain boundary.

The operators including the turbulence contribution are:
\[
\begin{split}
    &(\bm{C}_{\text{T1}})_{ijk}=(\boldsymbol{\phi}_i, \eta_j \nabla \cdot \nabla \boldsymbol{\phi}_k)_{L^2(\Omega)} \, ,\quad
    (\bm{C}_{\text{T2}})_{ijk}=(\boldsymbol{\phi}_i, \nabla \cdot \eta_j (\nabla \boldsymbol{\phi}_k)^T)_{L^2(\Omega)}\, ,\\
    &(\bm{C}_{\text{T3}})_{ijk}=(\nabla \chi_i, \eta_j \nabla \cdot \nabla \boldsymbol{\phi}_k)_{L^2(\Omega)}\, , \quad (\bm{C}_{\text{T4}})_{ijk}=(\nabla \chi_i, \nabla \cdot \eta_j(\nabla \boldsymbol{\phi}_k)^T)_{L^2(\Omega)}\,.
\end{split}
\]

Moreover, the additional term in the momentum equation is a \emph{weak}
enforcement of the non-homogeneous Dirichlet boundary conditions, according to
the \emph{penalty method}. The penalization factor, namely $\tau$, is tuned through
a sensitivity analysis on the specific problem considered \cite{hijazi2020data,
star2019extension}. In general, bigger values of the penalization factor lead to
a stronger enforcement of the boundary conditions. At the same time, large values of $\tau$ also lead to larger condition number for system \eqref{eq:ppe-rom-turb} and, hence, to less stability.

The matrices $\bm{E}^k$ and
vectors $\bm{D}^k$ are defined as follows:
\[(\bm{E}^k)_{ij}=(\boldsymbol{\phi}_i,
\boldsymbol{\phi}_j)_{L^2(\Gamma_{D_k})}, \quad
(\bm{D}^k)_{i}=(\boldsymbol{\phi}_i)_{L^2(\Gamma_{D_k})}, \text{ for all
}k=1,...,N_{\text{BC}},\]

where $\Gamma_{D_i}$ is the i-th boundary with non-homogeneous Dirichlet conditions, $N_{\text{BC}}$ is the total number of boundaries where non-homogeneous Dirichlet boundary conditions are imposed.

\subsubsection{On the eddy viscosity modeling}
The compact form of
system \eqref{eq:ppe-rom-turb} is:
\begin{equation}
\begin{cases}
    \bm{f}(\bm{a}, \bm{b}, \bm{g};\bm{\mu}) = \bm{0}\, ,\\
    \bm{h}(\bm{a}, \bm{b}, \bm{g};\bm{\mu}) = \bm{0}\, .
\end{cases}
    \label{eq:compact-ppe-rom-turb}
\end{equation}

In this dynamical system, the number of
unknowns is $N_u$ for the velocity, $N_p$ for the pressure, and $N_{\nu_t}$ for
the eddy viscosity. However, the number of equations is $N_u+N_p$. Thus, there
are more unknowns than equations and the system is not closed. In order to close
the systems, the eddy viscosity coefficients $[g_i(\bm{\mu})]_{i=1}^{N_{\nu_t}}$
can be computed considering the mapping $\bm{g}=\mathcal{G}(\bm{a}, \bm{\mu})$
through either \emph{interpolation}~\cite{lazzaro2002radial,
micchelli1986interpolation} or \emph{regression} techniques. An interpolation
technique was exploited in \cite{hijazi2020data}, following the POD-I approach
\cite{wang2012comparative,walton2013reduced,salmoiraghi2018free}.

In this paper, the mapping is a neural network whose loss function is the Mean Squared Error (MSE) between the output of the neural network
$\mathcal{G}(\bm{a}, \bm{\mu})$ and the known eddy viscosity coefficients
$\bm{g}^{proj}(\bm{\mu})$ found from the POD procedure.
The architecture of neural network $\mathcal{G}$ is a deep operator network and it is introduced in Subsection \ref{subsec:networks}.

\textbf{Remark 1: } It is worth specifying that in case of steady flows, the term $\bm{M} \dot{\bm{a}}$ in \eqref{eq:ppe-rom-turb} vanishes.
Moreover, the system \eqref{eq:ppe-rom-turb} is solved (at each time step, in case of unsteady flows) employing the Newton method.

\subsection{The data-driven closure approach for EV-ROMs: DD-EV-ROM}
\label{subsec:dd-ev-roms}
This part of the manuscript is dedicated to the description of the closure approach adopted to enhance the ROM results.

The main goal here is to identify a strategy to model the
\emph{correction/closure} terms already exploited in previous works like
\cite{mohebujjaman2019physically, mou2021data, san2013proper, ahmed2021closures,
ivagnes2023pressure, dar2023artificial, ivagnes2023hybrid}, but in a parametrized
setting.

The procedure used is as follows.
\begin{enumerate}
    
    \item[\textbf{1}.] Select a reduced dimension for the velocity, $N_u$, and for the pressure, $N_p$. The sum $N_u+N_p$ is the dimension
    of the reduced system \eqref{eq:ppe-rom-turb}.
    \item[\textbf{2}.] Select two bigger dimensions $\hat{N}_u>N_u$ and $\hat{N}_p>N_u$, for
    instance from literature we have $\hat{N}_u=\mathrm{k}{N}_u$, $\hat{N}_p=\mathrm{k}N_p$, with
    $\mathrm{k} \in \mathbb{N}$.
    \item[\textbf{3}.] Select one or more operators, that we name $\mathcal{C}$.
        For instance, if we consider the non-linear operators in System
        \eqref{eq:ppe-rom-turb} we obtain: \begin{equation} \mathcal{C}(\bm{a})=
        \begin{bmatrix}
            -\bm{a}^T \bm{C} \bm{a}\\
            \bm{a}^T \bm{G} \bm{a} \end{bmatrix}.
            \label{eq:nonturb-operator}
    \end{equation}
    
    \item[\textbf{4}.] Compute the term $\mathcal{C}(\hat{\bm{a}}^{proj}) \in \mathbb{R}^{\hat{N}_u+\hat{N}_p}$:
    \begin{equation}
    \mathcal{C}(\hat{\bm{a}}^{proj})=
    \begin{bmatrix}
            -(\hat{\bm{a}}^{proj})^T \hat{\bm{C}} \hat{\bm{a}}^{proj}\\
            (\hat{\bm{a}}^{proj})^T \hat{\bm{G}} \hat{\bm{a}}^{proj}
        \end{bmatrix},
        \label{eq:operator-big}
    \end{equation}
    where $\hat{\bm{a}}^{proj}$ is the velocity coefficients' vector found by
    directly projecting the field on the POD subspace of dimension $\hat{N}_u$. The
    operators $\hat{\bm{C}}$ and $\hat{\bm{G}}$ are expressed as:
    \[
    \begin{split}
    &(\hat{\bm{C}})_{ijk}=(\boldsymbol{\phi}_i,\nabla \cdot (\boldsymbol{\phi}_j \otimes \boldsymbol{\phi}_k))_{L^2(\Omega)}, \quad i, j, k=1, \dots, \hat{N}_u, \\
    &(\hat{\bm{G}})_{ijk}=(\nabla \chi_i,\nabla \cdot (\boldsymbol{\phi}_j \otimes \boldsymbol{\phi}_k))_{L^2(\Omega)}, \quad i=1, \dots, h; j, k=1, \dots, \hat{N}_p.
    \end{split}
    \]
    \item[\textbf{5}.] Compute the term $\mathcal{C}(\bm{a}^{proj}) \in \mathbb{R}^{N_u+N_p}$:
    \begin{equation}
    \mathcal{C}(\bm{a}^{proj})=
    \begin{bmatrix}
            -(\bm{a}^{proj})^T \bm{C} \bm{a}^{proj}\\
            (\bm{a}^{proj})^T \bm{G} \bm{a}^{proj}
        \end{bmatrix},
        \label{eq:operator-reduced}
    \end{equation}
    where $\bm{a}^{proj}$ is the velocity coefficients' vector found by
    directly projecting the field on the POD subspace of dimension $N_u$. The
    operators $\bm{C}$ and $\bm{G}$ are the reduced operators appearing in
    \eqref{eq:ppe-rom-turb} and already defined in \eqref{eq:operators-rom}.
    \item[\textbf{6}.] Compute the \emph{exact} correction term as follows:
    \begin{equation}
        \boldsymbol{\tau}^{exact} = \overline{\mathcal{C}(\hat{\bm{a}}^{proj})}^{N_u+N_p} - \mathcal{C}(\bm{a}^{proj}),
        \label{eq:exact-correction}
    \end{equation}
    where $\overline{(\cdot)}^r$ acts like a filter and indicates that only the
    first $r$ components should be retained. In our example:
     \begin{equation}
    \boldsymbol{\tau}^{exact}=
    \begin{bmatrix}
            -\overline{(\hat{\bm{a}}^{proj})^T \hat{\bm{C}} \hat{\bm{a}}^{proj}}^{N_u} + (\bm{a}^{proj})^T \bm{C} \bm{a}^{proj}\\
            \overline{(\hat{\bm{a}}^{proj})^T \bm{G}_p \hat{\bm{a}}^{proj}}^{N_p} - (\bm{a}^{proj})^T \bm{G} \bm{a}^{proj}
        \end{bmatrix},
        \label{eq:exact-nonlinear-correction}
    \end{equation}
    where we retain the first $N_u$ and $N_p$ components for the first and second
    term, respectively.
    \item[\textbf{7}.] Finally, one has to choose a \emph{model} $\mathcal{M}$ to approximate with good accuracy the exact correction term. In particular, $\boldsymbol{\tau}^{approx}=\boldsymbol{\tau}^{approx}(\bm{a},
   \bm{\mu})=\mathcal{M}(\bm{a}, \bm{\mu})$, where $\bm{\mu}$
    includes the parameters of the test case taken into account.
    Some of the
    possible choices to model the mapping $\mathcal{M}$ are the following.
    \begin{itemize}
        \item Choose an \emph{ansatz}, for example a quadratic ansatz with
        respect to the coefficients' vectors:
        \begin{equation}
            \mathcal{M}(\bm{a})= 
                \tilde{\bm{A}} \bm{a}+ \bm{a}^T\tilde{\bm{B}}  \bm{a}.
                \label{eq:quadratic-ansatz}
            \end{equation}
            where the matrix $\tilde{\bm{A}} \in \mathbb{R}^{(N_u+N_p)\times(N_u+N_p)}$
            and the operator $\tilde{\bm{B}} \in \mathbb{R}^{(N_u+N_p)\times(N_u+N_p)\times(N_u+N_p)}$. This approach is easily addressed in
            time-dependent problems by solving a minimization problem:
            \begin{equation}
                \tilde{\bm{A}}, \tilde{\bm{B}}=\text{arg} \min_{\bm{A}, \bm{B}}\left(\sum_{j=1}^{M} \left\lVert \Bigl(\bm{A} \,\bm{a}^{proj}(t_j)+ \bm{a}^{proj}(t_j)^T {\bm{B}} \,  \bm{a}^{proj}(t_j)  \Bigr)- \bm{\tau}^{exact}(t_j) \right\rVert^2 \right).
                \label{eq:min-least-squares}
            \end{equation}
            In previous works, like \cite{mou2021data,
            mohebujjaman2019physically, ivagnes2023hybrid, ivagnes2023pressure},
            the minimization problem \eqref{eq:min-least-squares} is typically
            re-written as a least squares problem. However, this can be done if
            time is the only parameter of the problem. Indeed, matrices
            $\tilde{\bm{A}}$ and $\tilde{\bm{B}}$ are parameter-constant and may not be appropriately generalized when multiple parameters are considered.
            Parametrized test cases require indeed more advanced mappings depending also on the remaining parameter(s).
            
        \item Train a \emph{neural network} (NN), which takes as input the
        coefficients and the parameters of the problem $(\bm{a},
        \bm{\mu})$ and gives as output the approximated correction coefficients
        $\bm{\tau}^{approx}$. In this case, the mapping $\mathcal{M}$ can be
        modeled considering different architectures, from multi-layer perceptron to recurrent NNs.
        A preliminary investigation comparing different models for test cases (a) and (b) can be found in \cite{ivagnes2024parametric}. 
        In the numerical results' Section \ref{sec:results}, we employ on a \emph{Multi-Input deep Operator NETwork} (MIONet), whose structure is detailed in Subsection \ref{subsec:networks}.
        This second approach is suitable for multi-parametrized problems, whose dependency is easily addressed by adding the
        parameters as input to the neural network.
    \end{itemize}
\end{enumerate}

When we employ a reduced eddy-viscosity model, as in
\eqref{eq:compact-ppe-rom-turb}, the operator $\mathcal{C}$ can also include the
turbulence contributions, namely terms $\bm{g}^T (\bm{C}_{\text{T1}} +
\bm{C}_{\text{T2}})\bm{a}$ and $\bm{g}^T (\bm{C}_{\text{T3}} +
\bm{C}_{\text{T4}})\bm{a}$. In that case, it is specialized as:
\begin{equation} \mathcal{C}(\bm{a}, \bm{g})=
    \begin{bmatrix}
            -\bm{a}^T \bm{C} \bm{a}+\bm{g}^T(\bm{C_{T1}}+\bm{C_{T2}})\bm{a}\\
            \bm{a}^T \bm{G} \bm{a} -\bm{g}^T(\bm{C_{T3}}+\bm{C_{T4}})\bm{a}\end{bmatrix}.
        \label{eq:turb-operator}
\end{equation}
In this specific case, the model will be: $\boldsymbol{\tau}^{approx}=\boldsymbol{\tau}(\bm{a}, \bm{g}, \bm{\mu})=\mathcal{M}(\bm{a}, \bm{g}, \bm{\mu})$ 

Moreover, as in the previous works \cite{ivagnes2023pressure,
ivagnes2023hybrid}, we identify the first $N_u$ elements of the correction term
$\bm{\tau}=\mathcal{M}(\bm{a}, \bm{b}, \bm{\mu})$ as the \emph{velocity}
correction, whereas the last $N_p$ components are identifies as \emph{pressure}
correction. Hence, the $\bm{\tau}$ vector can be decomposed as $(\bm{\tau}_u,
\bm{\tau}_p)$, where $\bm{\tau}_u \in \mathbb{R}^{N_u}$ and $\bm{\tau}_p \in \mathbb{R}^{N_p}$.

We are finally able to write the so-called \emph{data-driven eddy viscosity ROM} (\emph{DD-EV-ROM}), which
combines the physics-based and the purely data-driven based approaches. Its
final compact formulation reads as follows:
\begin{equation}
    \begin{cases}
    \bm{f}(\bm{a}, \bm{b}, \bm{g};\bm{\mu})+\bm{\tau}_u(\bm{a}, \bm{g},\bm{\mu}) = \bm{0}\, ,\\
    \bm{h}(\bm{a}, \bm{b}, \bm{g};\bm{\mu})+\bm{\tau}_p(\bm{a}, \bm{g}, \bm{\mu})= \bm{0}\, .
\end{cases}
\label{eq:hyb-dd-rom}
\end{equation}

To clarify the overall procedure, we provide the pseudo-code for the online stage of DD-EV-ROM in Algorithm \ref{alg-steady} (for the steady case), and in \ref{alg-unsteady} (for the unsteady case).
In Algorithm \ref{alg-steady}, we also provide an additional step for the geometrical parameterization, considered for test case (\textbf{c}). In this case, all the reduced operators are parameterized and have to be computed before solving the reduced system.

A more detailed description can be found in \ref{app:geometrical-params}.

\begin{algorithm*}[htpb!]
\caption{Pseudo-code for the online stage of algorithms DD-EV-ROM in the steady case. The part highlighted in gray corresponds to the case with geometrical parametrization in Subsection \ref{subsec:test-case-c}.}
\label{alg-steady}
\begin{algorithmic}[1]
\State $\bm{\mu}^{\star}$ test parameters;  \\
\highlight{compute operators $\bm{B}(\bm{\mu}^{\star}), \bm{B_T}(\bm{\mu}^{\star}), \bm{C}(\bm{\mu}^{\star}), \bm{H}(\bm{\mu}^{\star}), \bm{D}(\bm{\mu}^{\star}), \bm{G}(\bm{\mu}^{\star}), \bm{N}(\bm{\mu}^{\star}), \bm{C_{T1}}(\bm{\mu}^{\star}), \bm{C_{T2}}(\bm{\mu}^{\star}),$\\
$\qquad \qquad \qquad \qquad \bm{C_{T3}}(\bm{\mu}^{\star}), \bm{C_{T4}}(\bm{\mu}^{\star}), \bm{D}^k(\bm{\mu}^{\star}), \bm{E}^k(\bm{\mu}^{\star})$, with $k=1, \dots, N_{\text{BC}}$.}\\
initial guess $\bm{a}_0=\bm{a}_0^{proj}=[a_i^{proj}(\bm{\mu}^{\star})]_{i=1}^{N_u}$, $\bm{b}_0=\bm{a}_0^{proj}=[b_i^{proj}(\bm{\mu}^{\star})]_{i=1}^{N_p}$;
\State $i=0$;
\Repeat 
    \State $\bm{g}_i \gets \mathcal{G}(\bm{a}_i, \bm{\mu}^*)$;
    \State $\bm{\tau}_i=[\bm{\tau}_{u}, \bm{\tau}_p]_i \gets \mathcal{M}(\bm{a}_i, \bm{g}_i, \bm{\mu}^{\star})$;
    \State residual $\bm{r}_i \gets 
    \begin{bmatrix}
        \bm{f}(\bm{a}_i, \bm{b}_i, \bm{g}_i;\bm{\mu}^{\star})+\bm{\tau}_u\,\\
    \bm{h}(\bm{a}_i, \bm{b}_i, \bm{g}_i;\bm{\mu^{\star}})+\bm{\tau}_p
    \end{bmatrix}$;
    \State compute $\bm{a}_{i+1}$, $\bm{b}_{i+1}$ through the Newton Method;
    \State $i \gets i+1$;
    \Until{Convergence or maximum iterations reached}
    \State $\bm{a}^{sol}\gets \bm{a}_{i}$, $\bm{b}^{sol}\gets \bm{b}_{i}$.
\end{algorithmic}
\end{algorithm*}

\begin{algorithm*}[htpb!]
\caption{Pseudo-code for the online stage of algorithms DD-EV-ROM in the unsteady case.}
\label{alg-unsteady}
\begin{algorithmic}[1]
\State $\bm{\mu}^{\star}_{phys}$ test physical parameters, except for time;  $\bm{\mu}_0^{\star}= [t_0, \bm{\mu}^{\star}_{phys}]$;\\
initial guess $\bm{a}^{sol}(t_0)=[a_i^{proj}(\bm{\mu}_0^{\star})]_{i=1}^{N_u}$, $\bm{b}^{sol}(t_0)=[b_i^{proj}(\bm{\mu}_0^{\star})]_{i=1}^{N_p}$;\\
time steps $t_n$, $n=1, \dots, N_T$;

\For{$n \in [0, \dots, N_T - 1]$}
\State $\bm{\mu}_n^{\star}\gets [t_n, \bm{\mu}^{\star}_{phys}]$;
\State $\bm{a}_0\gets \bm{a}^{sol}(t_{n+1})$, $\bm{b}_0\gets\bm{b}^{sol}(t_{n+1})$;
\State $i=0$;
\Repeat 
    \State $\bm{g}_i \gets \mathcal{G}(\bm{a}_i, \bm{\mu}_n^*)$;
    \State $\bm{\tau}_i=[\bm{\tau}_{u}, \bm{\tau}_p]_i \gets \mathcal{M}(\bm{a}_i, \bm{g}_i, \bm{\mu}_n^{\star})$;
    \State residual $\bm{r}_i \gets 
    \begin{bmatrix}
        \bm{f}(\bm{a}_i, \bm{b}_i, \bm{g}_i;\bm{\mu}_n^{\star})+\bm{\tau}_u\,\\
    \bm{h}(\bm{a}_i, \bm{b}_i, \bm{g}_i;\bm{\mu_n^{\star}})+\bm{\tau}_p
    \end{bmatrix}$;
    \State compute $\bm{a}_{i+1}$, $\bm{b}_{i+1}$ through the Newton Method;
    \State $i \gets i+1$;
    \Until{convergence or maximum iterations reached}
    \State $\bm{a}^{sol}(t_{n+1}) \gets \bm{a}_{i}$, $\bm{b}^{sol}(t_{n+1}) \gets \bm{b}_{i}$.
\EndFor 
\end{algorithmic}
\end{algorithm*}

\subsection{Deep Operator networks}
\label{subsec:networks}

This Section is dedicated to describe the architecture of the machine learning models used for the turbulence map $\mathcal{G}(\bm{a}, \bm{\mu})$ and for the closure map $\mathcal{M}(\bm{a}, \bm{g}, \bm{\mu})$.

\begin{figure}[htpb!]
    \centering
    \includegraphics[width=0.6\linewidth]{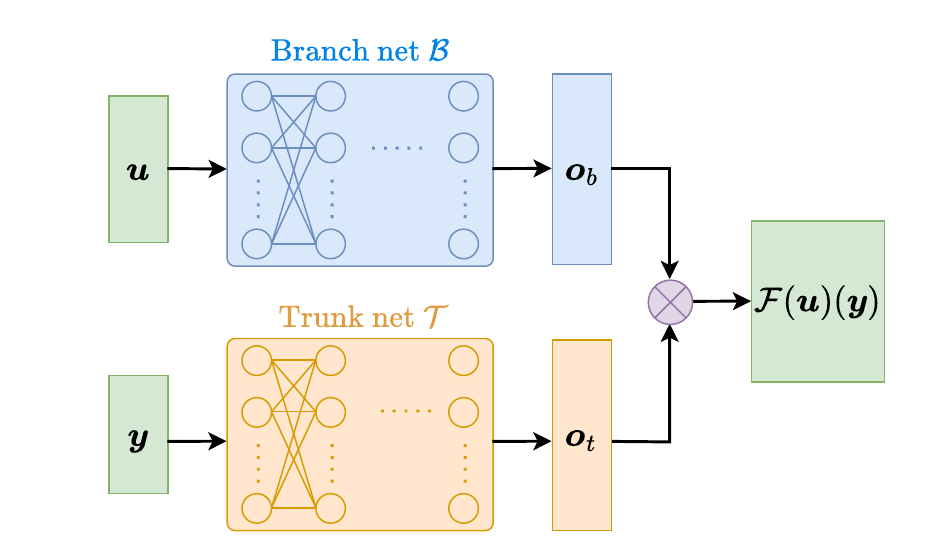}
    \caption{Standard version of the DeepONet, as in \cite{lu2019deeponet}.}
    \label{fig:deeponet-simple}
\end{figure}

For both mappings we employ \emph{deep operator networks}, a recently-proposed architecture specifically designed to learn operators.
The first example of deep operator network is the DeepONet \cite{lu2019deeponet} displayed in Figure \ref{fig:deeponet-simple}. Given the observations $\bm{u}(x_1), \bm{u}(x_2), \dots$, and the parameters $\bm{y}$, the goal is to model the operator $\mathcal{F}(\bm{u}(\bm{y}))$. The peculiarity of such model and the difference with respect to fully-connected neural networks is the presence of multiple sub-networks that separately process multiple inputs. In the notation we adopt, the \emph{branch network} $\mathcal{B}$ handles input $\bm{u}$, while the \emph{trunk network} $\mathcal{T}$ takes as input $\bm{y}$. The outputs of the two networks, namely $\bm{o}_b=\mathcal{B}(\bm{u})$ and $\bm{o}_t=\mathcal{T}(\bm{y})$, that should have the same dimension, are then combined with a cross product to obtain the final output, aimed at reproducing $\mathcal{F}(\bm{u}(\bm{y}))$.

Such class of neural networks may also have more than two sub-networks, allowing to separately handle more than two inputs.
We indeed distinguish two type of deep operator networks, one for the turbulence map and one for the closure map.

\begin{figure}[htpb!]
    \centering
    \includegraphics[width=0.9\linewidth]{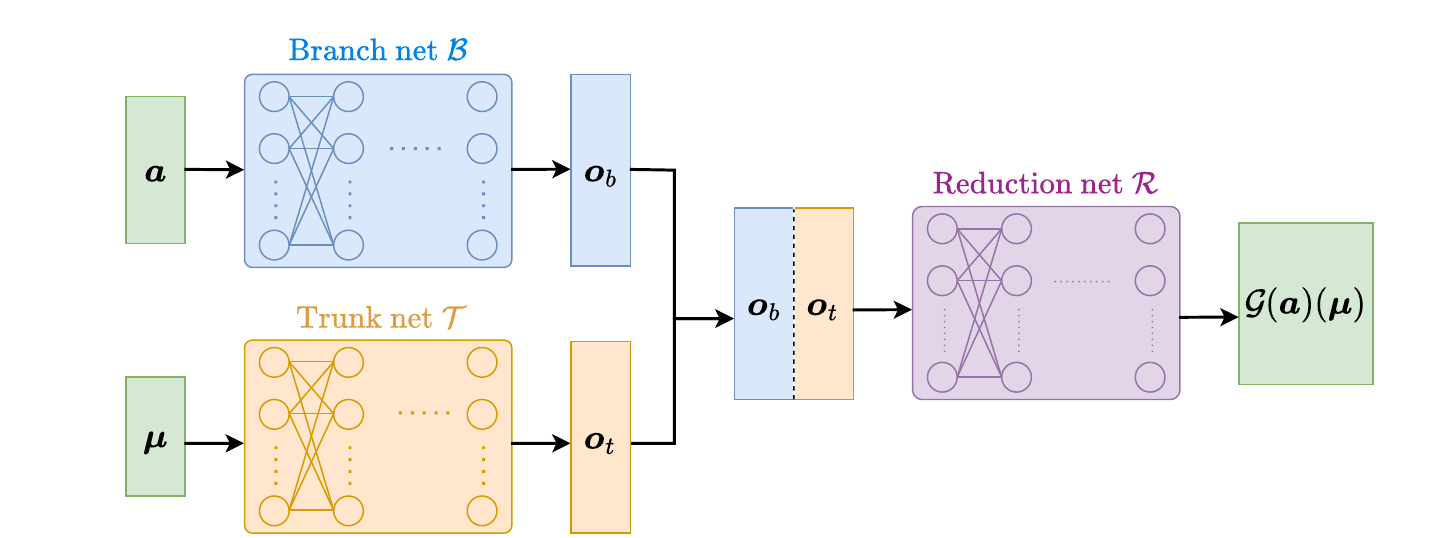}
    \caption{DeepONet architecture used for the eddy viscosity coefficients modeling $\mathcal{G}$.}
    \label{fig:deeponet}
\end{figure}

For the turbulence map $\mathcal{G}(\bm{a}, \bm{\mu})$ we consider a modified version of the \textbf{DeepONet}, whose structure is detailed in Figure \ref{fig:deeponet}. The two inputs are in this case the velocity reduced variables $\bm{a}=\{a_i\}_{i=1}^{N_u}$, and the problem parameters $\bm{\mu}=\{\mu_i\}_{i=1}^{N_P}$. The difference with respect to the DeepONet in Figure \ref{fig:deeponet-simple} is that the two outputs $\bm{o}_b$ and $\bm{o}_t$ are then concatenated and processed through another sub-network, here called \emph{reduction network} $\mathcal{R}$. Such procedure allows to have a multi-dimensional output, in our case $\bm{g}(\bm{a}(\bm{\mu}))=\mathcal{R}([\bm{o}_b, \bm{o}_t])$.
The loss function that the network minimizes is expressed as follows:
\begin{equation}
    \mathcal{L}_{\mathcal{G}} = \dfrac{1}{N_{train}}\sum_{i=1}^{N_{train}} \Bigl(\|\bm{g}_i^{proj}-\mathcal{G}(\bm{a}_i, \bm{\mu}_i)\|^2_2 \Bigr)
    \label{eq:loss-deeponet-turb}
\end{equation}

\begin{figure}[htpb!]
    \centering
    \includegraphics[width=\linewidth]{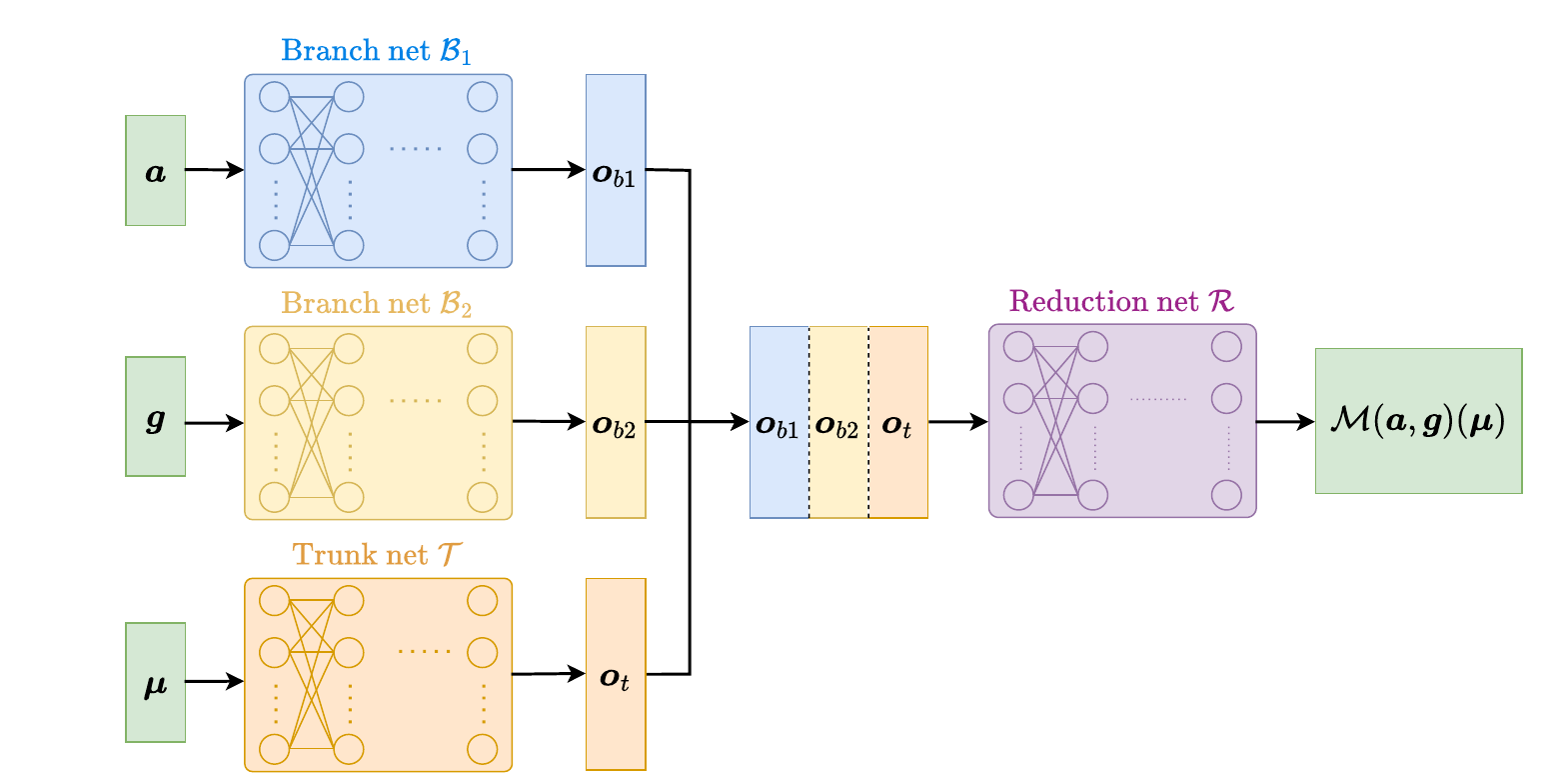}
    \caption{MIONet architecture used for the closure modeling $\mathcal{M}$.}
    \label{fig:mionet}
\end{figure}

The closure map $\mathcal{M}(\bm{a}, \bm{g}, \bm{\mu})$ is instead modeled through a modified version of the Multi-Input deep Operator Network (\textbf{MIONet}), whose structure is displayed in Figure \ref{fig:mionet}.
The architecture relies on three sub-networks separately processing the inputs $\bm{a}=\{a_i\}_{i=1}^{N_u}$, $\bm{g}=\{g_i\}_{i=1}^{N_{\nu_t}}$, and $\bm{\mu}=\{\mu_i\}_{i=1}^{N_P}$. Three outputs $\bm{o}_{b1}$, $\bm{o}_{b2}$ and $\bm{o}_{t}$ may be combined in different ways. In the original MIONet model \cite{jin2022mionet}, the temporary outputs are combined using a tensor product, in our case we use a simplified version and just consider a concatenation of the outputs as input to the reduction network $\mathcal{R}$. The final output is $\bm{\tau}(\bm{a}, \bm{g})(\bm{\mu})=\mathcal{R}([\bm{o}_{b1}, \bm{o}_{b2}, \bm{o}_{t}])$.

For the closure mapping we provide two different training procedures, characterized by different contributions in the loss function. 

Tha standard loss function is:
\begin{equation}
    \mathcal{L}_{\mathcal{M}} =\dfrac{1}{N_{train}}\sum_{i=1}^{N_{train}} \Bigl(\|\bm{\tau}^{exact}_i-\mathcal{M}(\bm{a}_i, \bm{g}_i, \bm{\mu}_i)\|^2_2 \Bigr).
    \label{eq:loss-mionet-stand}
\end{equation}

We consider a modified training with the following loss function:
\begin{equation}
\mathcal{L}^{\star}_{\mathcal{M}} =\mathcal{L}_{\mathcal{M}}+\mathcal{L}_{\mathcal{G}}+\mathcal{L}_{\mathcal{MG}},
    \label{eq:loss-mionet-mod}
\end{equation}
where $\mathcal{L}_{\mathcal{M}}$ is in Equation \eqref{eq:loss-mionet-stand}, $\mathcal{L}_{\mathcal{G}}$ in in Equation \eqref{eq:loss-deeponet-turb}, and $\mathcal{L}_{\mathcal{MG}}$ is computed as follows:
\begin{equation}
    \mathcal{L}_{\mathcal{MG}}=\dfrac{1}{N_{train}}\sum_{i=1}^{N_{train}} \Bigl(\|\bm{\tau}^{exact}_i(\bm{a}_i^{proj}, \mathcal{G}(\bm{a}_i^{proj}, \bm{\mu}_i), \bm{\mu}_i)-\mathcal{M}(\bm{a}_i, \mathcal{G}(\bm{a}_i, \bm{\mu}_i), \bm{\mu}_i)\|^2_2 \Bigr).
    \label{eq:loss-mixed}
\end{equation}
In this case, the weights of the pre-trained network $\mathcal{G}$ and of network $\mathcal{M}$ are optimized at the same time.

In the numerical results' Section, we compare the performance of the two training procedures. 

We introduce here the following notation:
\begin{itemize}
    \item We call \textbf{DD-EV-ROM} the ROM in Equation \eqref{eq:compact-ppe-rom-turb}, where:
    \begin{itemize}
        \item $\mathcal{G}$ minimizes the loss in Equation \eqref{eq:loss-deeponet-turb};
        \item $\mathcal{M}$ minimizes the loss in Equation \eqref{eq:loss-mionet-stand}.
    \end{itemize}
    \item We call \textbf{DD-EV-ROM}$^{\star}$ the ROM in Equation \eqref{eq:compact-ppe-rom-turb}, where:
    \begin{itemize}
        \item $\mathcal{G}$ is pre-trained with loss function \eqref{eq:loss-deeponet-turb}; 
        \item the pre-trained $\mathcal{G}$ and $\mathcal{M}$ are trained together to minimize the modified loss function in Equation \eqref{eq:loss-mionet-mod}. 
    \end{itemize}
\end{itemize}

We also investigated a different training procedure in order to guarantee higher accuracy, especially for test parameters and in extrapolation settings. A more detailed discussion on the training procedure can be found in \ref{subsec:app-training}. Additionally, the hyperparameters' setting of mappings $\mathcal{G}$ and $\mathcal{M}$ are reported in \ref{app:hyperparams}.

\section{Numerical results}
\label{sec:results}
In this paper, we test all the numerical methods presented in \ref{sec:methods} on two different test cases:
\begin{enumerate}
    \item[$(\bm{a})$] the periodic flow past a circular cylinder;
    \item[$(\bm{b})$] the channel-driven cavity flow;
    \item[$(\bm{c})$] the steady backward-facing step flow.
\end{enumerate}

The first two test cases are characterized by physical parameters, namely $t$ and $\nu$, while the third test case has three geometrical parameters.

The present Section is divided into three parts, one for each test case, and in all cases we focus on the following aspects:
\begin{enumerate}
    \item[$(\bm{i})$] Description of the test case, and POD eigenvalues decay;
    \item[$(\bm{ii})$] Comparison among EV-ROM, DD-EV-ROM and DD-EV-ROM$^{\star}$ in terms of relative $L^2$ errors in train and test configurations, in different selected regimes $(N_u, N_p, N_{\nu_t})$, depending on the energy retained by the POD modes.
    In test cases $(\bm{a})$ and $(\bm{b})$ the time trend of the relative $L^2$ error in an unseen configuration is also described.
    \item[$(\bm{iii})$] Graphical representation of the results at the final online time instance for a selected modes combination, for an unseen configuration.
\end{enumerate}

Regarding the offline stage, described in step $(\bm{ii})$, the simulations have been carried out using the open-source software OpenFOAM \cite{jasak1996error, ofsite}, which employs the finite-volume discretization method. The POD operators are extracted using the library ITHACA-FV \cite{ithacasite}, while the reduced-order model integrated with the neural networks is solved using Python and PyTorch \cite{paszke2019pytorch}. 

For what concerns step $(\bm{ii})$, the expression of the time-dependent relative errors used to compare the ROMs performances is computed as follows.
\begin{equation}
    \mathcal{E}_s^{\bm{\mu}_{phys}^{\star}} (t) = \dfrac{\| \tilde{\bm{s}}^{\bm{\mu}_{phys}^{\star}} (t) - \bm{s}^{\bm{\mu}_{phys}^{\star}} (t) \|_{L^2(\Omega)}}{\| \bm{s}^{\bm{\mu}_{phys}^{\star}}(t) \|_{L^2(\Omega)}}\, ,
    \label{eq:l2-errs-exp}
\end{equation}
where $\tilde{\bm{s}}$ is the ROM solution, while $\bm{s}$ is the corresponding FOM reference.
In step $(\bm{ii})$, we will also refer to a global metric, used to compare the performances of DD-EV-ROM, DD-EV-ROM$^{\star}$ compared with the baseline EV-ROM.

In particular, we consider as global metric the relative average gain in the error of the novel data-driven approaches with respect to the standard EV-ROM. Specifically, we consider a field $\bm{s}$, and we call:
\begin{itemize}
    \item $\mathcal{E}_s^{\text{EV-ROM}}(\bm{\mu}^{\star})$ the relative error of the EV-ROM approximation with respect to the FOM field, for parameter $\bm{\mu}^{\star}$, 
    \item $\mathcal{E}_s^{\text{DD-EV-ROM}}(\bm{\mu}^{\star})$ the relative error of the DD-EV-ROM approximation with respect to the FOM field, for parameter $\bm{\mu}^{\star}$,
    \item $\mathcal{E}_s^{\text{DD-EV-ROM}^{\star}}(\bm{\mu}^{\star})$ the relative error of the DD-EV-ROM$^{\star}$  approximation with respect to the FOM field, for parameter $\bm{\mu}^{\star}$.
\end{itemize}
The metrics of interest in the case of unsteady simulations are:
\begin{equation}
\begin{split}
   & g_{\bm{s}} = \dfrac{1}{N_{samples}} \sum_{i=1}^{N_{samples}}\dfrac{(\overline{\mathcal{E}_s^{\text{EV-ROM}}(\bm{\mu}_i)}^t-\overline{\mathcal{E}_s^{\text{DD-EV-ROM}}(\bm{\mu}_i)}^t)}{\overline{\mathcal{E}_s^{\text{EV-ROM}}(\bm{\mu}_i)}^t},\\
  & g_{\bm{s}}^{\star} = \dfrac{1}{N_{samples}} \sum_{i=1}^{N_{samples}}\dfrac{(\overline{\mathcal{E}_s^{\text{EV-ROM}}(\bm{\mu}_i)}^t-\overline{\mathcal{E}_s^{\text{DD-EV-ROM}^{\star}}(\bm{\mu}_i)}^t)}{\overline{\mathcal{E}_s^{\text{EV-ROM}}(\bm{\mu}_i)}^t},
\end{split}
    \label{eq:metric-heatmaps}
\end{equation}
where $\overline{(\cdot)}^t$ denotes a time average in the time window considered, and $N_{samples}$ is the number of train or test samples considered.

In the case of steady simulations, the metrics are simplified as:
\begin{equation}
\begin{split}
   & g_{\bm{s}} = \dfrac{1}{N_{samples}} \sum_{i=1}^{N_{samples}}\dfrac{(\mathcal{E}_s^{\text{EV-ROM}}(\bm{\mu}_i)-{\mathcal{E}_s^{\text{DD-EV-ROM}}(\bm{\mu}_i)})}{{\mathcal{E}_s^{\text{EV-ROM}}(\bm{\mu}_i)}},\\
  & g_{\bm{s}}^{\star} = \dfrac{1}{N_{samples}} \sum_{i=1}^{N_{samples}}\dfrac{({\mathcal{E}_s^{\text{EV-ROM}}(\bm{\mu}_i)}-{\mathcal{E}_s^{\text{DD-EV-ROM}^{\star}}(\bm{\mu}_i)})}{{\mathcal{E}_s^{\text{EV-ROM}}(\bm{\mu}_i)}},
\end{split}
    \label{eq:metric-heatmaps-steady}
\end{equation}

\subsection{Test case (\textbf{a}): periodic flow past a cylinder}
\label{subsec:test-case-a}

The periodic flow past a circular cylinder is a wide-known benchmark in fluid dynamics and in the field of reduced order modeling. In this Section we will focus on a parametric version of this test case, where the parameters considered are time and the dynamic viscosity, that indirectly parametrizes the Reynolds number.

\subsubsection{Offline stage}
\label{subsubsec:fom-a}
\begin{figure}[htpb!]
    \centering
    \subfloat[Domain with notation]{\includegraphics[width=0.5\textwidth]{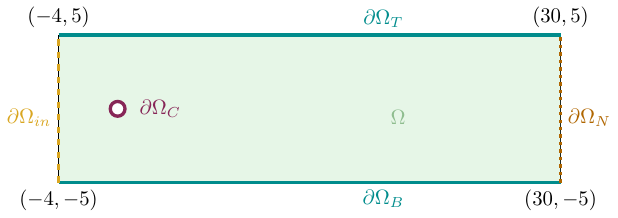}}
    \subfloat[Full order mesh]{\includegraphics[width=0.5\textwidth, trim={0cm 10cm 0cm 10cm}, clip]{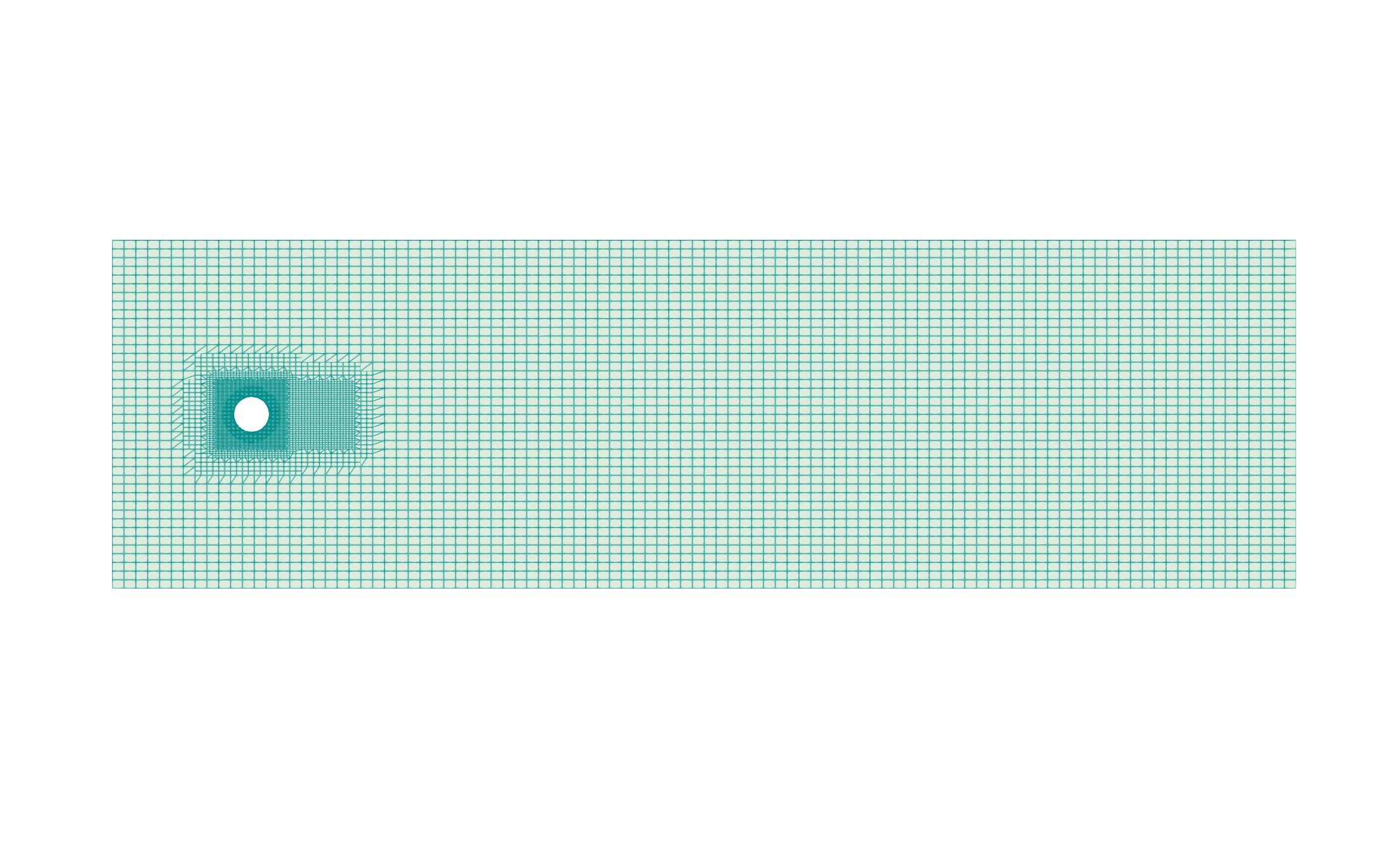}}
    \caption{The domain and full order mesh considered for the periodic flow around a circular cylinder.}
    \label{fig:cyl-domain}
\end{figure}

The domain and mesh for this test case are represented in Figure \ref{fig:cyl-domain}. Following the notation in the above-mentioned Figure, the boundary conditions read as follows.

\begin{equation*}
\text{On }\partial \Omega_{in}:
    \begin{cases}
        \bm{u} = (U_{in}, 0),\\
        \nabla p \cdot \bm{n} = 0;
    \end{cases}
    \quad
    \text{On }\partial \Omega_T \cup \partial \Omega_B:
    \begin{cases}
        \bm{u} \cdot \bm{n} = 0,\\
        \nabla p \cdot \bm{n} = 0;
        
    \end{cases}
\end{equation*}

\begin{equation*}
        \text{On }\partial \Omega_N:
    \begin{cases}
        \nabla \bm{u} \cdot \bm{n} = 0,\\
        p = 0;
        
    \end{cases}
    \quad
            \text{On }\partial \Omega_C:
    \begin{cases}
       \bm{u} = \bm{0},\\
        \nabla p \cdot \bm{n} = 0.
    \end{cases}
\end{equation*}

\begin{figure}[htpb!]
\centering{\includegraphics[width=0.6\textwidth]{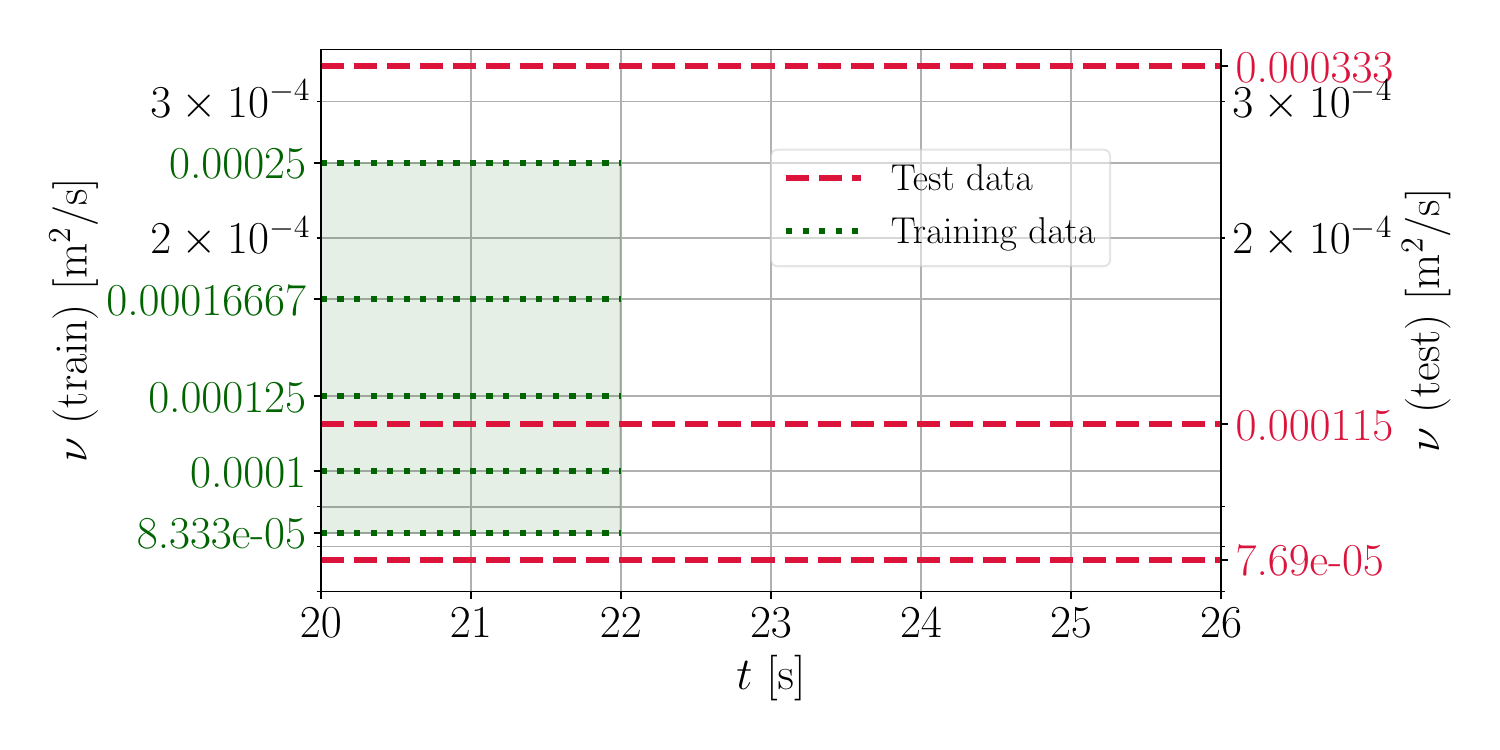}}
    \caption{Graphical representations of the sets of parameters used in the offline stage for test case ($\bm{a}$), in green, and in the online stage, in red. The green shadow describes the accessible offline/training area.}
    \label{fig:params-plot-a}
\end{figure}

The time and parameter ranges considered to run the FOM and collect the POD snapshots are displayed in Figure \ref{fig:params-plot-a}. Moreover, we collect a total amount of $500$ snapshots for each viscosity value, with a total of $2500$ snapshots.
The optimal penalty coefficient $\tau$ considered is $\tau=10^3$, and in the ROMs we consider a second-order time integration scheme.

\begin{figure}[htpb!]
    \centering
    \subfloat[POD cumulative eigenvalues]{\includegraphics[height=4.7cm]{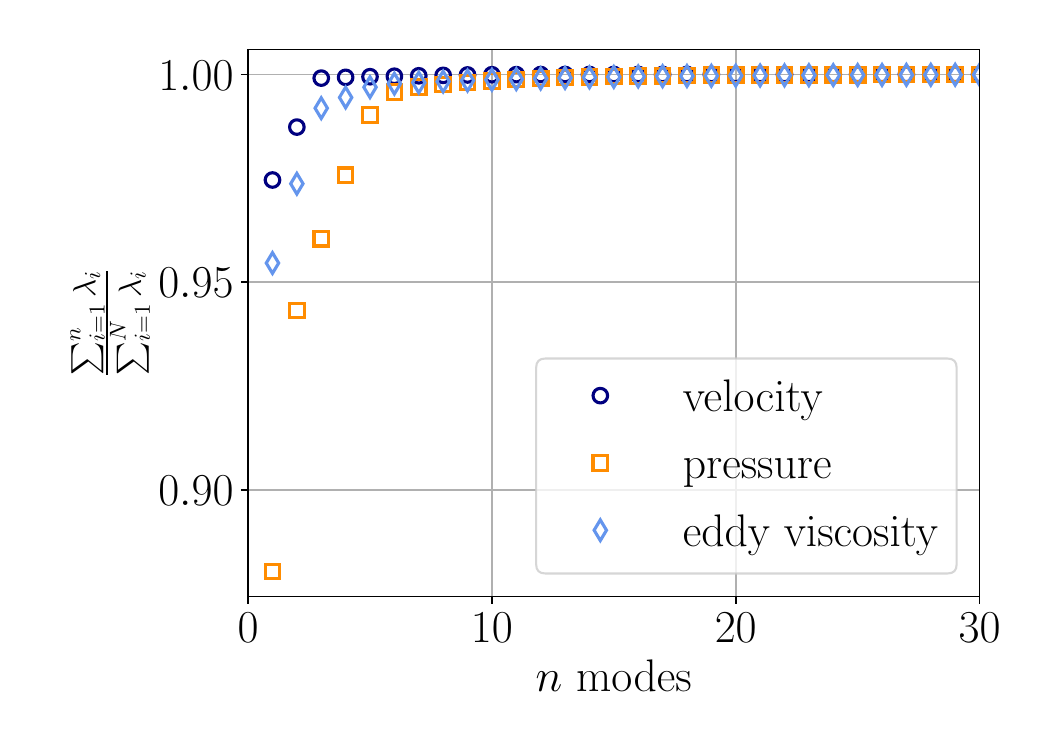}}
    \subfloat[POD eigenvalues decay]{\includegraphics[height=4.7cm]{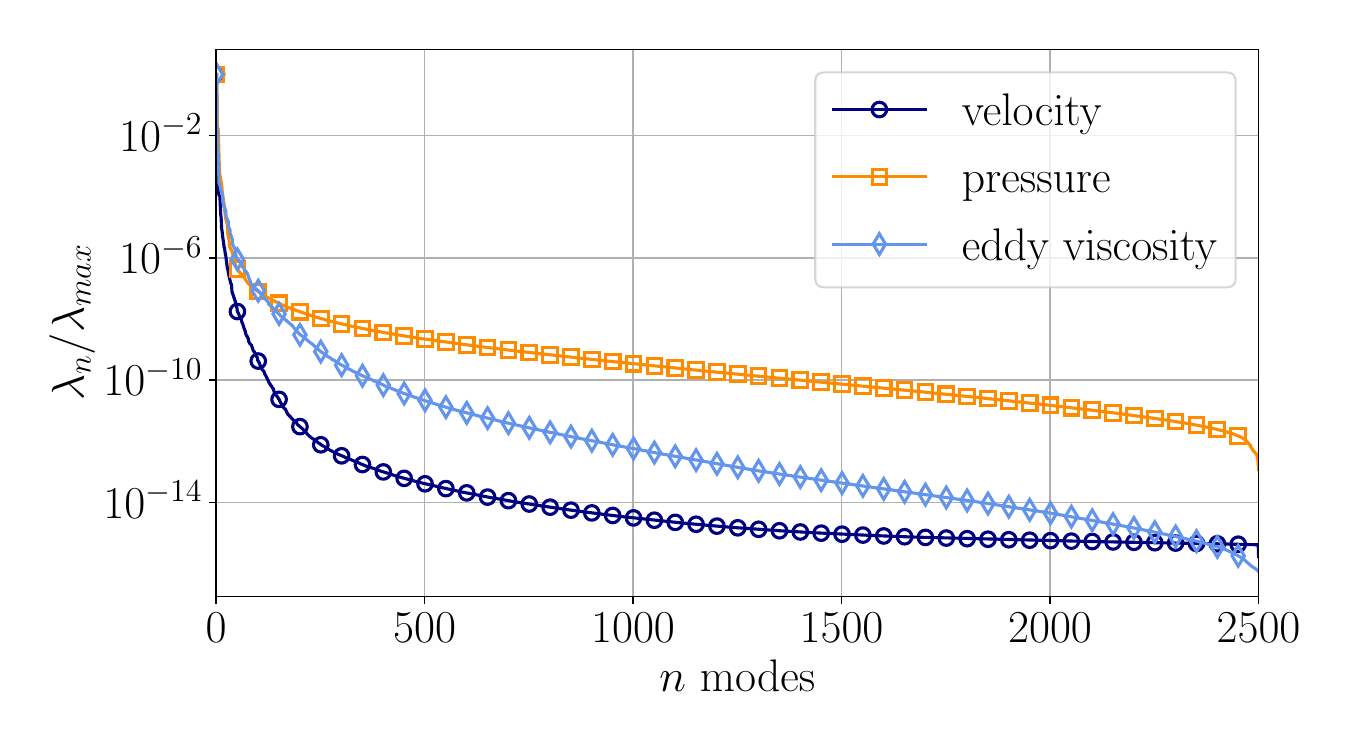}}
    \caption{Cumulative eigenvalues and eigenvalues decay for the test case \textbf{a}.}
    \label{fig:eig-a}
\end{figure}

The snapshots are taken every $\SI{0.004}{\second}$, the POD cumulative eigenvalues and the corresponding decay are represented in Figure \ref{fig:eig-a}. We can see a fast decay, especially for the velocity field. However, as we will see in the following part, the standard POD-Galerkin approach is characterized by lack of accuracy in the \emph{marginally-resolved} regime, even if the number of modes considered is enough to capture the dynamics of the problem.

The following part focus on a quantitative analysis on the relative error for the EV-ROM, DD-EV-ROM, and DD-EV-ROM$^{\star}$ for different combinations of the modes, namely different triples $(N_u, N_p, N_{\nu_t})$. The specific combinations are chosen starting from the eigenvalues' decay, namely from an energy analysis, as specified in Table \ref{tab:modal-regimes-a}.

\begin{table}[htpb!]
    \centering
 \caption{Combinations of the modes for the velocity, pressure and eddy viscosity fields in test case (\textbf{a}).}
    \label{tab:modal-regimes-a}
    \begin{tabular}{>{\centering\arraybackslash}p{0.5\linewidth}
    >{\centering\arraybackslash}p{0.1\linewidth}
    >{\centering\arraybackslash}p{0.1\linewidth}
    >{\centering\arraybackslash}p{0.1\linewidth}
    }
    \toprule
    {\textbf{\emph{Minimum retained energy}} [$\%$]} &$N_u$& $N_p$ & $N_{\nu_t}$\\
    \midrule
    95&1&3&1\\
    99&3&5&3\\
    99.5&3&6&5\\
    99.6&3&7&5\\
    99.7&3&8&6\\
    99.8&3&9&7\\
    99.9&3&12&13\\
    99.99&10&24&30\\
    \bottomrule
    \end{tabular}
\end{table}

\subsubsection{Error analysis}
\label{subsubsec:errs-a}

The first error analysis compares the baseline EV-ROM method with the novel DD-EV-ROM and DD-EV-ROM$^{\star}$. 

\begin{figure}[htpb!]
    \centering
    \includegraphics[width=\textwidth]{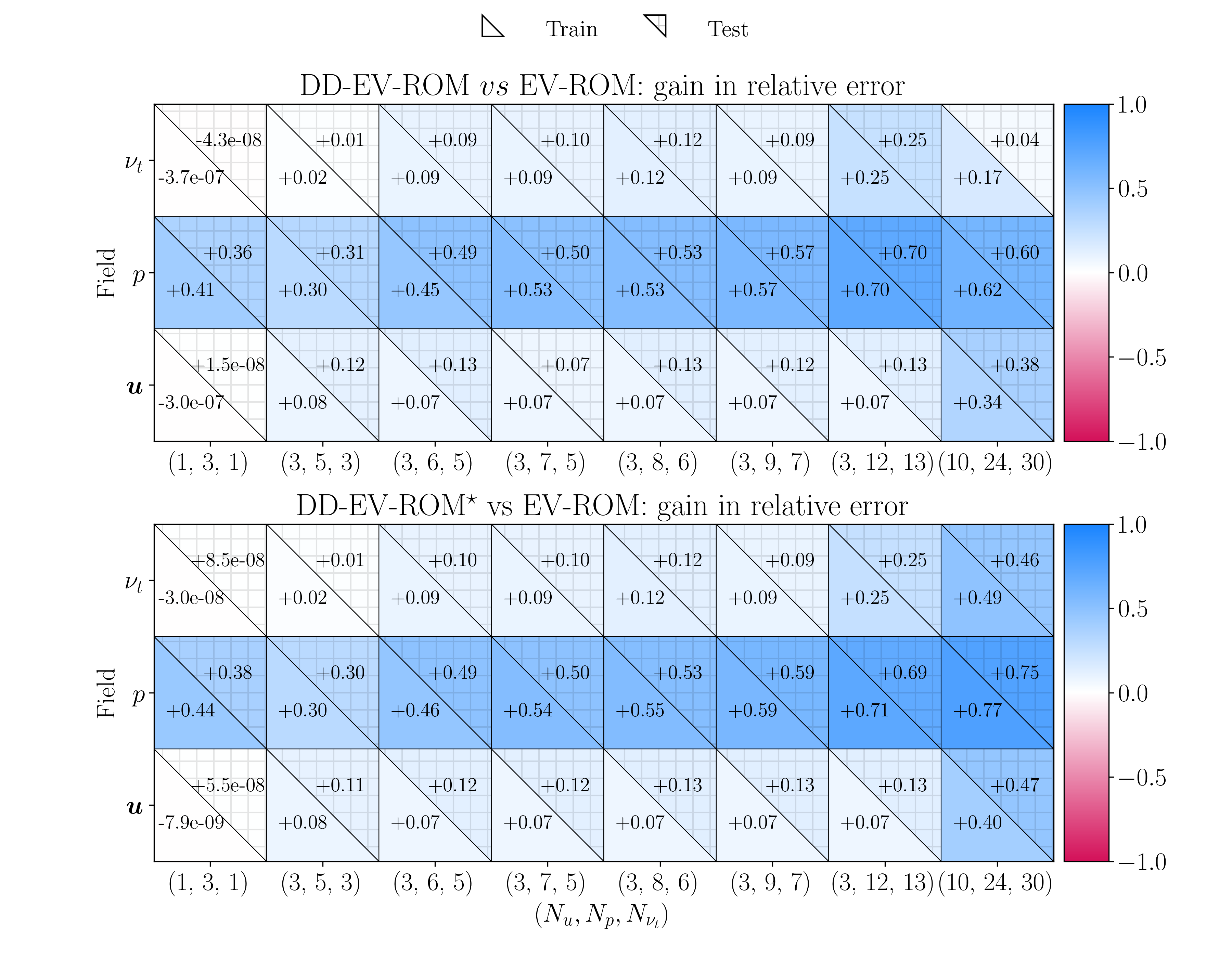}
    \caption{Comparison in the average performance of the DD-EV-ROM and DD-EV-ROM$^{\star}$ with respect to the baseline EV-ROM, for test case \textbf{a}.
    The metric of interest is the time average gain in relative error for train and test viscosities (lower/upper triangles). The error is represented for different modes combination and for the three fields velocity, pressure and eddy viscosity.}
    \label{fig:cylinder-heatmaps}
\end{figure}

The average gain in the relative error \eqref{eq:metric-heatmaps} over the time window $[0, 2]$ is represented in the heatmaps of Figure \ref{fig:cylinder-heatmaps} for the three fields of interest, and for the regimes of Table \ref{tab:modal-regimes-a}. The performances of both train and test configurations are analysed (lower and upper triangles in the plots).
The results show that in general both data-driven approaches outperform the baseline EV-ROM in all the modal regimes considered. In particular:
\begin{itemize}
    \item[$\bullet$] The average gain in the relative error is higher for the pressure fields. This happens because we employ a PPE approach and we introduce a dedicated pressure correction term. Moreover, the pressure approximation with the baseline method is typically worse with respect to the velocity one, and this leads to larger pressure closure terms $\bm{\tau}_p$, which are then easier to learn by a neural network.
    \item[$\bullet$] The values for train and test are similar. This is a good indicator of the robustness of the closure mapping $\mathcal{M}$. It indeed generalizes well, also in extrapolation settings and does not perform overfitting phenomena.
    \item[$\bullet$] The average gain increases as the number of modes (or, alternatively the retained energy of the reduced system) increases. This is not intuitive, since one would expect better results with less modes, say in a regime where the correction term is larger. The reason why this happens is that the baseline EV-ROM method becomes typically more unstable as the reduced dimension increases, leading to inaccurate approximations. This aspect is deeply described in \ref{app:case-a} by Figures \ref{fig:cylinder-errs-1} and \ref{fig:cylinder-errs-2}, representing the relative error trend in time, for a test parameter in two different modal regimes.
\end{itemize}

A different statistical analysis on the proposed methods is also described in the Appendix (Figures \ref{fig:cylinder-violin-1} and \ref{fig:cylinder-violin-2}).

\subsubsection{Graphical results}
\label{subsubsec:graph-a}

This Subsection is dedicated to show the graphical approximations of the EV-ROM, DD-EV-ROM and DD-EV-ROM$^{\star}$, comparing them with the FOM reference fields, when $(N_u, N_p, N_{\nu_t})=(10, 24, 30)$. In particular, we show the absolute errors for an unseen parameter at the final time instance of the online simulation for the pressure field (Figure \ref{fig:cyl-graphical-p-err}).

The graphical errors for the other fields are showed in the Appendix (Figures \ref{fig:cyl-graphical-u-err} and \ref{fig:cyl-graphical-nut-err}). 
The results highlight both the gain in the accuracy of the DD-EV-ROM$^{\star}$, which is graphically very close to the projection error. Moreover, for all the fields of interest it can be immediately that the DD-EV-ROM$^{\star}$ is significantly more accurate than the DD-EV-ROM$^{\star}$.

\begin{figure}[htpb!]
    \centering
    \subfloat[]{
    \includegraphics[width=0.5\linewidth]{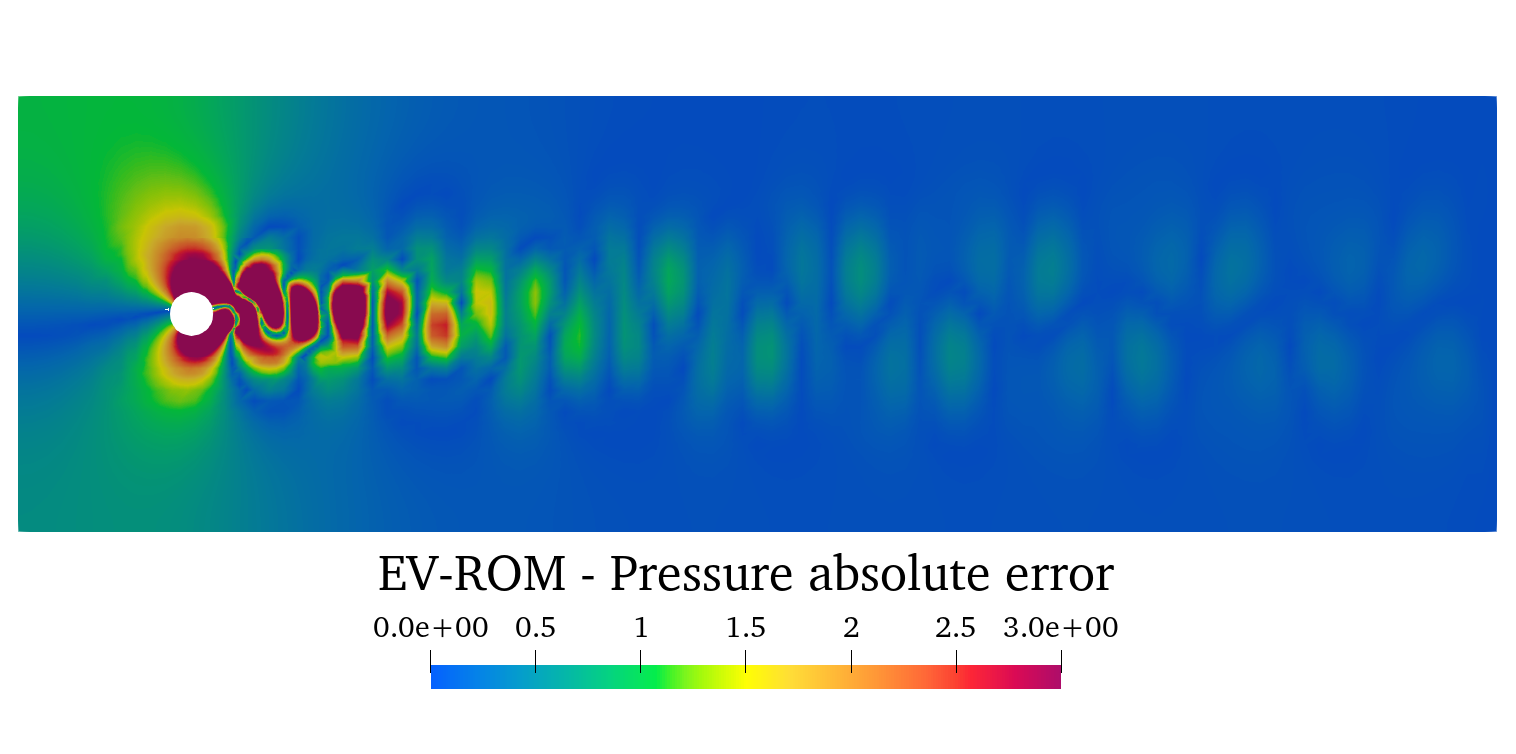}}
    \subfloat[]{
    \includegraphics[width=0.5\linewidth]{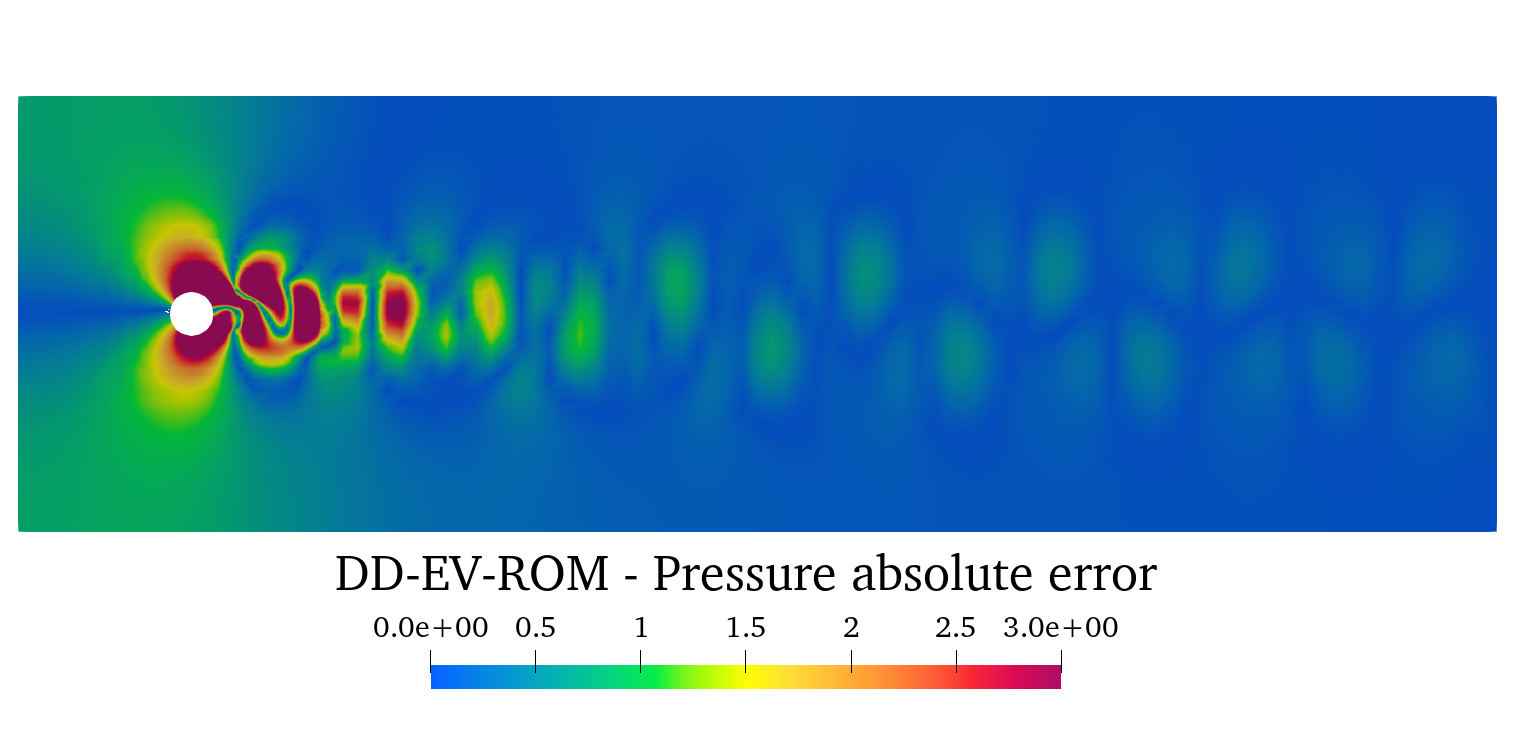}}\\
    \subfloat[]{
    \includegraphics[width=0.5\linewidth]{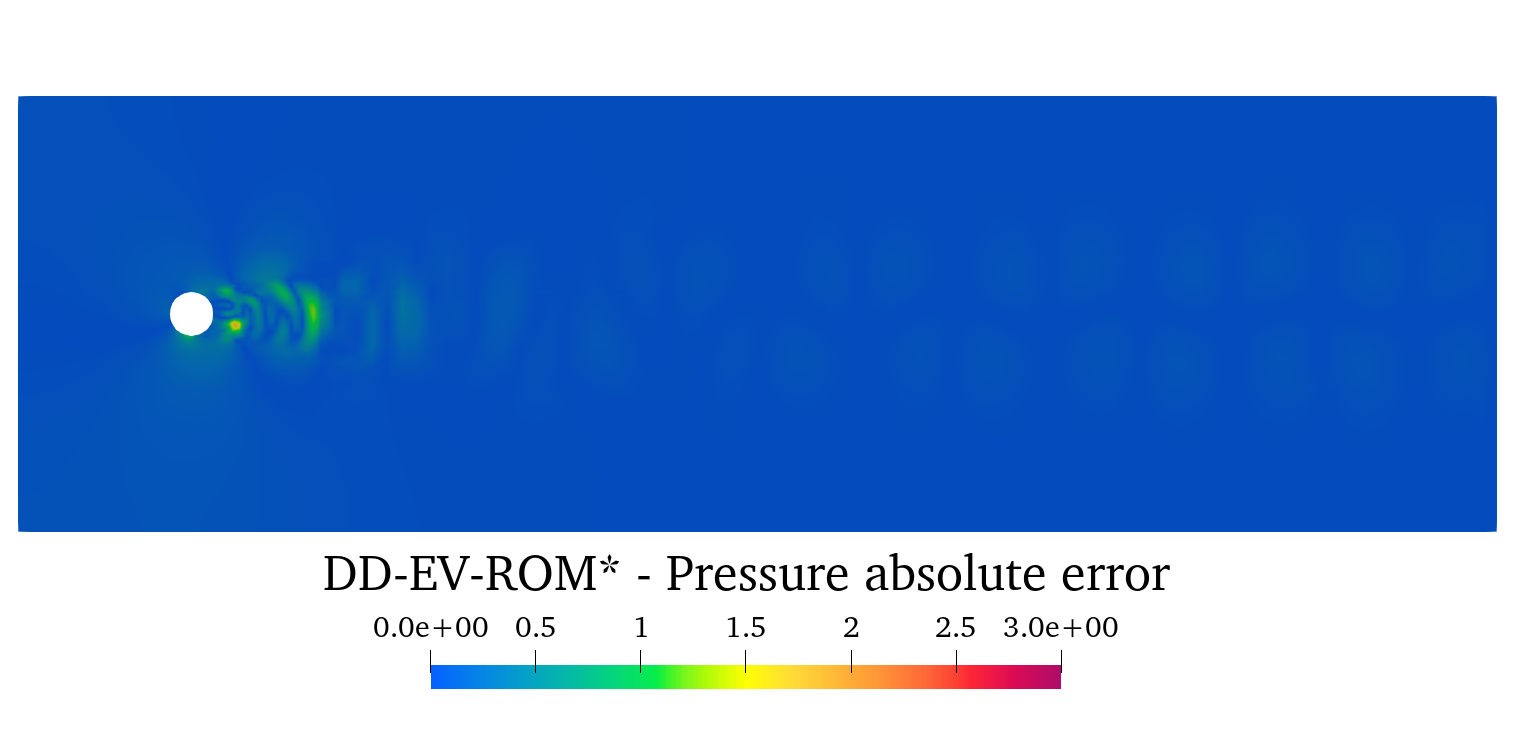}}
    \subfloat[]{
    \includegraphics[width=0.5\linewidth]{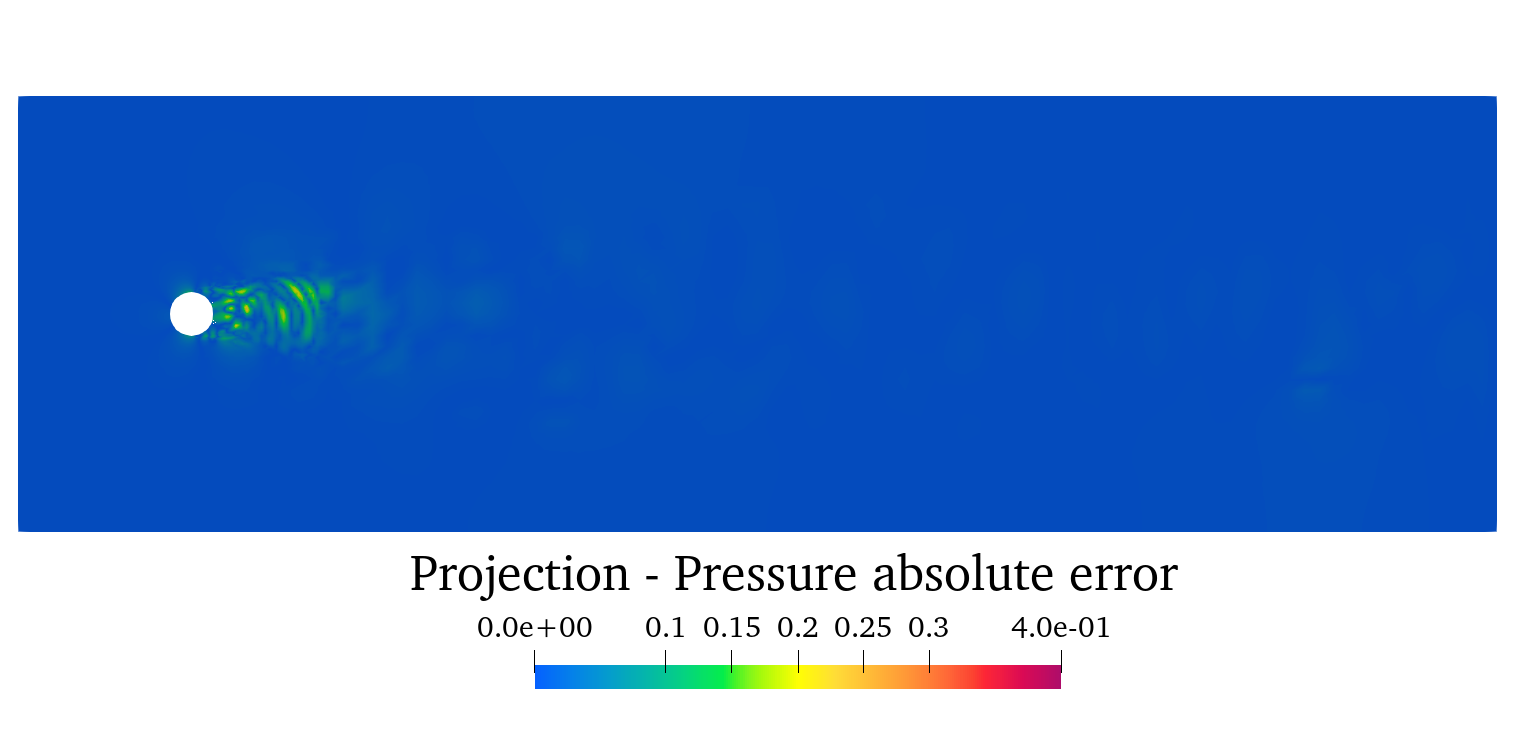}}
    \caption{Graphical pressure absolute error at $\nu=1.15e-4 \, \frac{m^2}{s}$ at the final time instance ($t=8$), of EV-ROM, DD-EV-ROM, DD-EV-ROM$^{\star}$, and for the projected field. The modal regime is $(N_u, N_p, N_{\nu_t})=(10, 24, 30)$.}
    \label{fig:cyl-graphical-p-err}
\end{figure}

\subsection{Test case (\textbf{b}): flow in a channel-driven cavity}
\label{subsec:test-case-b}
The unsteady case of the flow inside a cavity-shaped channel is characterized by a different behavior with respect to the first test case described in \ref{subsec:test-case-a}. Also in this case, we consider a parametrized setup, with parameters time and the viscosity $\nu$.

\subsubsection{Offline stage}
\label{subsubsec:fom-b}

\begin{figure}[htpb!]
    \centering
    \subfloat[Domain with notation]{\includegraphics[width=0.5\textwidth]{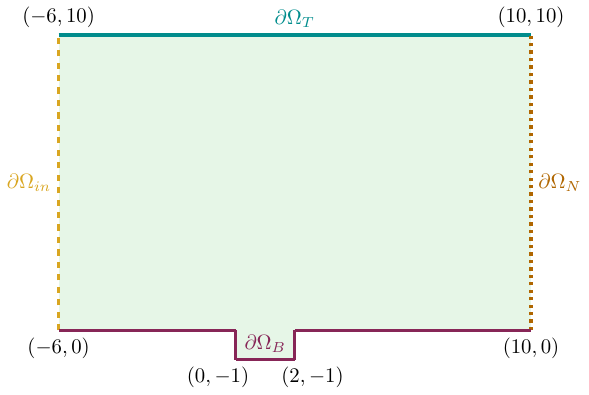}}
    \subfloat[Full order mesh]{\includegraphics[width=0.5\textwidth, trim={4.5cm 0 4.5cm 0cm}, clip]{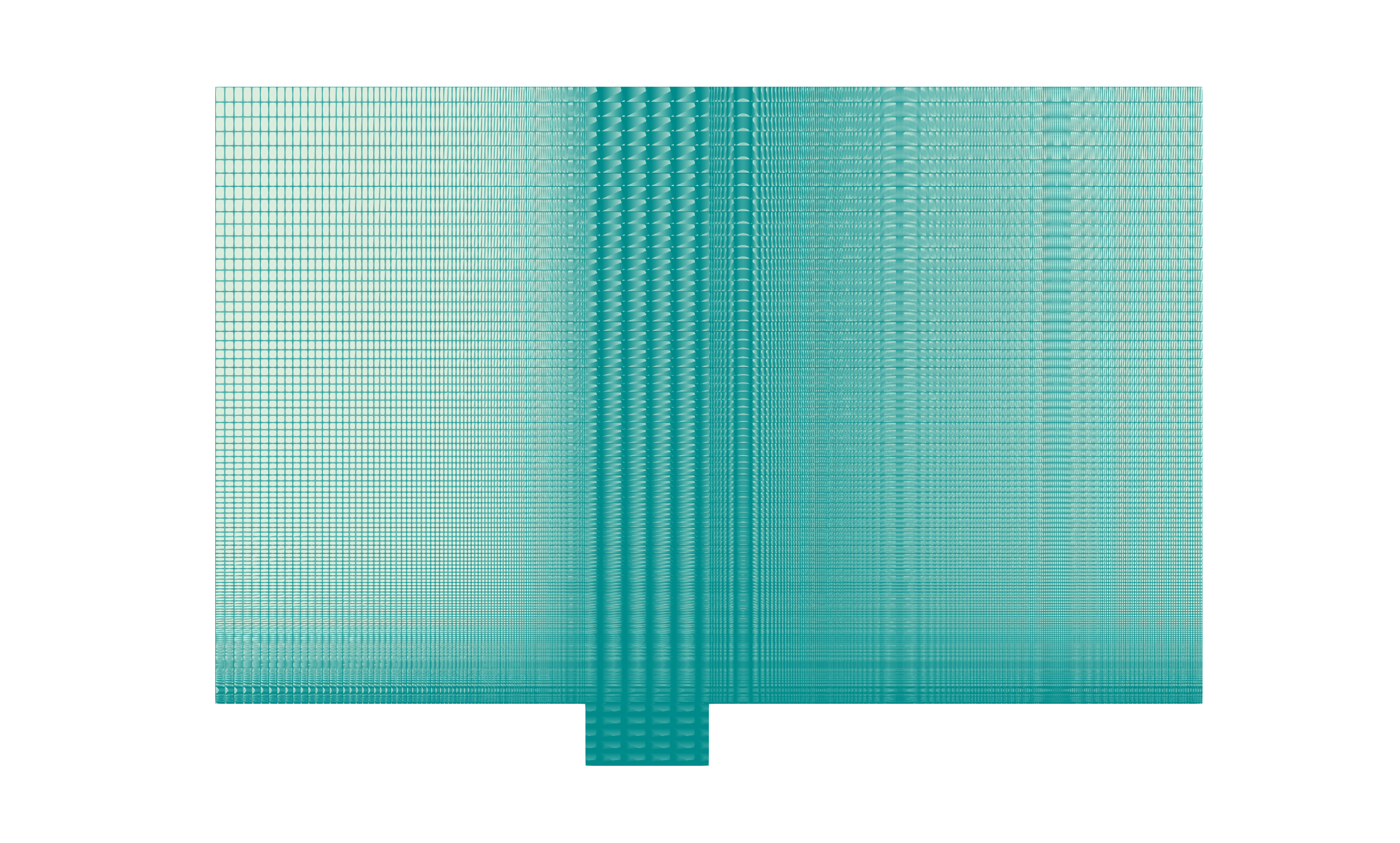}}
    \caption{The domain and full order mesh considered for the channel-driven flow test case.}
    \label{fig:cavity-domain}
\end{figure}

In Figure \ref{fig:cavity-domain} we represent the domain and the mesh considered in the FOM analysis.

Employing the notation in Figure \ref{fig:cavity-domain}, we consider the following boundary conditions: 
\begin{equation*}
    \text{On }\partial \Omega_T:
    \begin{cases}
        \nabla \bm{u} \cdot \bm{n} = 0,\\
        \nabla p \cdot \bm{n} = 0;
    \end{cases}
    \quad
        \text{On }\partial \Omega_B:
    \begin{cases}
       \bm{u}  = \bm{0},\\
       \nabla p \cdot \bm{n} = 0
        
    \end{cases}
\end{equation*}

\begin{equation*}
\text{On }\partial \Omega_{in}:
    \begin{cases}
        \bm{u} = (U^{b}_{in}, 0),\\
        \nabla p \cdot \bm{n} = 0;
    \end{cases}
    \quad
            \text{On }\partial \Omega_N:
    \begin{cases}
       \nabla \bm{u} \cdot \bm{n} = 0,\\
         p = 0.
    \end{cases}
\end{equation*}

\begin{figure}[htpb!]
\centering{\includegraphics[width=0.6\textwidth]{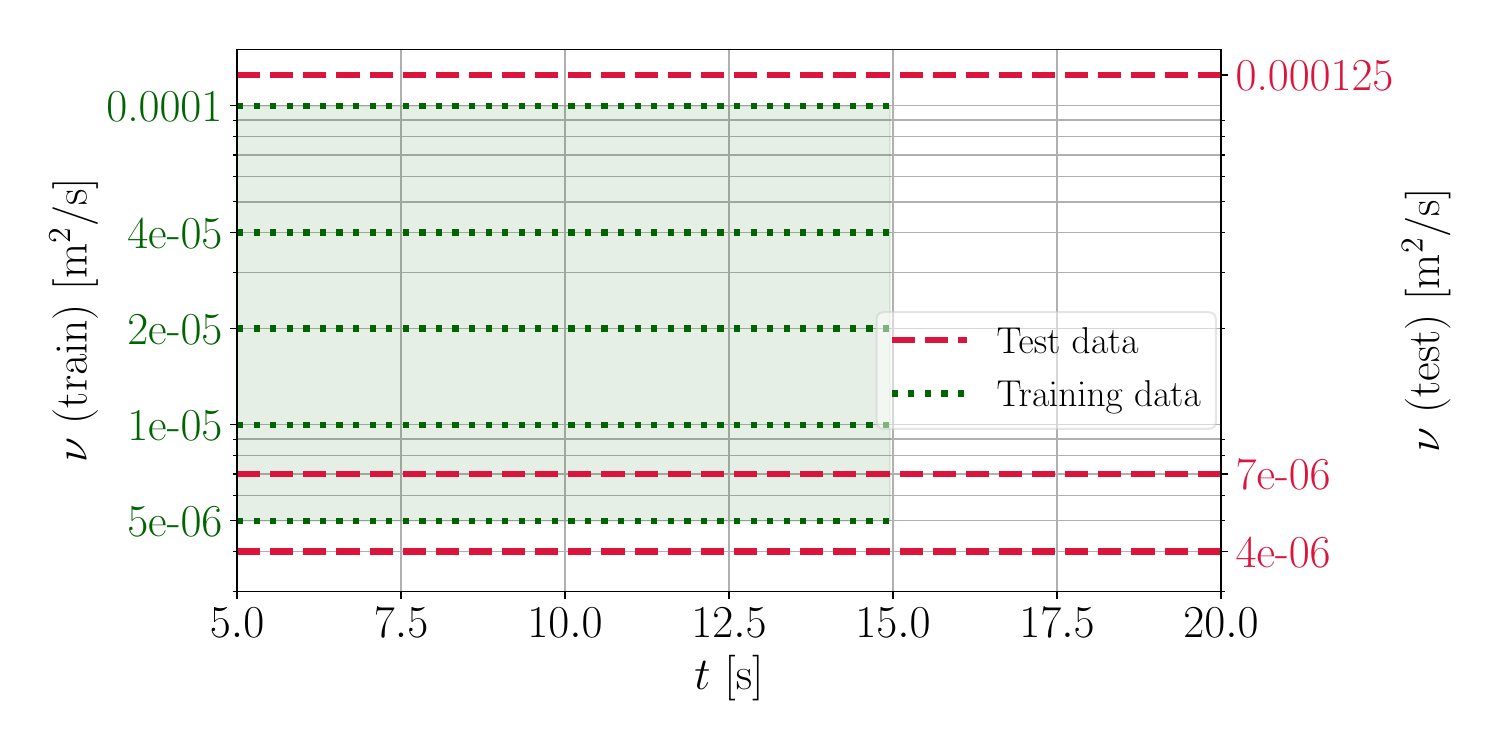}}
    \caption{Graphical representations of the sets of parameters used in the offline stage for test case ($\bm{b}$), in green, and in the online stage, in red. The green shadow describes the accessible offline/training area.}
    \label{fig:params-plot-b}
\end{figure}

For the POD we consider the viscosity range specified in Figure \ref{fig:params-plot-b}, with a corresponding Reynolds number in range $[\num{1e3}, \num{2e4}]$.
The time snapshots are retained every $\num{0.05}$ seconds in the time interval of 10 seconds, for a total amount of snapshots of $1000$, namely $200$ for each FOM simulation.

The penalty coefficient is $\tau=1$ for this test case, while the time integration scheme considered is of first-order, consistent with the one employed in OpenFOAM.

\begin{figure}[htpb!]
    \centering
    \subfloat[POD cumulative eigenvalues]{\includegraphics[height=4.7cm]{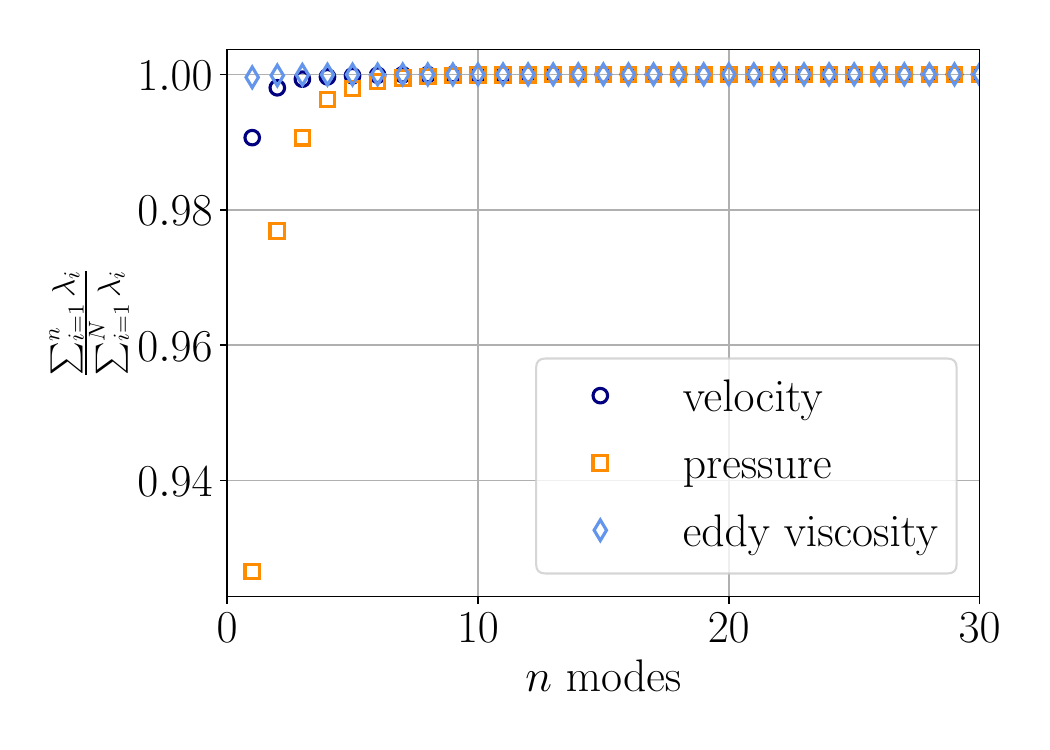}}
    \subfloat[POD eigenvalues decay]{\includegraphics[height=4.7cm]{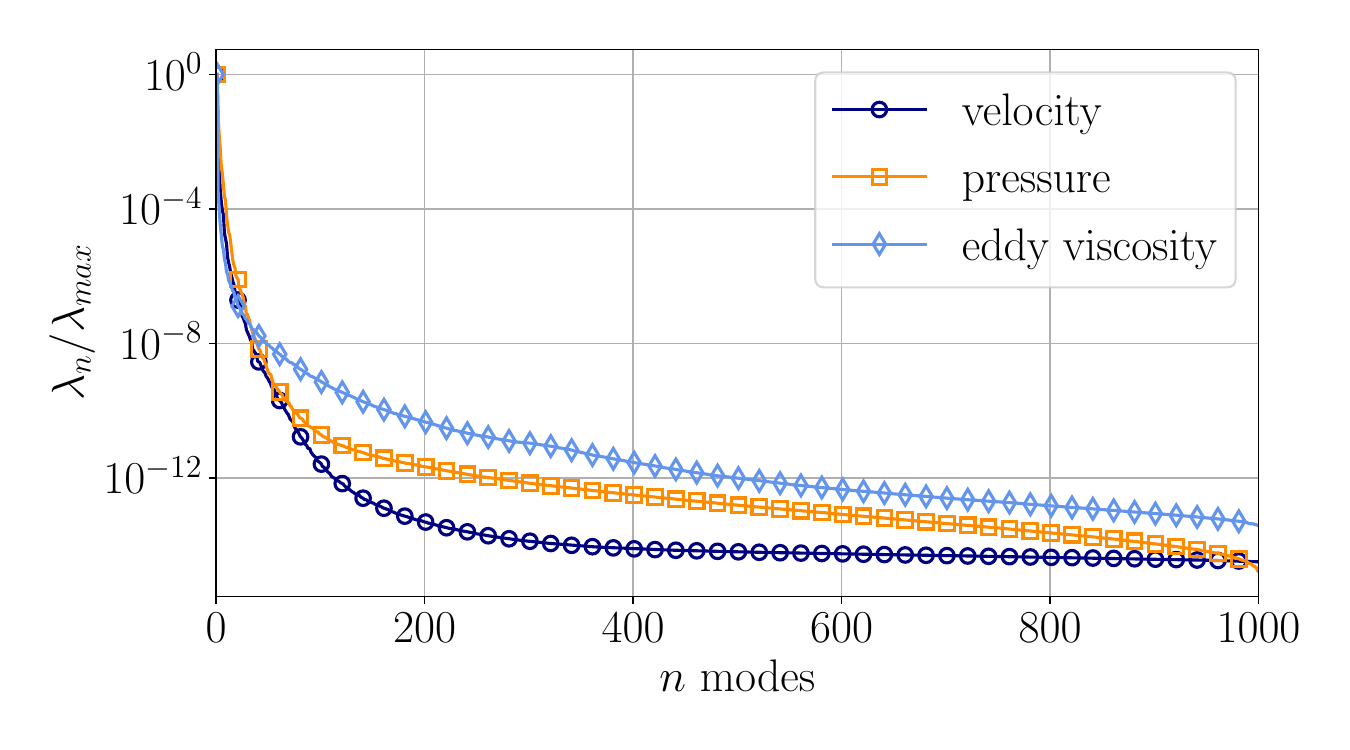}}
    \caption{Cumulative eigenvalues and eigenvalues decay for test case \textbf{b}.}
    \label{fig:eig-b}
\end{figure}

Figure \ref{fig:eig-b} shows the cumulative eigenvalues and the eigenvalues decay of the POD for all the fields considered. The Figure shows a slower decay than in the first test case (Figure \ref{fig:eig-a}). Indeed, this second test case is not periodic leading to a more challenging setting.

The following part is dedicated to a quantitative error analysis, comparing the results obtained with different modal combinations, select through energy criteria and specified in Table \ref{tab:modal-regimes-b}.

\begin{table}[htpb!]
    \centering
 \caption{Combinations of the modes for the velocity, pressure and eddy viscosity fields in test case (\textbf{b}).}
    \label{tab:modal-regimes-b}
    \begin{tabular}{>{\centering\arraybackslash}p{0.5\linewidth}
    >{\centering\arraybackslash}p{0.1\linewidth}
    >{\centering\arraybackslash}p{0.1\linewidth}
    >{\centering\arraybackslash}p{0.1\linewidth}
    }
    \toprule
    {\textbf{\emph{Minimum retained energy}} [$\%$]} &$N_u$& $N_p$ & $N_{\nu_t}$\\
    \midrule
    95&1&2&1\\
    99&1&3&1\\
    99.5&2&4&1\\
    99.7&2&5&1\\
    99.9&3&7&1\\
    \bottomrule
    \end{tabular}
\end{table}

\subsubsection{Error analysis}
\label{subsubsec:errs-b}

\begin{figure}[htpb!]
    \centering
    \includegraphics[width=0.9\textwidth]{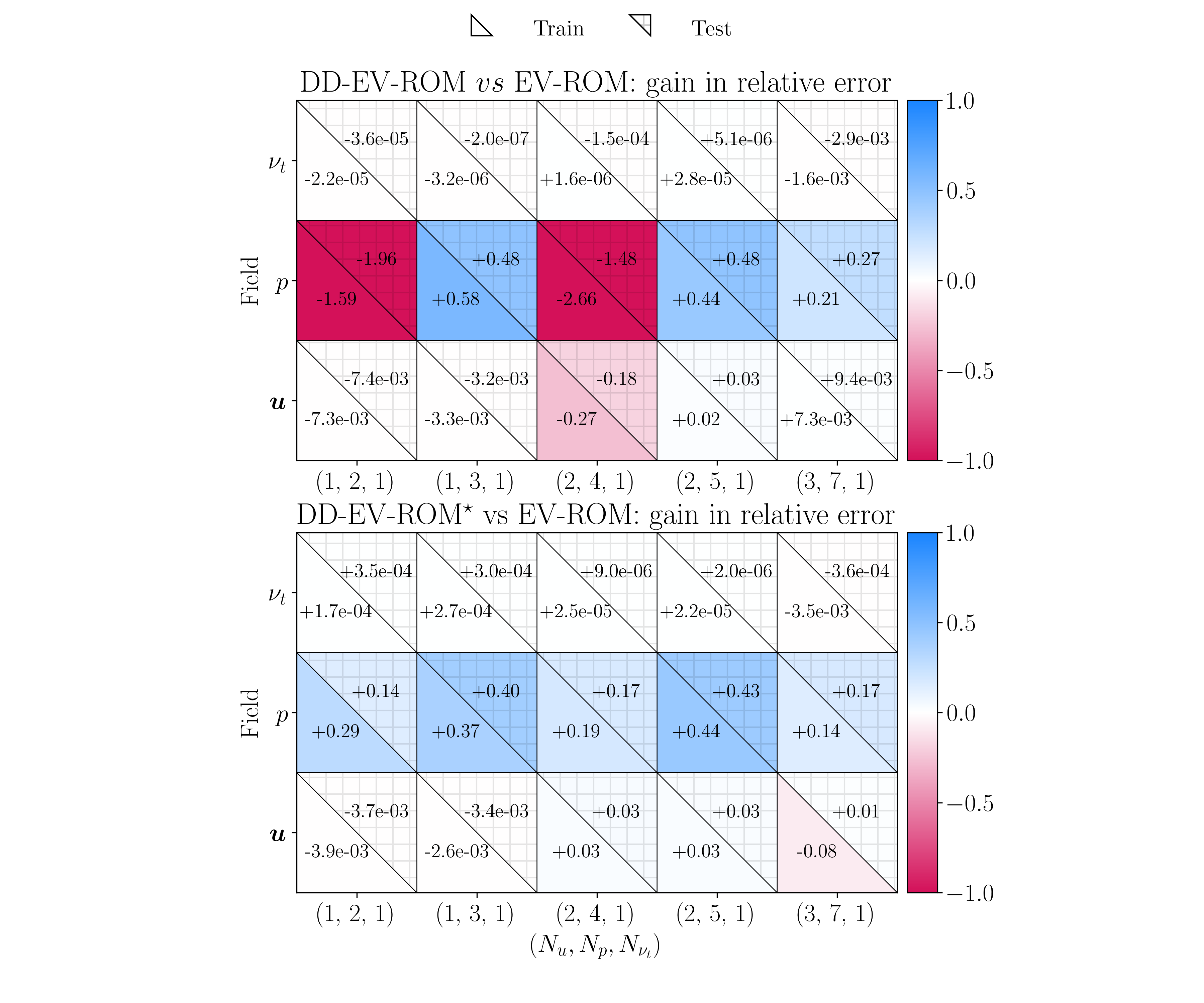}
    \caption{Comparison in the average performance of the DD-EV-ROM and DD-EV-ROM$^{\star}$ with respect to the baseline EV-ROM, for test case (\textbf{b}).
    The metric of interest is the time average gain in relative error for train and test viscosities (lower/upper triangles). The error is represented for different modes combination and for the three fields velocity, pressure and eddy viscosity.}
    \label{fig:cavity-heatmaps}
\end{figure}

The global error analysis is performed in Figure \ref{fig:cavity-heatmaps}, which show the gain metric of Equation \eqref{eq:metric-heatmaps} for the velocity, pressure and eddy viscosity fields, for different modal combinations and for train and test configuration.

As above-mentioned, this test case is more challenging than the cylinder one, since it is unsteady, turbulent, and not periodic.
We can indeed draw the following considerations:
\begin{itemize}
    \item DD-EV-ROM$^{\star}$ always provides improved approximations with respect to EV-ROM.
    \item The main difference with respect to the previous test case is that methods DD-EV-ROM and DD-EV-ROM$^{\star}$ perform significantly differently. In particular, DD-EV-ROM$^{\star}$ always provides improved approximations, while the DD-EV-ROM shows negative error gains for some under-resolved modal regimes, namely its prediction accuracy is lower than the EV-ROM one.
    \item The results obtained for the velocity and for the eddy viscosity are similar to the results of EV-ROM, since the gain is close to 0. This is confirmed by the time trends and by a statistical analysis, that can be found in \ref{app:case-b}. The reason is that the EV-ROM predictions are very close to the projection, which is the best result we can achieve. Hence, there is no hope to improve the EV-ROM result, if it coincides with the projection.
    \item Train and test results are always similar, also when negative gains occur, as in the case of the standard DD-EV-ROM. It happens indeed that the train errors are large even if the mapping $\mathcal{M}$ reaches small loss values. The reason is that, even if the reduced velocity solution is close to the projected field, the problem may be ill-conditioned. Therefore, small perturbations $\|\bm{a}^{sol}-\bm{a}^{proj}\|_2$ may lead to large output differences $\|\tau^{exact}(\bm{a}^{sol}, \mathcal{G}(\bm{a}^{sol}, \bm{\mu}), \bm{\mu})-\tau^{exact}(\bm{a}^{proj}, \bm{g}^{proj}, \bm{\mu})\|_2$, and the \emph{standard} training procedure of $\mathcal{M}$ is not robust enough to capture the correct mapping.
\end{itemize} 

From Figure \ref{fig:cavity-heatmaps}, we can conclude the positive performance of DD-EV-ROM$^{\star}$ over DD-EV-ROM.
Moreover, there are modal regimes where the DD-EV-ROM slightly outperforms the DD-EV-ROM$^{\star}$, which are deeply analysed in \ref{app:case-b}.

\subsubsection{Graphical results}
\label{subsubsec:graph-b}

\begin{figure}[htpb!]
    \centering
    \subfloat[]{
    \includegraphics[width=0.5\linewidth]{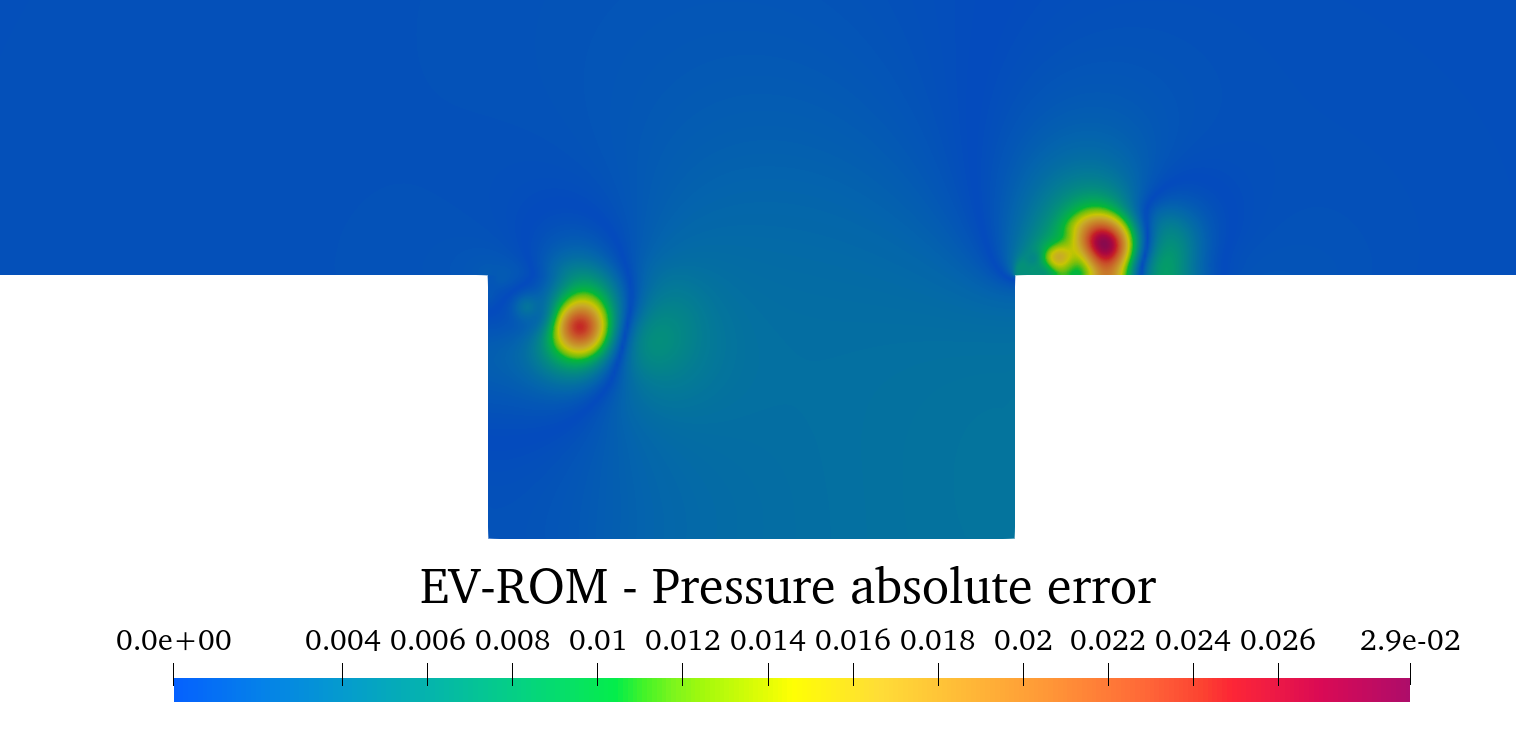}}
    \subfloat[]{
    \includegraphics[width=0.5\linewidth]{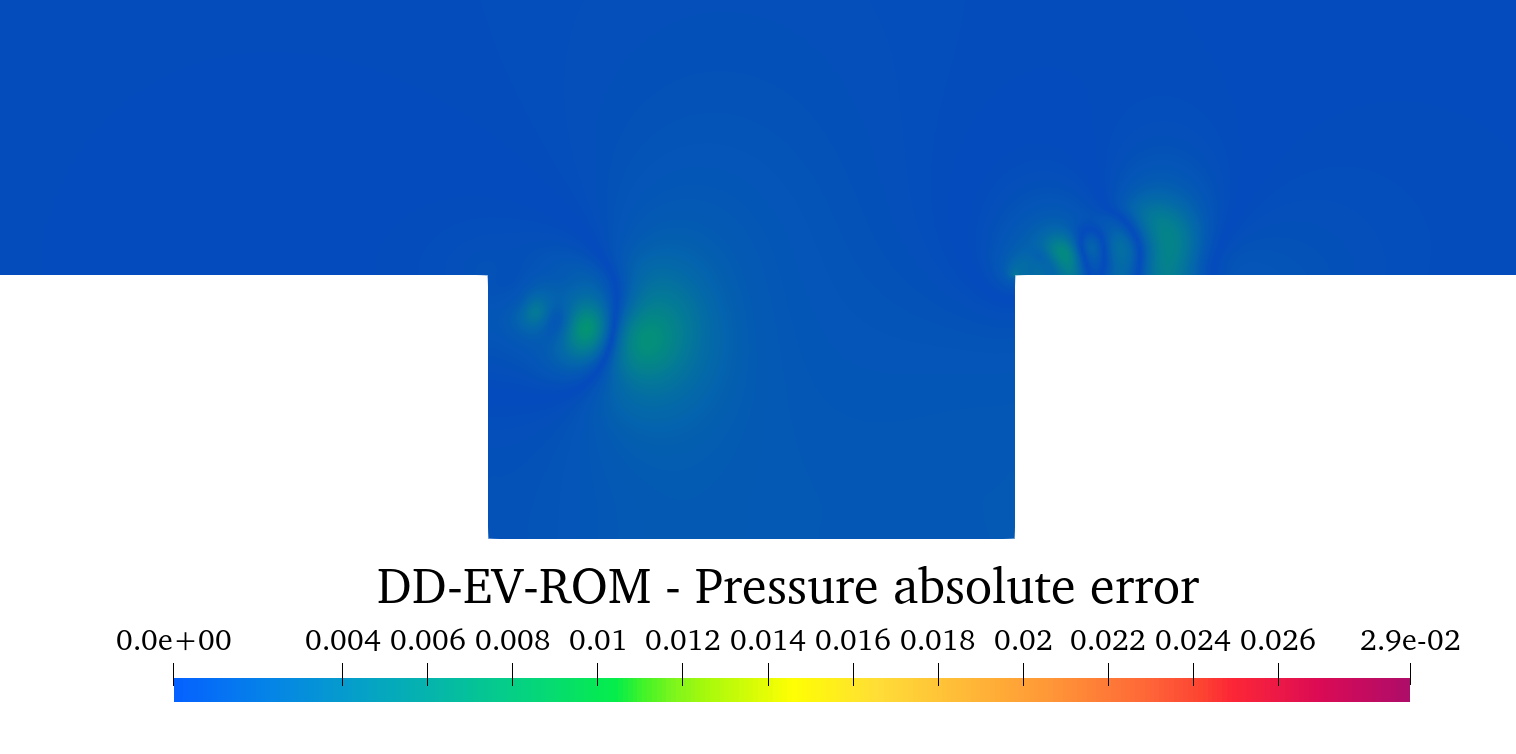}}\\
    \subfloat[]{
    \includegraphics[width=0.5\linewidth]{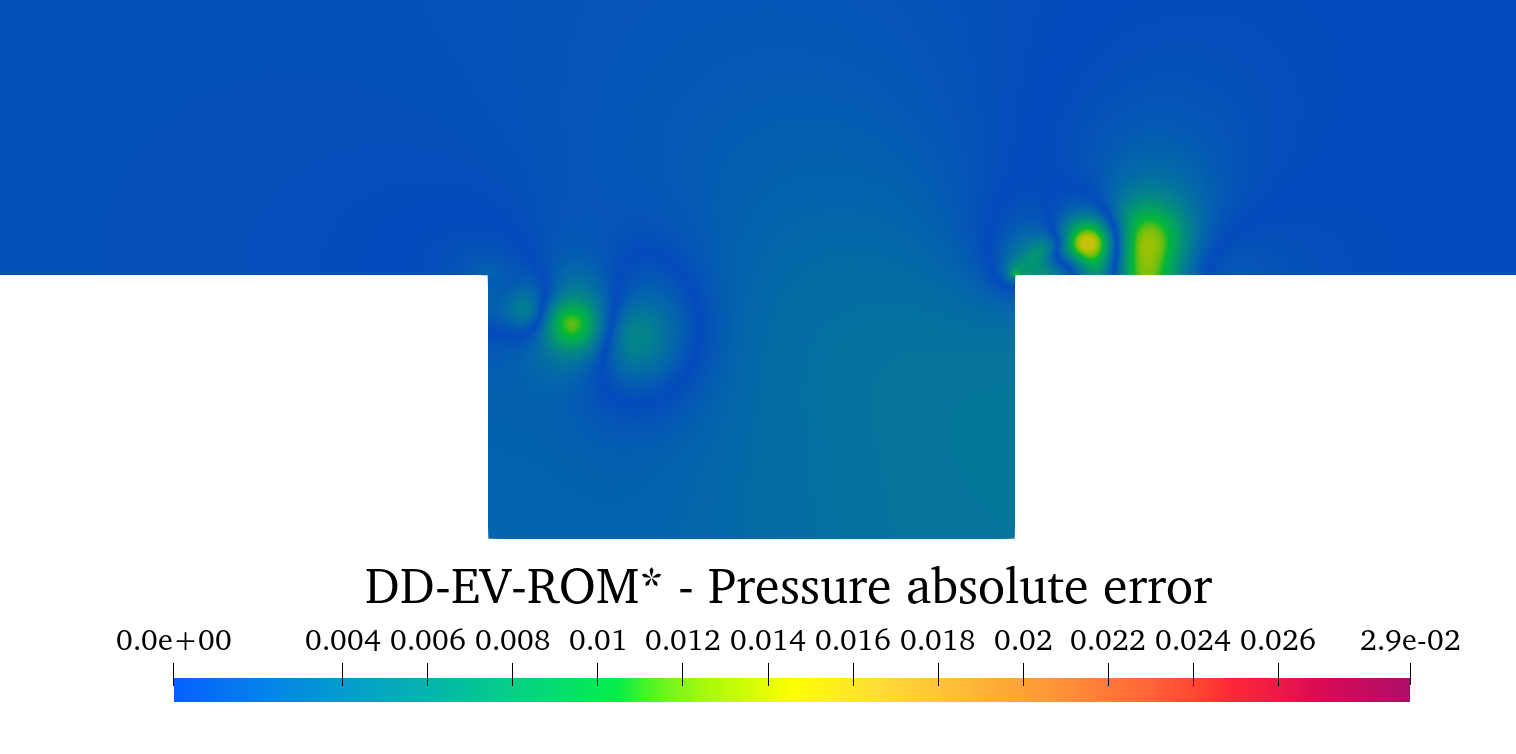}}
    \subfloat[]{
    \includegraphics[width=0.5\linewidth]{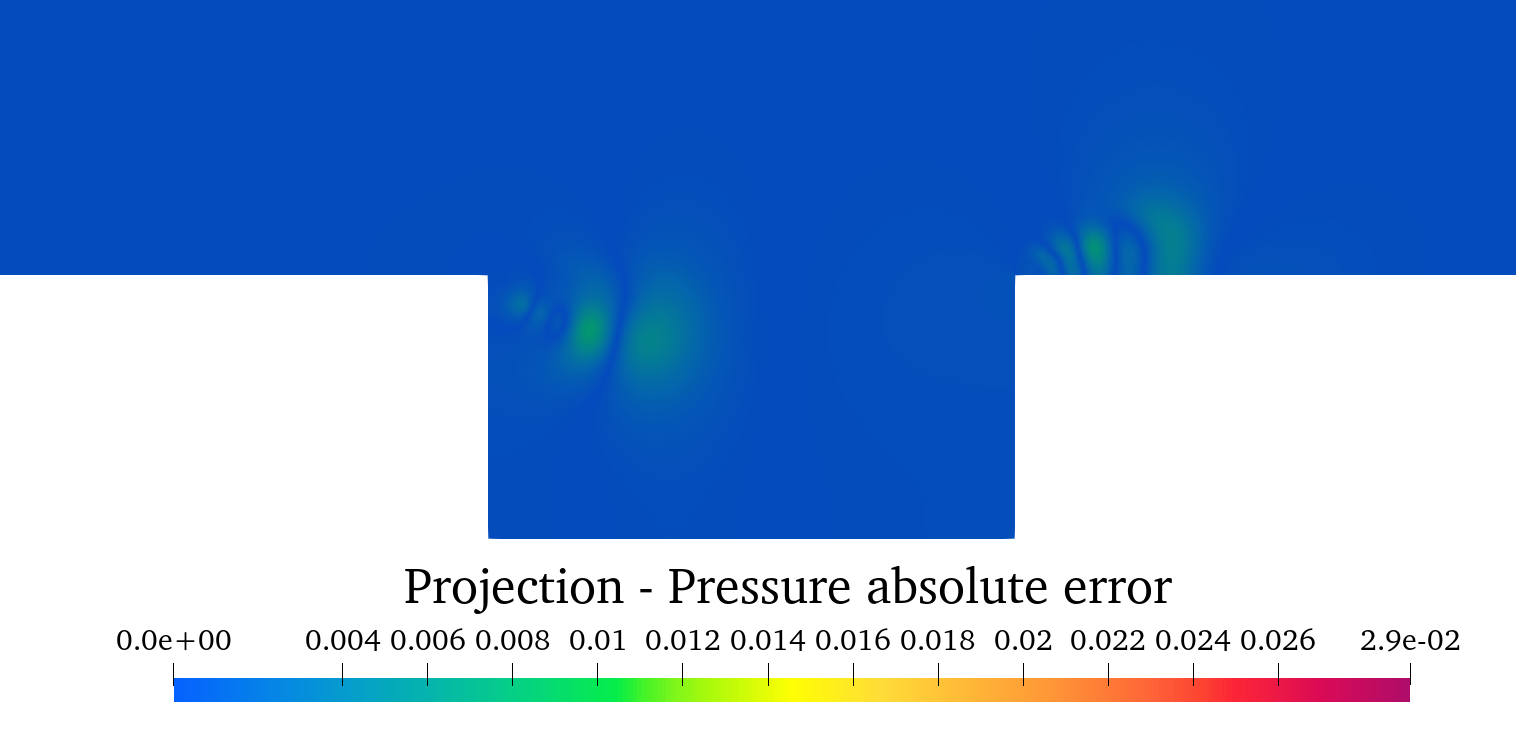}}
    \caption{Graphical pressure absolute error at $\nu=7e-6 \, \frac{m^2}{s}$ at the final time instance ($t=15$), of EV-ROM, DD-EV-ROM, DD-EV-ROM$^{\star}$, and for the projected field. The modal regime is $(N_u, N_p, N_{\nu_t})=(2, 5, 1)$. The Figure shows only a ``zoom" close to the cavity, since the remaining part of the domain exhibits error 0.}
    \label{fig:cav-graphical-p-err}
\end{figure}

We focus here only on the pressure graphical performance since the other fields have similar reconstructions in all the regimes considered.
Figure \ref{fig:cav-graphical-p-err} represents the absolute error between the ROM simulations and the FOM reference, for one test viscosity at the final time instance of the online simulation.
The Figure shows a zoomed region of the domain nearby the cavity for $(N_u, N_p, N_{\nu_t})=(2, 5, 1)$. It shows that the DD-EV-ROM and DD-EV-ROM$^{\star}$ improve the results of EV-ROM.

In particular, the approximation error of DD-EV-ROM is qualitatively similar to the reconstruction error, which is the best result we can achieve with our method.

\begin{comment}
\begin{figure}[htpb!]
    \centering
    \subfloat[]{
    \includegraphics[width=0.5\linewidth]{img/cavity/pturb_1_3_1.png}}
    \subfloat[]{
    \includegraphics[width=0.5\linewidth]{img/cavity/pcorrsimple_1_3_1.png}}\\
    \subfloat[]{
    \includegraphics[width=0.5\linewidth]{img/cavity/pcorrturb_1_3_1.png}}
    \subfloat[]{
    \includegraphics[width=0.5\linewidth]{img/cavity/pproj_1_3_1.png}}\\
    \subfloat[]{
    \includegraphics[width=0.5\linewidth]{img/cavity/pfom.png}}
    \caption{Graphical pressure field at $\nu=7e-6 \, \frac{m^2}{s}$ at the final time instance ($t=15$), of EV-ROM, DD-EV-ROM, DD-EV-ROM$^{\star}$, and for the projected field. The modal regime is $(N_u, N_p, N_{\nu_t})=(1, 3, 1)$. The Figure shows a ``zoom" close to the cavity.}
    \label{fig:cav-graphical-p}
\end{figure}
\end{comment}

\subsection{Test case (\textbf{c}): flow over a geometrically-parameterized backstep}
\label{subsec:test-case-c}
The backward-facing step test case is characterized by a geometrical parameterization, with three parameters, namely the step slop $\alpha$, the inlet height $h_1$, and the global height of the domain $h_2$, represented in Figure \ref{fig:backstep-domain}.
Differently from the other two test cases, the flow is steady. Therefore, we focus on the error analysis of the stationary solution with respect to the FOM counterpart.

\subsubsection{Offline stage}
\label{subsubsec:fom-c}

\begin{figure}[htpb!]
    \centering
    \subfloat[Domain with notation]{\includegraphics[width=0.5\textwidth]{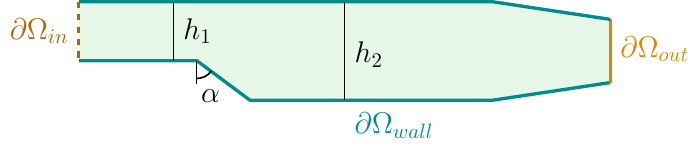}}
    \subfloat[Full order mesh]{\includegraphics[width=0.5\textwidth, trim={55cm 0 55cm 0cm}, clip]{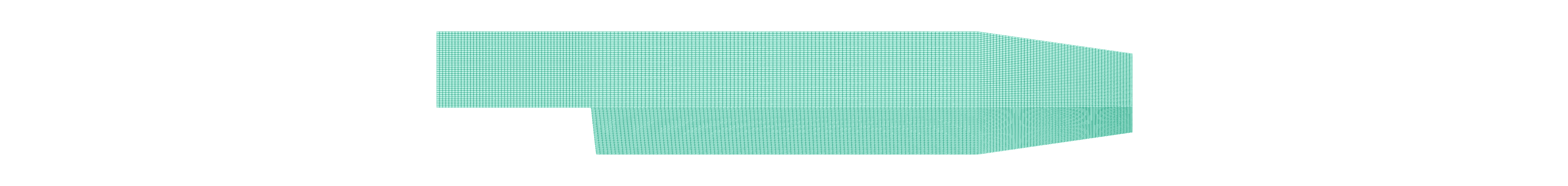}}
    \caption{An example of the domain and full order mesh considered for the backstep flow, for one set of parameters.}
    \label{fig:backstep-domain}
\end{figure}

Figure \ref{fig:backstep-domain} displays the domain and the grid used for a specific set of geometrical parameters.
The offline simulations are run considering as initial conditions $\bm{u}=(1, 0, 0)$ and $p=0$ in the whole domain, while the boundary conditions read as follows.

\begin{equation*}
    \text{On }\partial \Omega_{in}:
    \begin{cases}
        \bm{u} = (1, 0, 0),\\
        \nabla p \cdot \bm{n} = 0;
    \end{cases}
    \quad
        \text{On }\partial \Omega_{out}:
    \begin{cases}
       \nabla \bm{u}  \cdot \bm{n} = \bm{0},\\
       p = 0;
        
    \end{cases}
\end{equation*}

\begin{equation*}
    \text{On }\partial \Omega_{wall}:
    \begin{cases}
        \bm{u} = \bm{0},\\
        \nabla p \cdot \bm{n} = 0;
    \end{cases}
\end{equation*}

\begin{figure}[htpb!]
\centering{\includegraphics[width=0.7\textwidth]{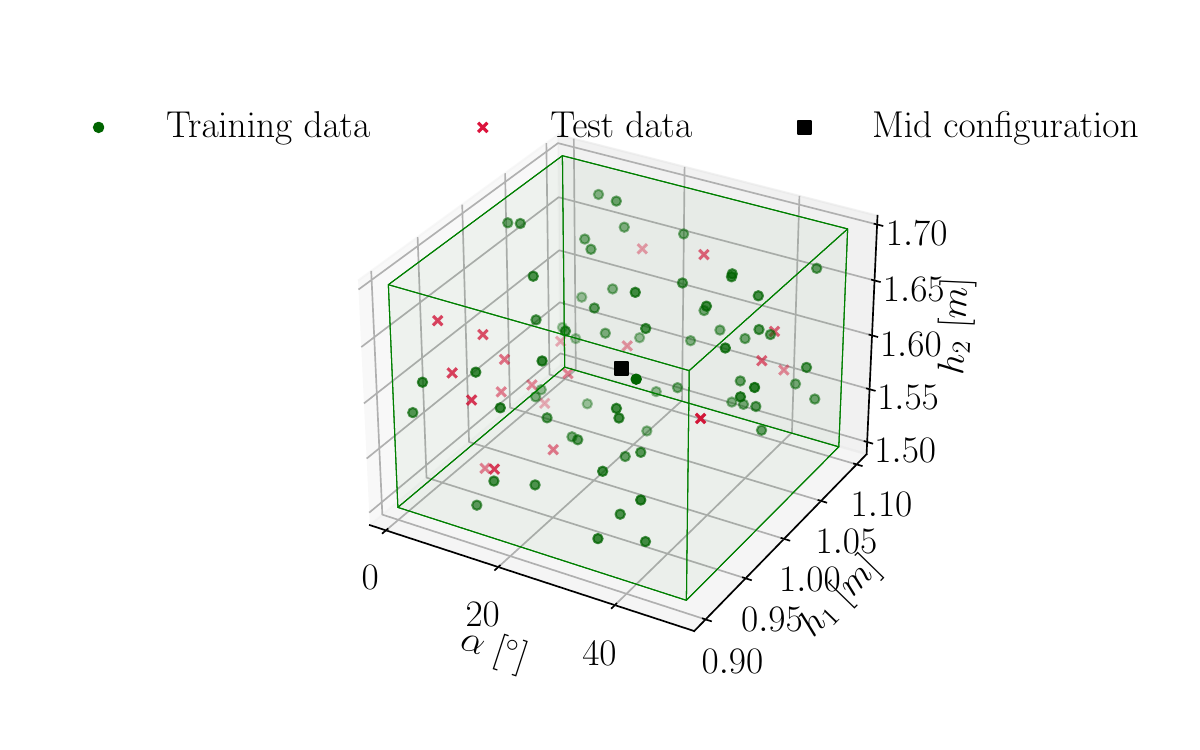}}
    \caption{Graphical representations of the sets of parameters used in the offline stage for test case ($\bm{c}$), in green, and in the online stage, in red. The green shadow describes the accessible offline/training area, while the black square is the mid-configuration, where the modes are projected.}
    \label{fig:params-plot-c}
\end{figure}

We consider a specific set of training parameters for the POD $\{\bm{\mu}_i\}_{i=1}^{N_{\mu}}$, where $\bm{\mu}_i=(\alpha_i, h_1, h_2)$. $N_{\mu}=70$ is the number of simulations with deformed mesh, while we test the approach on $20$ test parameters.
The train, test, and mid-configuration parameters are displayed in Figure \ref{fig:params-plot-c}.

However, we specify that the POD matrices $\mathcal{S}_u$, $\mathcal{S}_p$, $\mathcal{S}_{\nu_t}$ include not only the stationary fields, but also the intermediate states, following the same logic of \cite{zancanaro2021hybrid}. Hence, considering a generic field $\bm{s}$, we have the corresponding snapshots' matrix:
\begin{equation*}
\mathcal{S}_s=\{\bm{s}^{(1)}(\bm{x}_1,\bm{\mu}_1), ...,\bm{s}^{(T_1)}(\bm{x}_1,\bm{\mu}_1), \bm{s}(\bm{x}_1,\bm{\mu}_1), \dots, \bm{s}^{(1)}(\bm{x}_{N_{\mu}},\bm{\mu}_{N_{\mu}}), ...,\bm{s}^{(T_{N_{\mu}})}(\bm{x}_{N_{\mu}},\bm{\mu}_{N_{\mu}}), \bm{s}(\bm{x}_{N_{\mu}},\bm{\mu}_{N_{\mu}}))\},
\end{equation*}
where $\bm{s}^{(j)}(\bm{x}_i,\bm{\mu}_i)$ indicates the $j$-th intermediate state for the field $\bm{s}$ at the deformed mesh $\bm{x}_i$ corresponding to parameter $\bm{\mu}_{i}$.
This strategy is here employed to ensure good ROM accuracy. 
Indeed, considering only the final steady states may lead to instability at the ROM level, since the reduced system (as the full order counterpart) visits different intermediate solutions.
On the other hand, if we consider all the intermediate solutions to compute the modes and the reduced operator, we increase the robustness and stability of the reduced system.
The global number of snapshots, including all the intermediate states, is $\sim 11700$.

\begin{figure}[htpb!]
    \centering
    \subfloat[POD cumulative eigenvalues]{\includegraphics[height=4.7cm]{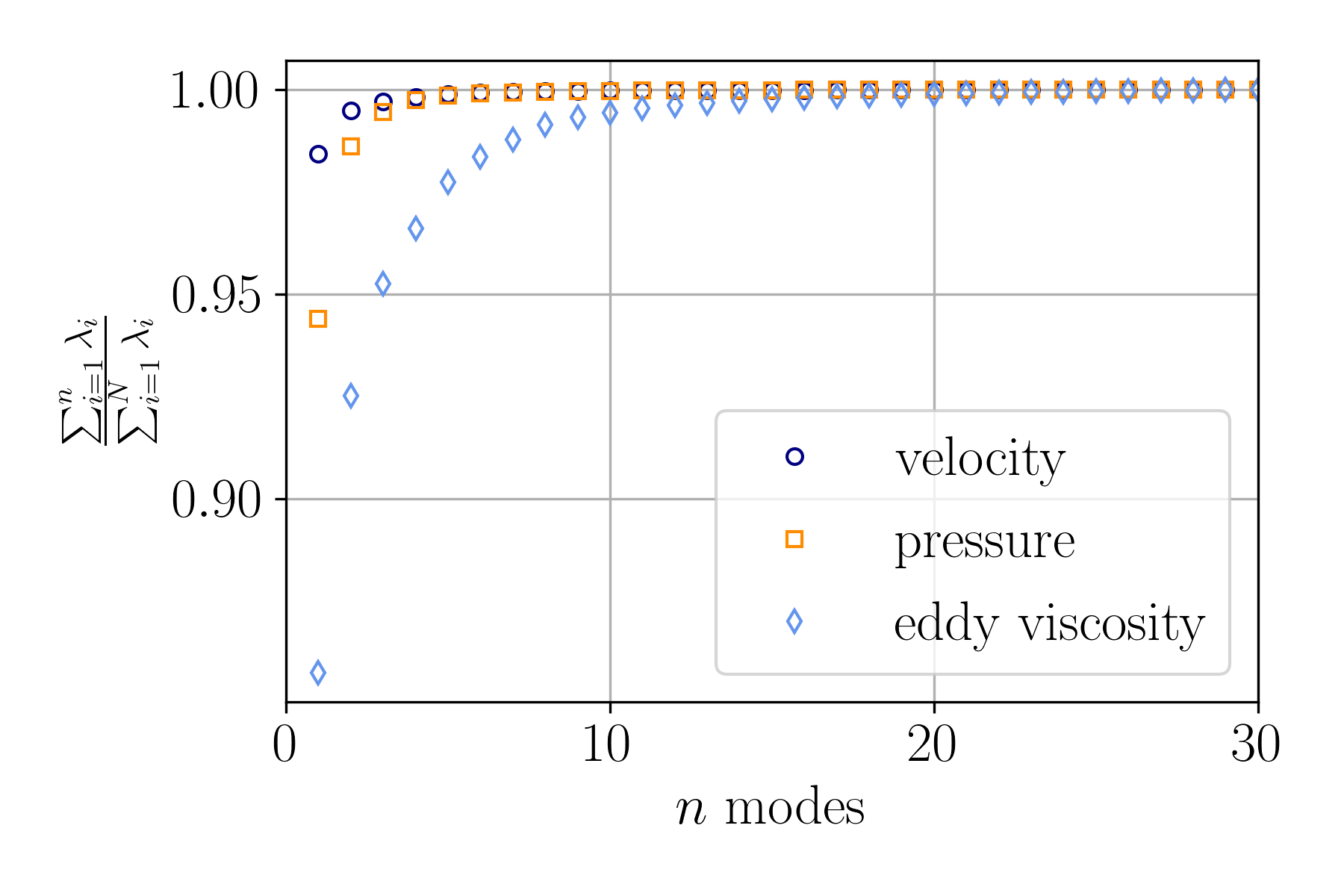}}
    \subfloat[POD eigenvalues decay]{\includegraphics[height=4.7cm]{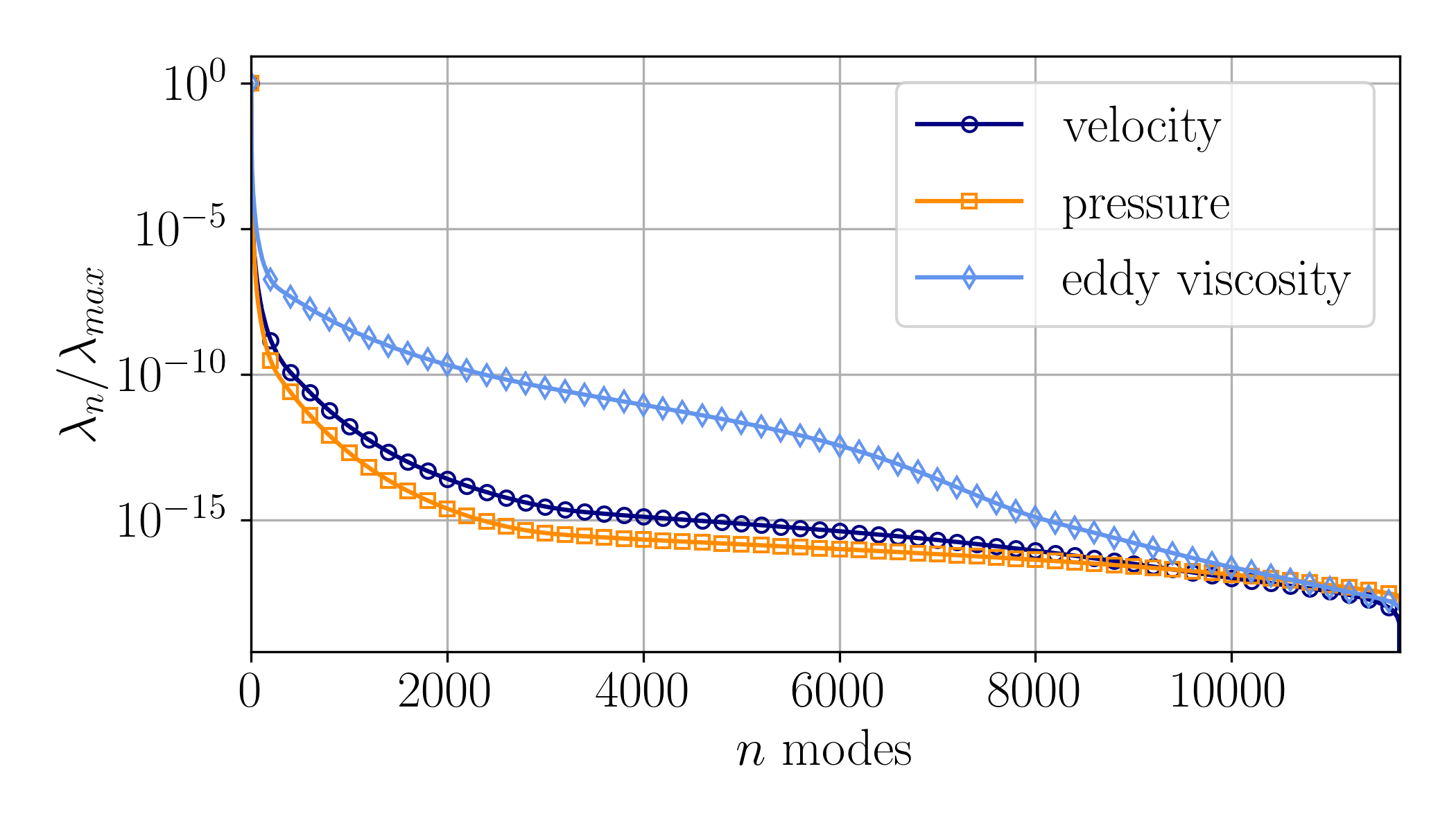}}
    \caption{Cumulative eigenvalues and eigenvalues decay for test case (\textbf{c}).}
    \label{fig:eig-c}
\end{figure}

As described in Subsection \ref{app:geometrical-params}, the POD basis is computed in the mid-configuration, corresponding to parameters $\bm{\mu}_{mid}=(\alpha_{mid}, h_1, h_2)=(25^{\circ}, 1 \text{ m}, 1.6\text{ m})$.
The POD cumulative eigenvalues' and the eigenvalues decay are displayed in Figure \ref{fig:eig-c}. This test case exhibits a significantly slower decay for the eddy viscosity field. We will analyze the ROM performance in the modes' combinations of Table \ref{tab:modal-regimes-c}. As for the other test cases, the criterion considered to select the modes $(N_u, N_p, N_{\nu_t})$ is the modes' retained energy. Indeed, as the Table shows, the number of $N_{\nu_t}$ is always bigger than $N_u$ and $N_p$.

\begin{table}[htpb!]
    \centering
 \caption{Combination of the modes for the velocity, pressure and eddy viscosity fields in test case (\textbf{c}).}
    \label{tab:modal-regimes-c}
    \begin{tabular}{>{\centering\arraybackslash}p{0.5\linewidth}
    >{\centering\arraybackslash}p{0.1\linewidth}
    >{\centering\arraybackslash}p{0.1\linewidth}
    >{\centering\arraybackslash}p{0.1\linewidth}
    }
    \toprule
    {\textbf{\emph{Minimum retained energy}} [$\%$]} &$N_u$& $N_p$ & $N_{\nu_t}$\\
    \midrule
    95&1&2&3\\
    99&2&3&8\\
    99.5&3&4&12\\
    99.6&3&4&14\\
    99.7&3&4&16\\
    99.8&4&5&20\\
    \bottomrule
    \end{tabular}
\end{table}

\subsubsection{Error analysis}
\label{subsubsec:errs-c}

For this test case, the penalty coefficient considered for the boundary conditions is $\tau=1$.

\begin{figure}[htpb!]
    \centering
    \includegraphics[width=0.9\textwidth]{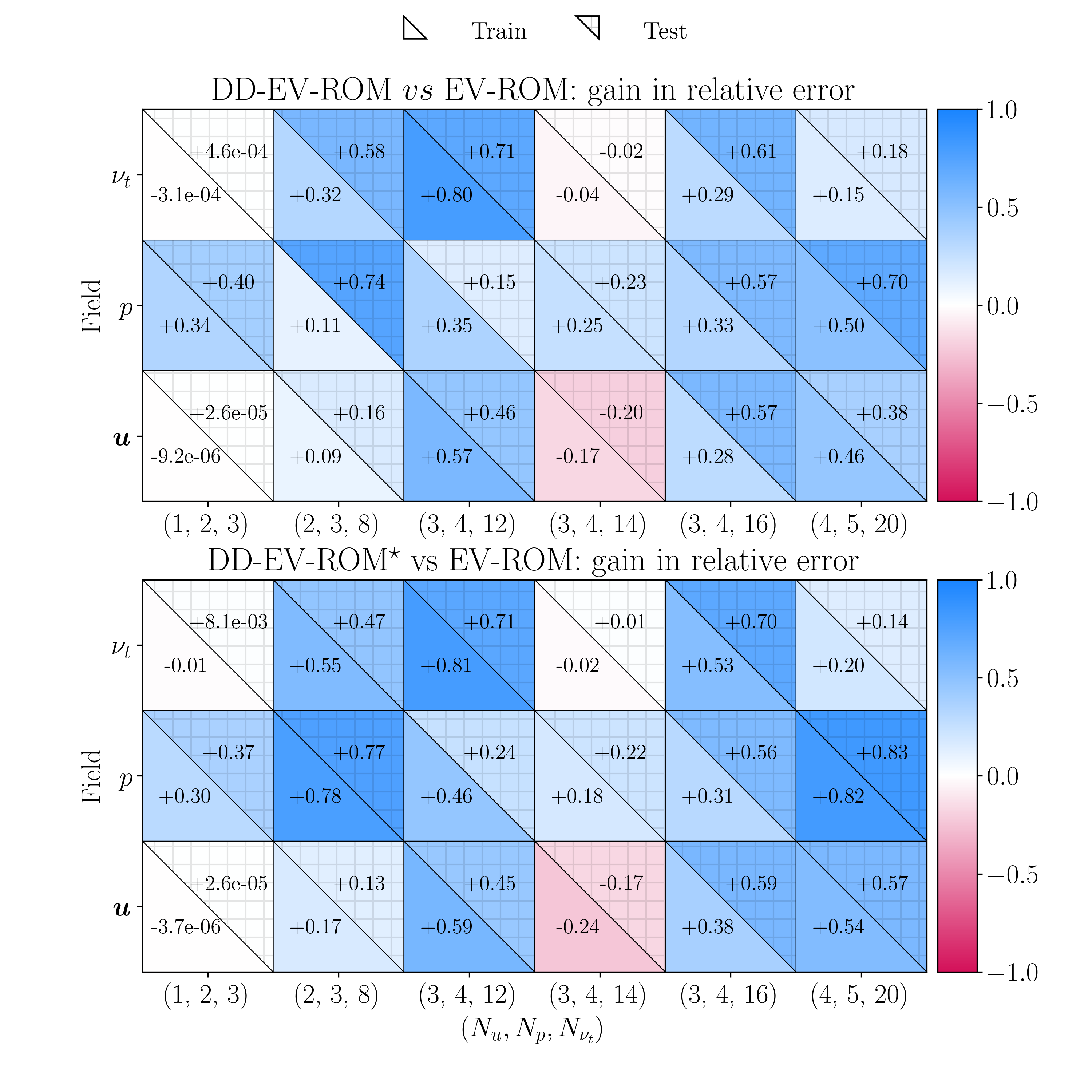}
    \caption{Comparison in the average performance of the DD-EV-ROM and DD-EV-ROM$^{\star}$ with respect to the baseline EV-ROM, for test case (\textbf{c}).
    The metric of interest is the average gain in relative error for train and test parameters (lower/upper triangles). The error is represented for different modes combination and for the three fields velocity, pressure and eddy viscosity.}
    \label{fig:backstep-heatmaps}
\end{figure}

Figure \ref{fig:backstep-heatmaps} represents the average gain in train and test configurations for the three fields of interest, in different regimes (also specified in Table \ref{tab:modal-regimes-c}).
We can draw the following considerations:
\begin{itemize}
    \item[$\bullet$] There is no clear trend in the average gain, differently from the first two test cases.
    For instance, in the first test case, we noticed that the average gain increases as the number of modes increases and we justified this because of the worse performance of the baseline EV-ROM.
    In the second test case, we noticed that only the pressure has a consistent gain, while the other fields have zero gain, since the EV-ROM was already close to the projection and cannot be improved more.
    Here, the fact that there is not a clear trend is linked to the \textbf{complexity of the snapshots' manifold} we are trying to approximate.
    Moreover, the POD basis is computed in a mid-configuration and this may lead to lower ROM accuracy for parameters far from the mid-configuration (even in the train set).
    \item[$\bullet$] The gain obtained for the eddy viscosity field is in general linked to the velocity one. For instance, when $(N_u, N_p, N_{\nu_t})=(1, 2, 3)$, the gain is close to zero for both velocity and eddy viscosity, while for $(N_u, N_p, N_{\nu_t})=(3, 4, 16)$ the two fields have similar gain values.
    This is due to the mapping $\mathcal{G}$, that has inputs $(\bm{a}, \bm{\mu})$. Therefore, $\bm{a}^{sol}\simeq\bm{a}^{proj}$ leads to $\mathcal{G}(\bm{a}^{sol}, \bm{\mu})\simeq\bm{g}^{proj}$.
    \item[$\bullet$] In general, we obtain positive or zero gain in all the cases considered, except for the velocity when $(N_u, N_p, N_{\nu_t})=(3, 4, 14)$. However, even in this case, this is compensated by a positive gain for the pressure and a zero gain for the eddy viscosity. Furthermore, we noticed that the negative gain for the velocity does not reflect in a negative gain for the eddy viscosity.
    \item[$\bullet$] The results of DD-EV-ROM and DD-EV-ROM$^{\star}$ are similar, with slightly better performances of the DD-EV-ROM$^{\star}$ method, as in $(N_u, N_p, N_{\nu_t})=(4,5, 20)$.
\end{itemize}

A statistical analysis can be found also for this test case in the Appendix, in \ref{app:case-c}.

\subsubsection{Graphical results}
\label{subsubsec:graph-c}
This part shows the qualitative results in terms of absolute error and of the fields of interest, for a test parameter, namely $(\alpha, h_1, h_2)=(6.3^{\circ}, 0.99\text{ m}, 1.58\text{ m})$, in the regime $(N_u, N_p, N_{\nu_t})=(3, 4, 12)$.

\begin{figure}[htpb!]
    \centering
    \subfloat[]{
    \includegraphics[width=0.5\linewidth, trim={3cm 5cm 3cm 5cm}, clip]{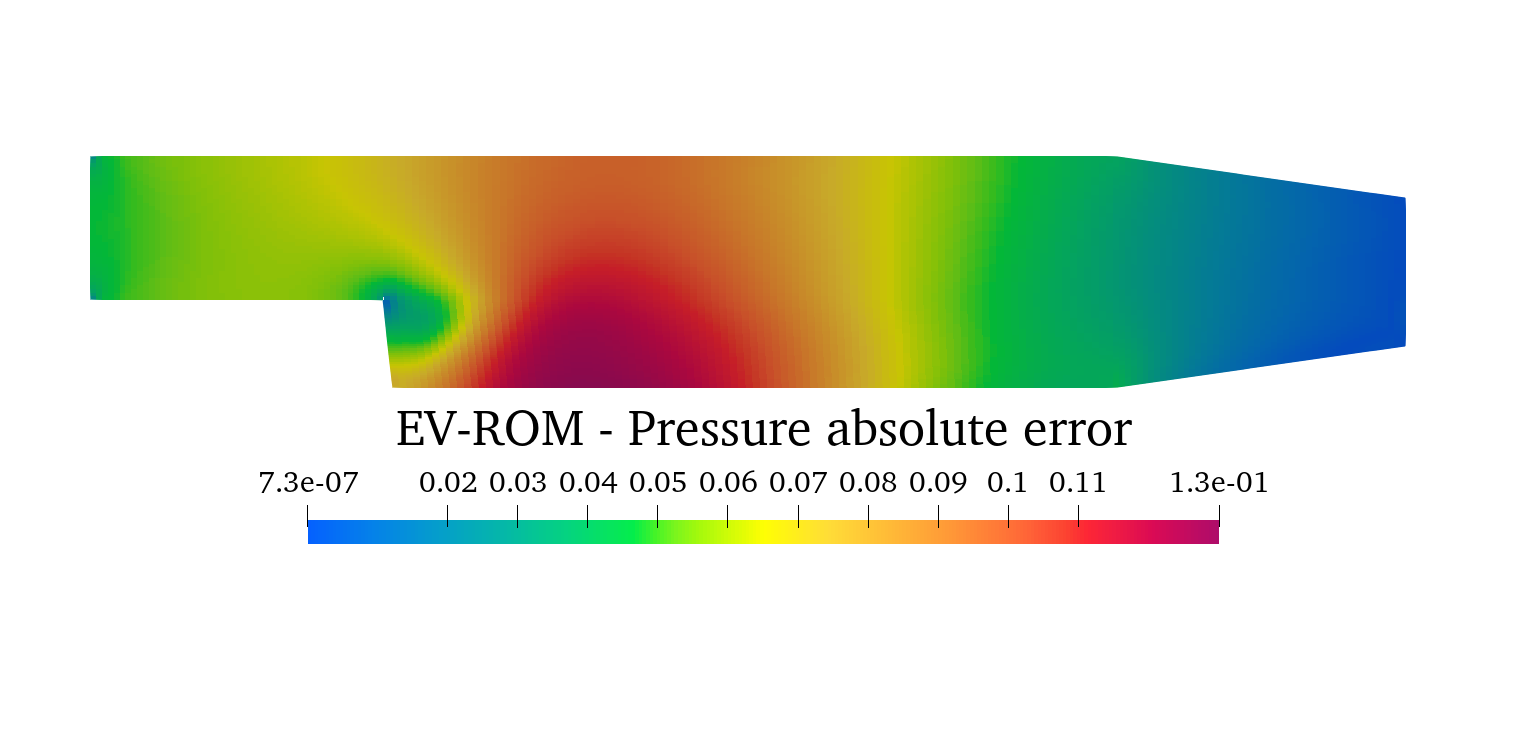}}
    \subfloat[]{
    \includegraphics[width=0.5\linewidth, trim={3cm 5cm 3cm 5cm}, clip]{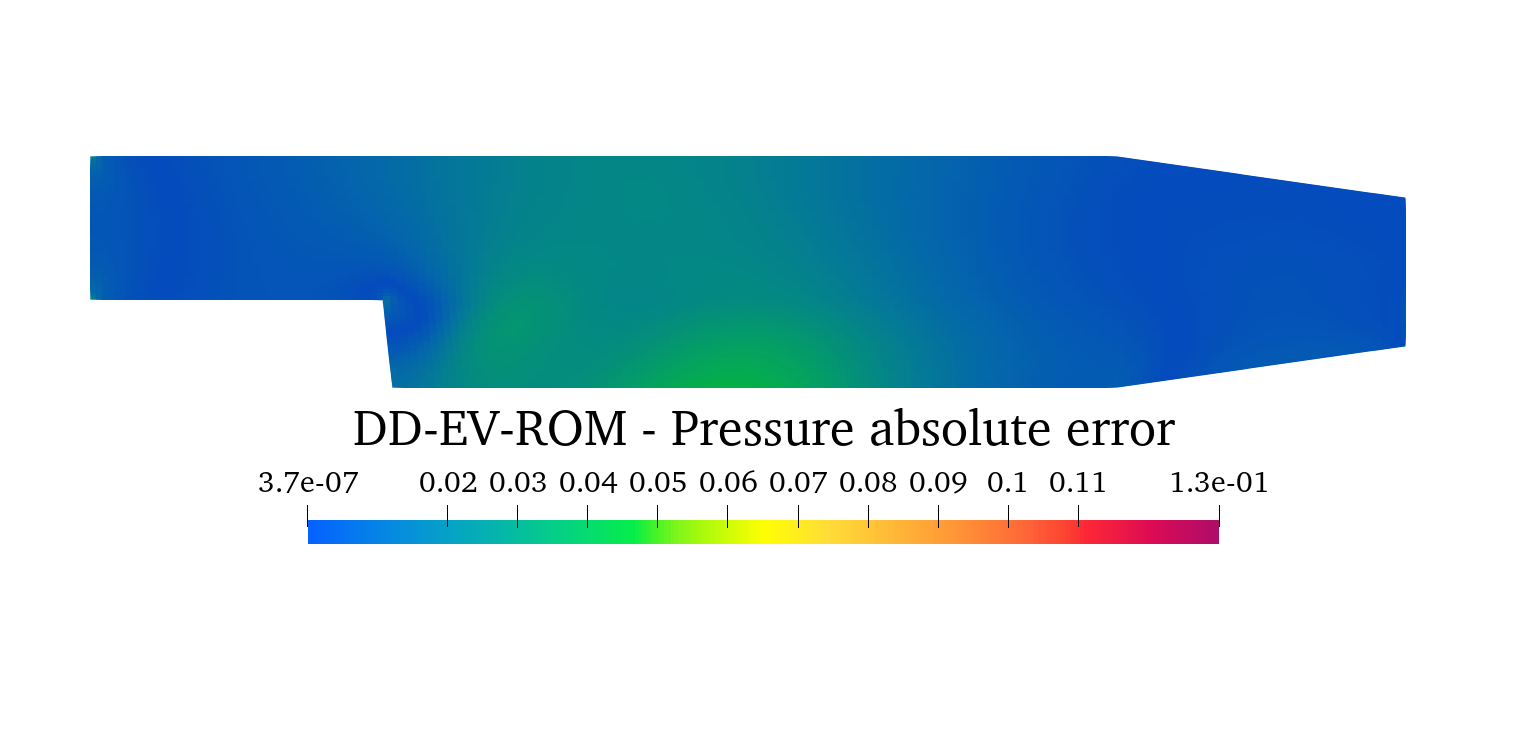}}\\
    \subfloat[]{
    \includegraphics[width=0.5\linewidth, trim={3cm 5cm 3cm 5cm}, clip]{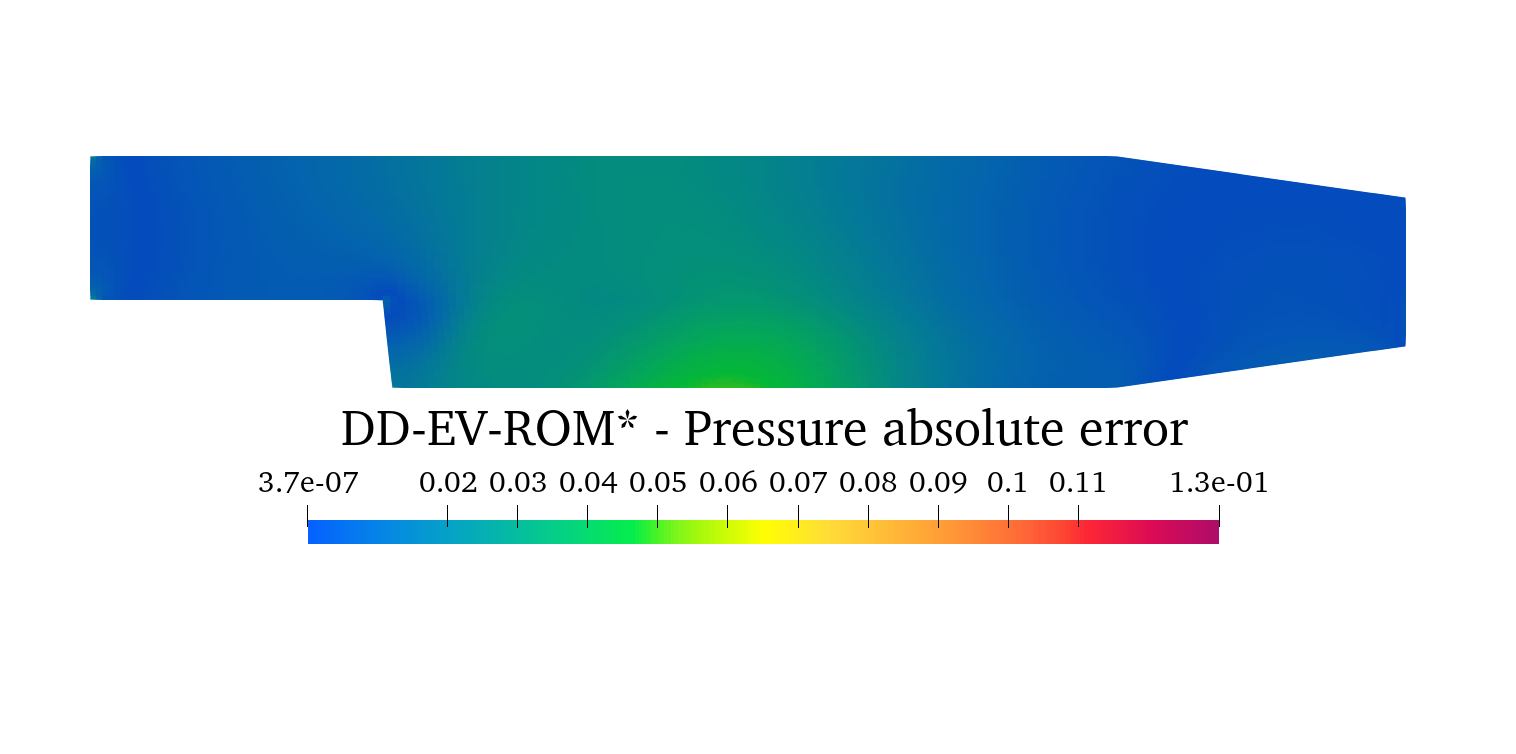}}
    \subfloat[]{
    \includegraphics[width=0.5\linewidth, trim={3cm 5cm 3cm 5cm}, clip]{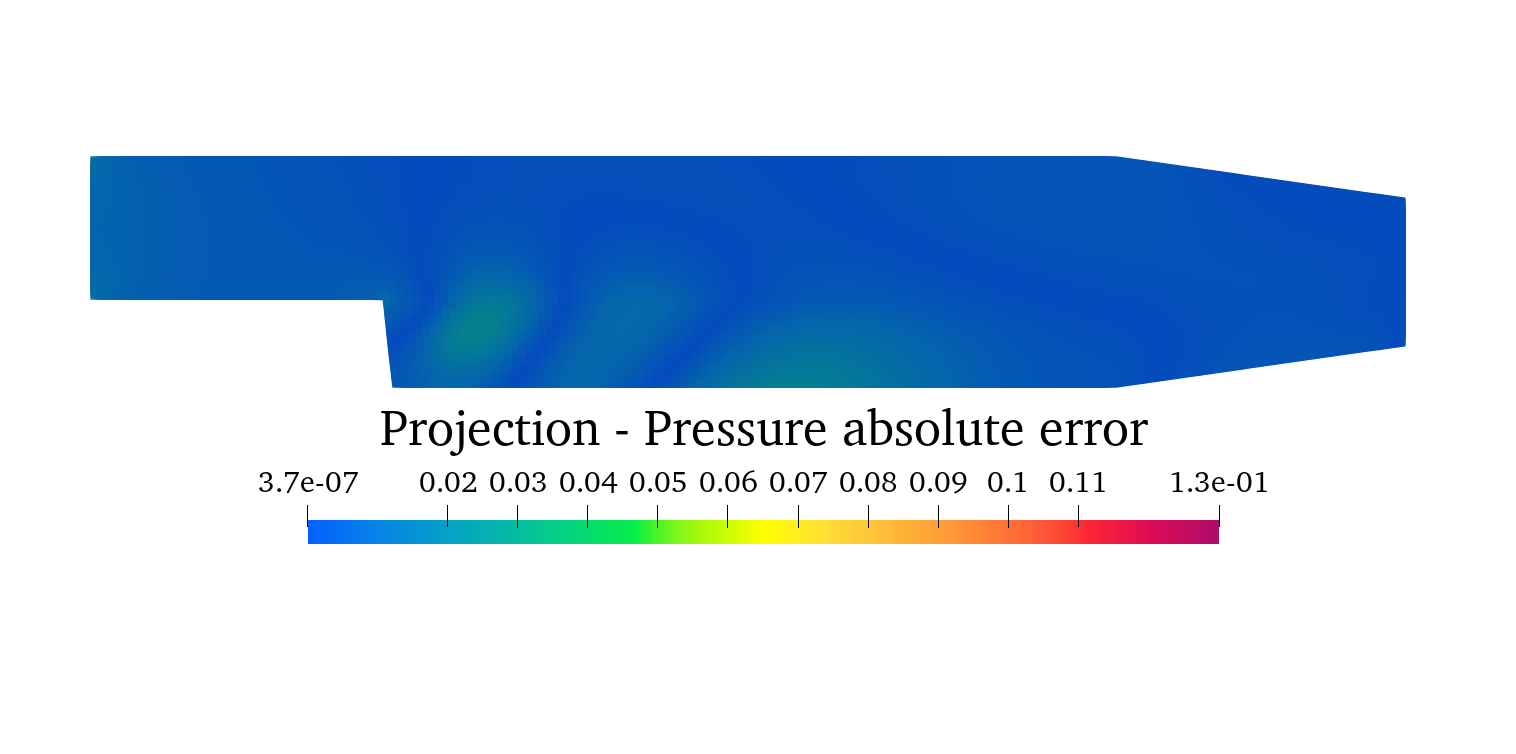}}
    \caption{Graphical pressure absolute error at the test parameters $(\alpha, h_1, h_2)=(6.3^{\circ}, 0.99 \text{ m}, 1.58\text{ m})$ of EV-ROM, DD-EV-ROM, DD-EV-ROM$^{\star}$, and for the projected field. The modal regime is $(N_u, N_p, N_{\nu_t})=(3, 4, 12)$.}
    \label{fig:backstep-graphical-p-err}
\end{figure}

Figure \ref{fig:backstep-graphical-p-err} displays the errors for the pressure field, and highlight the enhancement obtained with the DD-EV-ROM and DD-EV-ROM$^{\star}$, with respect to the baseline EV-ROM.
The velocity and eddy viscosity error fields are showed in \ref{app:case-c}. However, we stress the fact that the improvement is particularly evident especially for the pressure field, as for the other test cases.
The standard EV-ROM fails in accurately approximating this field, but the addition of the pressure closure $\bm{\tau}_p$ allows to almost reach the projection accuracy.

\subsection{Discussion}
\label{subsec:discussion}

This part of the manuscript is dedicated to a final comparison among the three test cases, and between the novel DD-EV-ROM and the baseline EV-ROM.

\begin{itemize}
    \item \textbf{POD performance}: Test cases (\textbf{b}) and (\textbf{c}) are more challenging than the first test case (\textbf{a}). Indeed, the POD eigenvalues decay is slower (Figures \ref{fig:eig-a}, \ref{fig:eig-b}, \ref{fig:eig-c}) and the reconstruction errors are much higher, as discussed in Subsections \ref{subsec:test-case-b} and \ref{subsec:test-case-c}.
    \item \textbf{EV-ROM}:
    The physics-based approach used in the EV-ROM leads to a good accuracy, especially when the number of modes is small. As already noticed, the reduced systems are in general more stable when characterized by a smaller number of degrees of freedom, and the Newton method converges faster. However, larger and ill-conditioned EV-ROM systems may lead to inaccurate results.
    \item \textbf{DD-EV-ROM} and \textbf{DD-EV-ROM}$^{\star}$:
    \begin{itemize}
        \item[$\circ$] As a general consideration, the data-driven closure approaches improve almost in all the modal regimes considered the baseline EV-ROM results, and also in a predictive setting and in unseen configurations. Moreover, the field that is most affected by the addition of the closure term is the pressure, especially for test cases (\textbf{a}) and (\textbf{b}). This is likely due to the addition of a dedicated pressure closure term into the PPE.
        \item[$\circ$] In the cylinder test case, the two training strategies for the neural operators act in a similar manner. From the time trends in Figures \ref{fig:cylinder-errs-1} and \ref{fig:cylinder-errs-2}, it can be noticed that the DD-EV-ROM$^{\star}$ is more stable in time, even in extrapolation settings.
        \item[$\circ$] In the cavity test case, the difference among the two training strategies is more evident. In specific modes' combinations, the alternative training procedure is necessary to enhance the ROM accuracy. Such strategy acts not only as a \emph{results-enhancer}, as for the cylinder, but also as a \emph{system stabilizer}.
        \item[$\circ$] In the backstep test case, the two approaches perform similarly. This is the most challenging test case because of the complexity of the snapshots' manifold, characterized by three geometrical parameters.
    \end{itemize}
    
\end{itemize}

\section{Conclusions and Outlook}
\label{sec:conclusions}
The project presented in this paper aims to enhance the classical ROM approaches using machine learning tools. The general paradigm of DD-ROMs was already presented in \cite{ivagnes2023hybrid}, but it is here extended to a more general parametric setup.

Section \ref{sec:intro} introduces the problem and the issues of standard ROMs in capturing the evolution of the system's dynamics.

Section \ref{sec:methods} is dedicated to the presentation of the methodologies used in this work. In particular, we introduce the FOM used to collect the snapshots in \ref{subsec:fom}, the POD-Galerkin ROM approach in presence of turbulence (EV-ROM) in \ref{subsec:ev-roms}, and the closure approaches in \ref{subsec:dd-ev-roms}. Such closure approaches have the goal of reintroducing the contribution of the neglected modes in the reduced system. Such strategies have already been investigated in other works \cite{xie2018data, ivagnes2023hybrid}, but not in a parameterized setting, which is the goal of this manuscript.
We provide a parametric extension through a machine learning approach based on neural operators.
Subsection \ref{subsec:networks} briefly recalls the logic of simple neural operators, namely the DeepONet, used to model the turbulence mapping $\mathcal{G}$, and the MIONet architecture, used to approximate the closure mapping $\mathcal{M}$. We introduce two alternative data-driven ROMs: the DD-EV-ROM, based on a standard training procedure, and the DD-EV-ROM$^{\star}$, based on an alternative training procedure. The second strategy is based on a coupled training of the two mappings $\mathcal{G}$ and $\mathcal{M}$.

Finally, in Section \ref{sec:results} we validate our method, showing the effects of the machine-learning approaches on the periodic flow past a cylinder (\ref{subsec:test-case-a}), on the channel-driven cavity flow (\ref{subsec:test-case-b}), and on the geometrically-parameterized backward-facing step flow (\ref{subsec:test-case-c}).

For each test case, we set up the offline stage in a parameterized setting, both with physical parameters (test cases (\textbf{a}) and (\textbf{b})) and with geometrical parameters (test case (\textbf{c})).

After that, a POD eigenvalues analysis is performed, after which we analyze the ROM performance in different modal regimes. In this framework, the baseline EV-ROMs may produce inaccurate results and poor reconstructions of the solutions.

The DD-EV-ROM and DD-EV-ROM$^{\star}$ systems provide successful results : (i) average gain of the DD-EV-ROMs with respect to the EV-ROM in train and test configurations; (ii) velocity and pressure relative errors with respect to the FOM counterpart in train and test settings (including the time trend for the unsteady cases (\textbf{a}) and (\textbf{b})); (iii) graphical and qualitative analysis of the fields and of the absolute errors with respect to the high-fidelity reference.

In conclusion, this work provides a strategy able to enhance the baseline EV-ROM results which consists in adding data-driven \emph{closure} terms into the reduced system. This machine-learning based addition not only improve the accuracy of the standard EV-ROM, but also acts as a stabilizer in highly unsteady setups and ill-conditioned reduced-systems. Moreover, the DD-EV-ROM$^{\star}$ method leads to more stable results, also in time extrapolation settings.

As pointed out for the geometrically-parameterized test case, the results' accuracy is influenced by the linearity of the POD approach. Hence, the method may be further be improved introducing \emph{nonlinearity} in the projection step, by replacing the POD with more advanced techniques, like autoencoders. This will be the focus of the authors' future research work.

\section*{Acknowledgements}
\label{sec:acknoledgements}
We acknowledge the support by the European Commission H2020 ARIA (Accurate ROMs for Industrial Applications, GA 872442) project, by MIUR (Italian Ministry for Education University
and Research) and by the European Research Council Consolidator Grant Advanced Reduced Order Methods with Applications in Computational Fluid Dynamics-GA 681447, H2020-ERC COG 2015 AROMA-CFD. 
This work has been conducted within the research activities of the consortium iNEST (Interconnected North-East Innovation Ecosystem), Piano Nazionale di Ripresa e Resilienza (PNRR) – Missione 4 Componente 2, Investimento 1.5 – D.D. 1058 23/06/2022, ECS00000043, supported by the European Union's NextGenerationEU program.
We also acknowledge the support by INdAM-GNCS: Istituto Nazionale di Alta Matematica –– Gruppo Nazionale di Calcolo Scientifico.

Giovanni Stabile acknowledges the financial support under the National Recovery and Resilience Plan (NRRP), Mission 4, Component 2, Investment 1.1, Call for tender No. 1409 published on 14.9.2022 by the Italian Ministry of University and Research (MUR), funded by the European Union – NextGenerationEU– Project Title ROMEU – CUP P2022FEZS3 - Grant Assignment Decree No. 1379 adopted on 01/09/2023 by the Italian Ministry of Ministry of University and Research (MUR).

The main computations in this work were carried out by the usage of ITHACA-FV \cite{ithacasite}, an open-source library and an implementation in OpenFOAM \cite{ofsite} for reduced order modeling techniques. Its developers and contributors are acknowledged.

\newpage

\bibliographystyle{abbrv}
\bibliography{main}

\appendix

\appendix

\section{Additional methodology details}
\label{appA}

\subsection{On the training procedure}
\label{subsec:app-training}
This paragraph is dedicated to highlight the potential of the DD-EV-ROM$^{\star}$ (described in Section \ref{subsec:networks}) and the reason why the standard DD-EV-ROM may be inaccurate.

Indeed, considering only the contribution in \eqref{eq:loss-mionet-stand} may cause lack of accuracy in the approximation of the closure coefficients, especially in test and extrapolation settings.
The reason is that network $\mathcal{M}$ is trained only relying on the \emph{exact} coefficient $\bm{a}^{proj}$ and $\bm{g}^{proj}$. However, it may happen that, due to bad conditioning, when solving the online system \eqref{eq:ppe-rom-turb}, the velocity solution $\bm{a}^{sol}$ differs from the projected coefficients $\bm{a}^{proj}$, then $$\|\bm{a}^{sol}-\bm{a}^{proj}\|_2 \neq 0.$$
Since network $\mathcal{G}$ is also trained relying on the exact coefficients $\bm{a}^{proj}$, we might have:
$$\|\mathcal{G}(\bm{a}^{sol}, \bm{\mu})-\bm{g}^{proj}\|_2 \neq 0,$$
and, hence, equivalently for network $\mathcal{M}$: 
\begin{equation}
    \begin{split}
        \|&\mathcal{M}\bigl(\bm{a}^{sol}, \mathcal{G}(\bm{a}^{sol}, \bm{\mu}), \bm{\mu}\bigr)-\bm{\tau}^{exact}(\bm{a}^{proj}, \bm{g}^{proj}, \bm{\mu})\|_2=\\
        \| &\Bigl[ \mathcal{M}\bigl(\bm{a}^{sol}, \mathcal{G}(\bm{a}^{sol}, \bm{\mu}), \bm{\mu}\bigr)- \bm{\tau}^{exact}\bigl(\bm{a}^{proj}, \mathcal{G}(\bm{a}^{sol}, \bm{\mu}), \bm{\mu}\bigr) \Bigr] 
        + \\
        &\Bigl[ \bm{\tau}^{exact}\bigl(\bm{a}^{proj}, \mathcal{G}(\bm{a}^{sol}, \bm{\mu}), \bm{\mu}\bigr) - \bm{\tau}^{exact}\bigl(\bm{a}^{proj}, \bm{g}^{proj}, \bm{\mu}\bigr) \Bigr]  \|_2
\gg0.
    \end{split}
    \label{eq:error-standard-mionet}
\end{equation}

The expression in Equation \eqref{eq:error-standard-mionet} contains two distinct contributions:
\begin{itemize}
    \item The first term $\mathcal{M}\bigl(\bm{a}^{sol}, \mathcal{G}(\bm{a}^{sol}, \bm{\mu}), \bm{\mu}\bigr)- \bm{\tau}^{exact}\bigl(\bm{a}^{proj}, \mathcal{G}(\bm{a}^{sol}, \bm{\mu}), \bm{\mu}\bigr)$ is only related with the generalization capability of network $\mathcal{M}$.
    \item The second term $\bm{\tau}^{exact}\bigl(\bm{a}^{proj}, \mathcal{G}(\bm{a}^{sol}, \bm{\mu}), \bm{\mu}\bigr) - \bm{\tau}^{exact}\bigl(\bm{a}^{proj}, \bm{g}^{proj}, \bm{\mu}\bigr)$ is related with the performance of network $\mathcal{G}$. It is also important to highlight that the \emph{exact} closure term (the MIONet's target) depends on the eddy viscosity variables $\bm{g}$, as in expression \eqref{eq:exact-nonlinear-correction}.
\end{itemize}

The difference between the two types of training can be seen in Subsection \ref{subsec:test-case-b}, where even small perturbations in the velocity solution lead to inaccurate results for mapping $\mathcal{M}$, as in \eqref{eq:error-standard-mionet}.
Moreover, in the DD-EV-ROM this may be reflected in inaccurate pressure approximations, as the reduced pressure $\bm{b}^{sol}$ comes after the residual minimization of the data-driven enriched reduced system, as in Algorithms \ref{alg-steady} and \ref{alg-unsteady}.

\subsection{On the geometrical parameterization}
\label{app:geometrical-params}

When a geometrical parameterization is employed, as in the third test case of Subsection \ref{subsec:test-case-c}, the snapshots belong to different domains and grids, that share the same number of degrees of freedom $N_{dof}$. In this case, we build the reduced order model considering a procedure similar to  that employed in \cite{zancanaro2021hybrid}.

In particular, we consider a reference grid, corresponding to a mid-configuration for each parameter:
$$\bm{\mu}_{mid}=\dfrac{1}{N_{\mu}}\sum_{i=1}^{N_{\mu}} \bm{\mu}_i.$$

Then, we collect the snapshots matrices for each field of interest $\bm{s}$, and we assemble the correlation matrix as:
$$K_{ij}=\bm{s}_i^T M_{mid}\bm{s}_j,$$
where $M_{mid}$ is the mass matrix defined for $\Omega (\bm{\mu}_{mid})$ \cite{zancanaro2021hybrid}.

After that, a unique set of basis functions is found $\mathbb{V}^s_{\text{POD}}$.

However, the key difference in the geometrical parameterization is that all the operators appearing in System \eqref{eq:ppe-rom-turb} are mesh-dependent, and hence, parameter-dependent, namely:
 \[
 \begin{split}
    &\bm{B}(\bm{\mu}^{\star}), \bm{B_T}(\bm{\mu}^{\star}), \bm{C}(\bm{\mu}^{\star}), \bm{H}(\bm{\mu}^{\star}), \bm{D}(\bm{\mu}^{\star}), \bm{G}(\bm{\mu}^{\star}), \bm{N}(\bm{\mu}^{\star}), \bm{C_{T1}}(\bm{\mu}^{\star}), \bm{C_{T2}}(\bm{\mu}^{\star}),\bm{C_{T3}}(\bm{\mu}^{\star}), \bm{C_{T4}}(\bm{\mu}^{\star}), \\&\bm{D}^k(\bm{\mu}^{\star}), \bm{E}^k(\bm{\mu}^{\star}) \text{ with }k=1, \dots, N_{\text{BC}}. 
 \end{split}
 \]
There exist different alternatives in literature to deal with this type of test case. One example is the Discrete Empirical Interpolation Method (DEIM). However, this goes beyond the scope of this project and we choose here to build the reduced operators for each parameter, also in the test setting, as highlighted in the gray part of Algorithm \ref{alg-steady}.
This choice does not affect consistently the computational time and cost of the online stage, since the reduced operators that we need to assemble have a small size, since $N_u, N_p, N_{\nu_t}\ll N_{dof}$.

\subsection{Hyperparameters of the turbulence and closure mappings}
\label{app:hyperparams}

Table \ref{tab:networks_operators}
displays the details of the architectures of the DeepONet used for the turbulence model $\mathcal{G}$ and of the MIONet used for the correction model $\mathcal{M}$.

Additionally, in both models, the learning rate decays by a factor $\gamma=0.2$ every $N_{\text{step}}=3000$ epochs. In the DD-EV-ROM$^{\star}$, after a pre-training (following the hyperparameters in Table \ref{tab:networks_operators}), both models are trained together for additionally $20000$ epochs, minimizing the loss function expressed in $\mathcal{L}_{\mathcal{M}}^{\star}$.

\begin{table}[htpb!]
    \centering
 \caption{Hyperparameters of DeepONet $\mathcal{G}$ and MIONet $\mathcal{M}$ networks, considered for the turbulence and closure mappings, respectively.}
    \label{tab:networks_operators}
    \begin{tabular}{>{\centering\arraybackslash}p{0.18\linewidth}
    >{\centering\arraybackslash}p{0.15\linewidth}
    >{\centering\arraybackslash}p{0.13\linewidth}
    >{\centering\arraybackslash}p{0.15\linewidth}
    >{\centering\arraybackslash}p{0.2\linewidth}
    >{\centering\arraybackslash}p{0.1\linewidth}
    }
    \toprule
    {\textbf{\emph{Mapping}}} &{\textbf{\emph{Network}}}& \textbf{\emph{Hidden layers}} & \textbf{\emph{Non-linearity}} &\textbf{\emph{Learning rate}} & \textbf{\emph{Stopping epoch}}\\
    \midrule
\multirow{3}{*}{$\bm{g}=\mathcal{G}(\bm{a}, \bm{\mu})$}& Branch net $\mathcal{B}$ &{$[20,20,20]$}&{Softplus}& \multirow{3}{*}{$(\num{1e-3})\gamma^{\lfloor \frac{t_{\text{epoch}}}{N_{\text{step}}} \rfloor}$} & \multirow{3}{*}{20000}\\ 
\cmidrule{2-4}
&Trunk net $\mathcal{T}$& $[20, 20, 20]$& Softplus&&\\
\cmidrule{2-4}
&Reduction net $\mathcal{R}$& $[20, 20, 20]$& Softplus&&\\
\midrule \multirow{3}{*}{$\bm{\tau}=\mathcal{M}(\bm{a}, \bm{g}, \bm{\mu})$}&Branch net $\mathcal{B}_1$&$[20, 20, 20]$& Softplus&\multirow{3}{*}{$(\num{1e-3})\gamma^{\lfloor \frac{t_{\text{epoch}}}{N_{\text{step}}} \rfloor}$} & \multirow{3}{*}{20000}\\
\cmidrule{2-4}
    &Branch net $\mathcal{B}_2$&$[20, 20, 20]$& Softplus&&\\
   \cmidrule{2-4}
   &Trunk net $\mathcal{T}$&$[20, 20, 20]$& Softplus&&
   \\
   \cmidrule{2-4}
   &Reduction net $\mathcal{R}$&$[20, 20, 20]$& Softplus&&\\
    \bottomrule
    \end{tabular}
\end{table}

\section{Additional numerical results}
\label{app:results}
This part is dedicated to show additional numerical results for the test cases presented in Sections \ref{sec:results}.

\subsection{Test case (\textbf{a})}
\label{app:case-a}
\begin{figure}[htpb!]
    \centering
    \includegraphics[width=\linewidth]{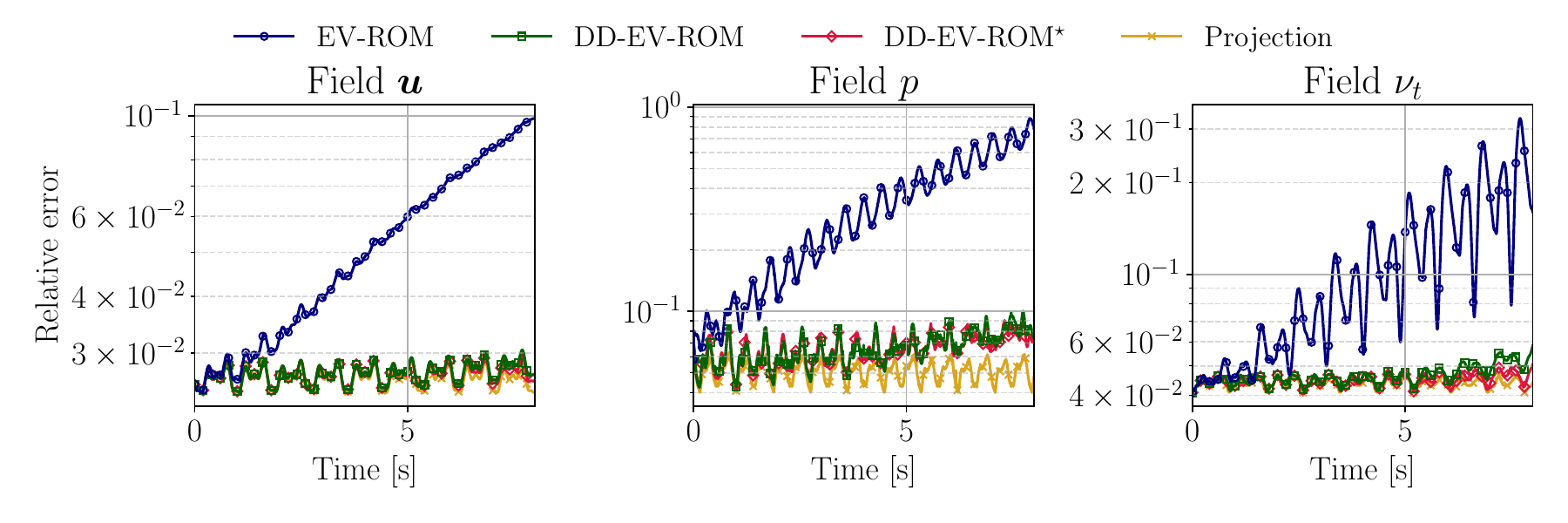}
    \caption{Time trends of the relative $L^2$ errors for the three fields of interest $\bm{u}$, $p$, and $\nu_t$, for one test viscosity $\nu=1.15e-4 \, \frac{m^2}{s}$. The modal regime is $(N_u, N_p, N_{\nu_t})=(3, 6, 5)$.}
    \label{fig:cylinder-errs-1}
\end{figure}

\begin{figure}[htpb!]
    \centering
    \includegraphics[width=\linewidth]{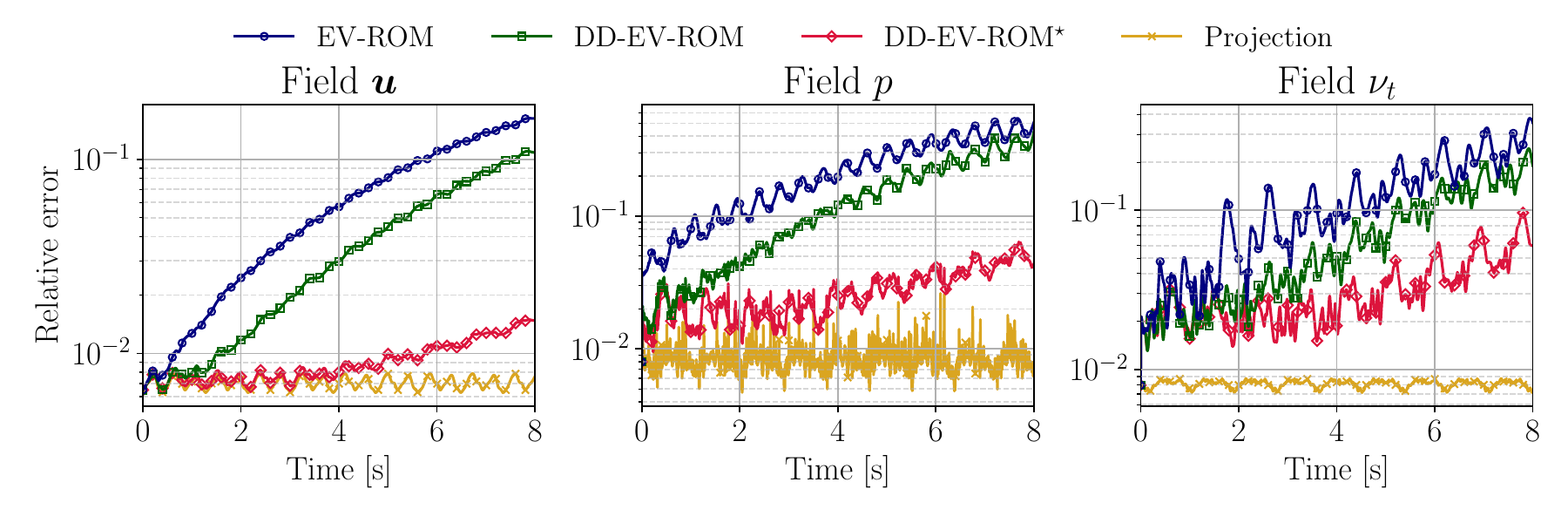}
    \caption{Time trends of the relative $L^2$ errors for the three fields of interest $\bm{u}$, $p$, and $\nu_t$, for one test viscosity $\nu=1.15e-4 \, \frac{m^2}{s}$. The modal regime is $(N_u, N_p, N_{\nu_t})=(10, 24, 30)$.}
    \label{fig:cylinder-errs-2}
\end{figure}

For the first test case, we focus here on how the error varies varying the number of modes. In particular, Figures \ref{fig:cylinder-errs-1} and \ref{fig:cylinder-errs-2} represent the relative error trend in time, for two modal regimes $(N_u, N_p, N_{\nu_t})=(3, 6, 5)$ and $(N_u, N_p, N_{\nu_t})=(10, 24, 30)$, and in unseen settings.

    \begin{itemize}
        \item[$\circ$] On the one hand, POD-Galerkin ROMs become more unstable when considering a large number of modes, and it may happen that the EV-ROM residuals blow up and/or provide inaccurate approximations. It can be immediately noticed from Figure \ref{fig:cylinder-errs-2} that the EV-ROM (the blue line in the plots) has a large error, increasing as time evolves, reaching the same order of magnitude as in the case $(N_u, N_p, N_{\nu_t})=(3, 6, 5)$ (Figure \ref{fig:cylinder-errs-1} in the Appendix).
        When this happens, the addition of the closure $\bm{\tau}$ acts as a \emph{stabilizer} for the reduced system, leading to a higher gain in the relative error.
        \item[$\circ$] On the other hand, the reconstruction error (yellow line in the plots) decreases as the number of modes increases. The closure terms aim to model the contribution that allows to reach such reconstruction error, that is the best result we can achieve with DD-EV-ROMs.
    \end{itemize}
     Hence, increasing the number of modes: (i) the EV-ROM shows the same performance and its error is large, (ii) the reconstruction error decreases, (iii) the DD-EV-ROMs error should be as close as possible to the reconstruction error, and it decreases as well. This leads to the larger gains in the error observed in Figure \ref{fig:cylinder-heatmaps}.

\begin{figure}[htpb!]
    \centering{\includegraphics[width=\textwidth]{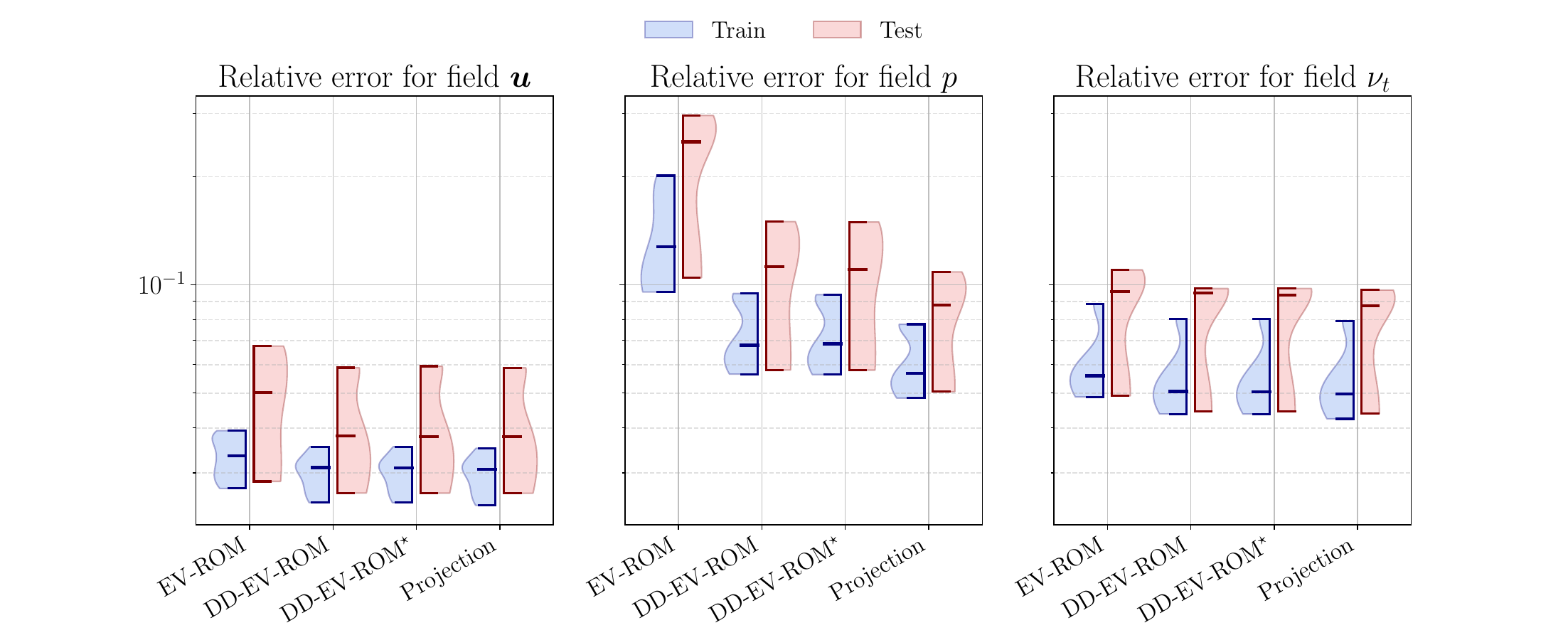}}
    \caption{Statistical performance for train and test viscosities in the time average of the relative $L^2$ error with respect to the high-fidelity solutions. The Figure represents the median and the error bounds for the three fields of interest in the case $(N_u, N_p, N_{\nu_t})=(3, 6, 5)$.}
    \label{fig:cylinder-violin-1}
\end{figure}

\begin{figure}[htpb!]
    \centering{\includegraphics[width=\textwidth]{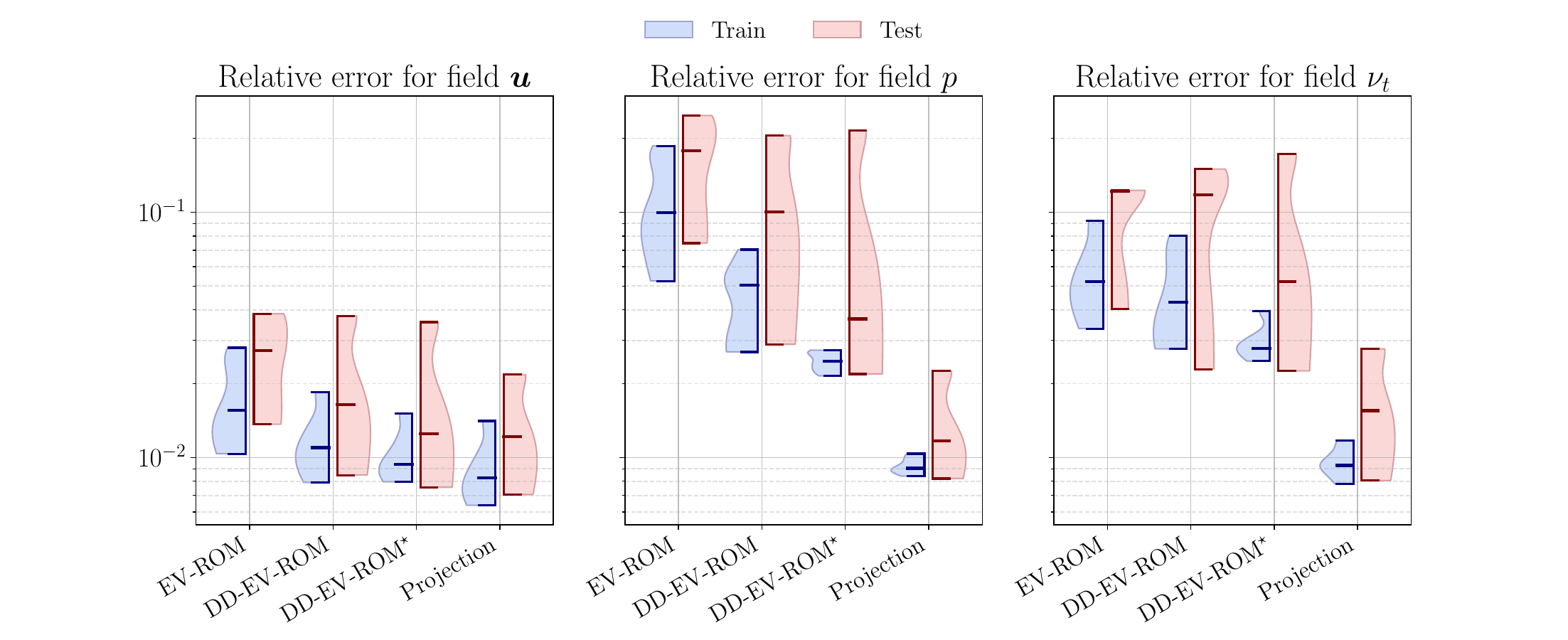}}
    \caption{Statistical performance for train and test viscosities in the time average of the relative $L^2$ error with respect to the high-fidelity solutions. The Figure represents the median and the error bounds for the three fields of interest in the case $(N_u, N_p, N_{\nu_t})=(10, 24, 30)$.}
    \label{fig:cylinder-violin-2}
\end{figure}

Similar considerations on the comparison of different regimes arise from Figures \ref{fig:cylinder-violin-1} and \ref{fig:cylinder-violin-2}, which show the statistical performance of the methods in different train and test viscosities, and for the two modes' combinations, respectively. It can be immediately seen that the ranges of the EV-ROM relative errors for the pressure are only slightly different in the two modal regimes, as described above.
Moreover, from plots in Figures \ref{fig:cylinder-violin-1} and \ref{fig:cylinder-violin-2}, we can draw the following conclusions:
\begin{itemize}
    \item In all cases a large difference between train and test errors for all the methods. Since this happens also for the projection error, it is related with the POD space itself and with the chosen sampling, but not with the closure strategies employed.
    \item Both DD-EV-ROM and DD-EV-ROM$^{\star}$ outperform the baseline EV-ROM when $(N_u, N_p, N_{\nu_t})=(3, 6, 5)$, and show statistics close to the projection. Moreover, their performances are similar in this modal regime.
    \item The DD-EV-ROM$^{\star}$ outperforms the DD-EV-ROM method when considering a larger number of modes $(10, 24, 30)$. As above-mentioned, ROM systems may be more unstable when the number of unknowns is larger, and the alternative training procedure described in Section \ref{subsec:networks} enhances the stabilization of the method and the accuracy of the predictions. This behavior is even more evident in Figure \ref{fig:cylinder-errs-2}, where the DD-EV-ROM$^{\star}$ provides more accurate approximations also in time extrapolation settings, when the time interval is \textbf{four times larger} than the training one.
\end{itemize}


\begin{figure}[htpb!]
    \centering
    \subfloat[]{
    \includegraphics[width=0.5\linewidth, trim={3cm 1cm 3cm 3cm}, clip]{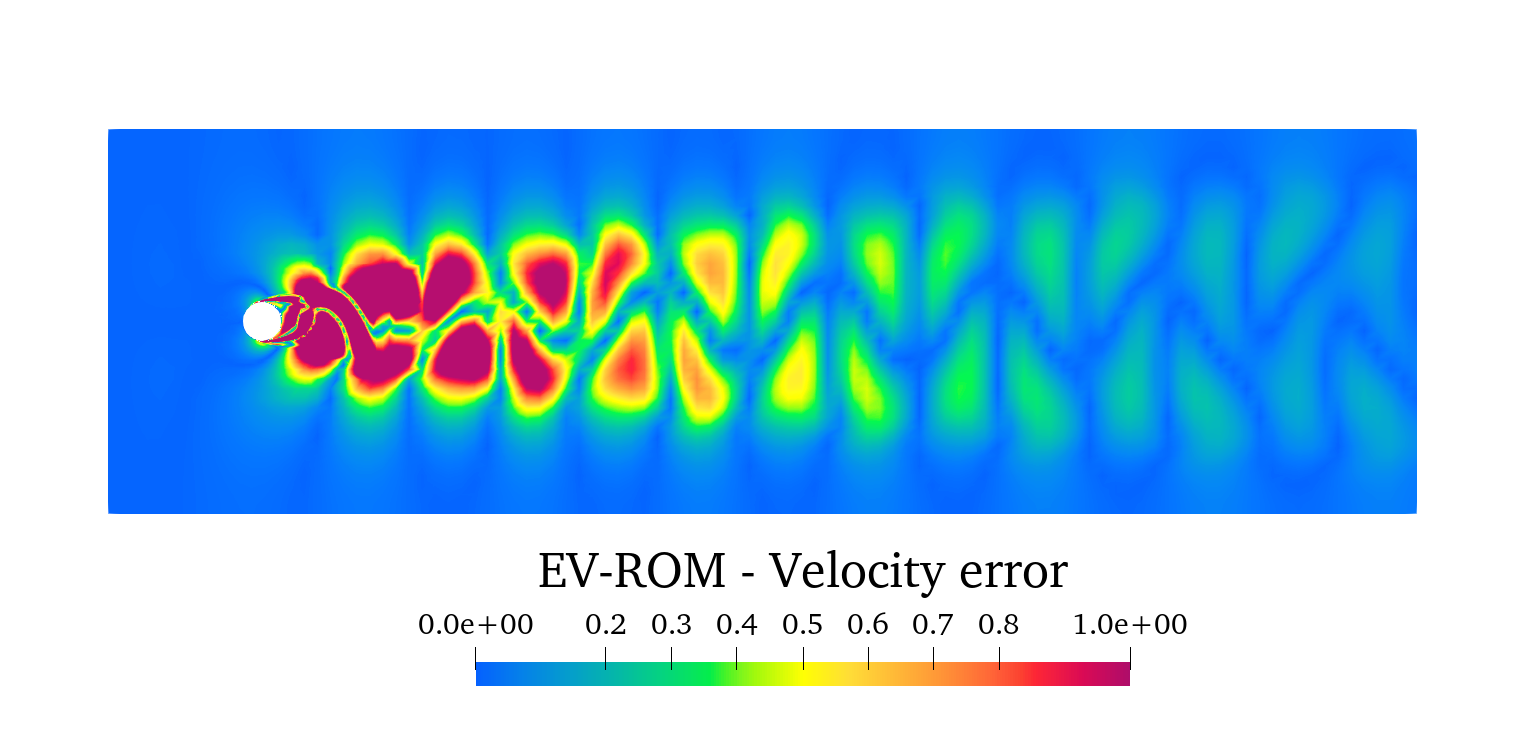}}
    \subfloat[]{
    \includegraphics[width=0.5\linewidth, trim={3cm 1cm 3cm 3cm}, clip]{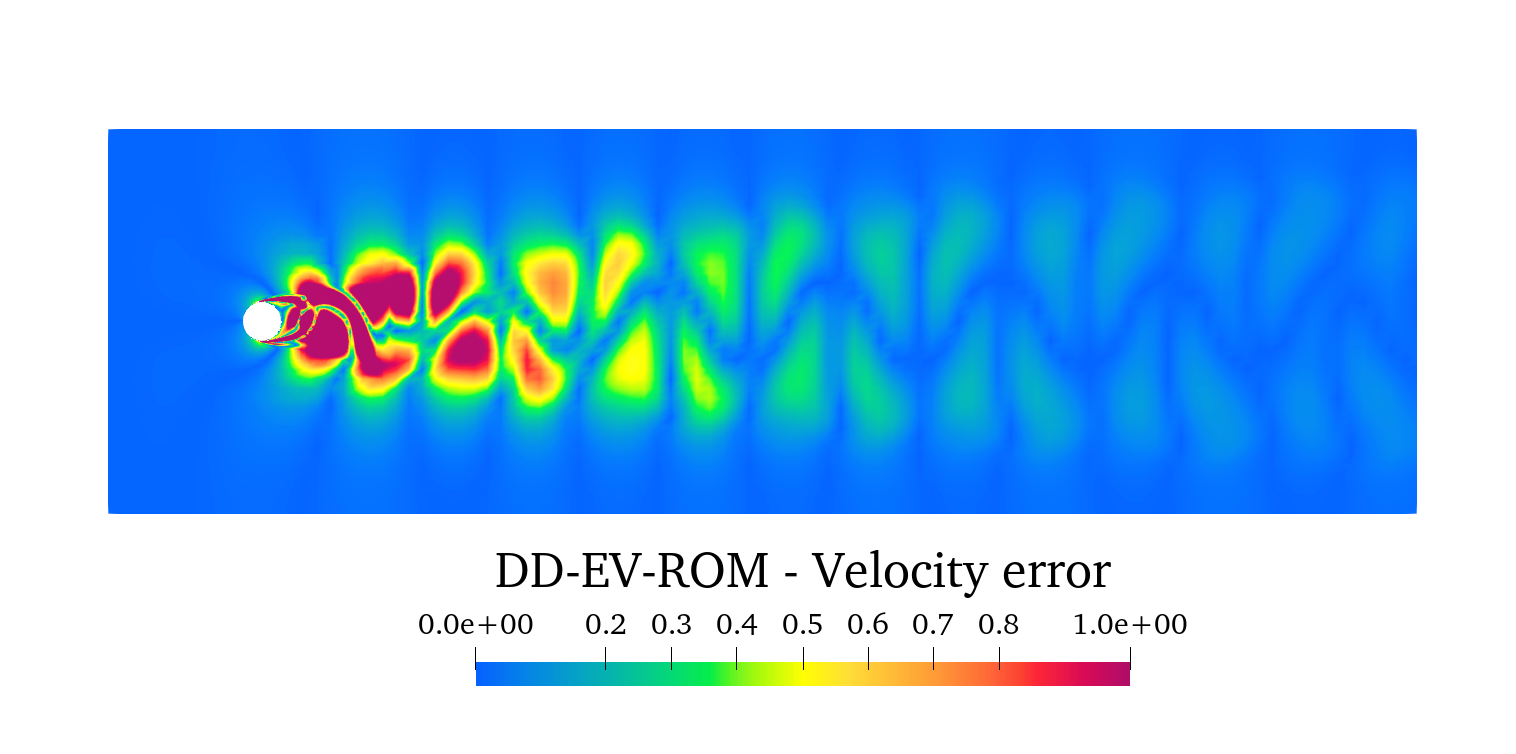}}\\
    \subfloat[]{
    \includegraphics[width=0.5\linewidth, trim={3cm 1cm 3cm 3cm}, clip]{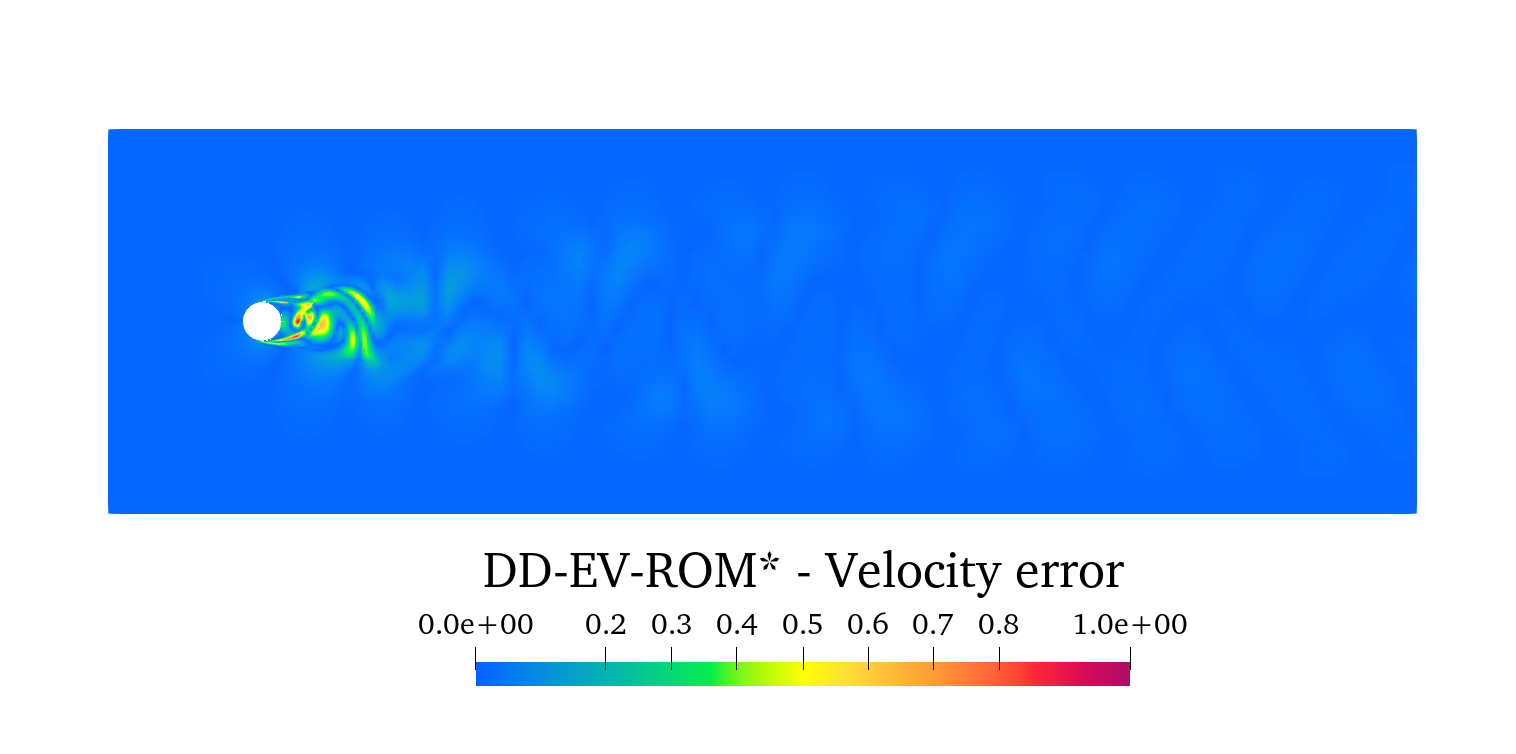}}
    \subfloat[]{
    \includegraphics[width=0.5\linewidth, trim={3cm 1cm 3cm 3cm}, clip]{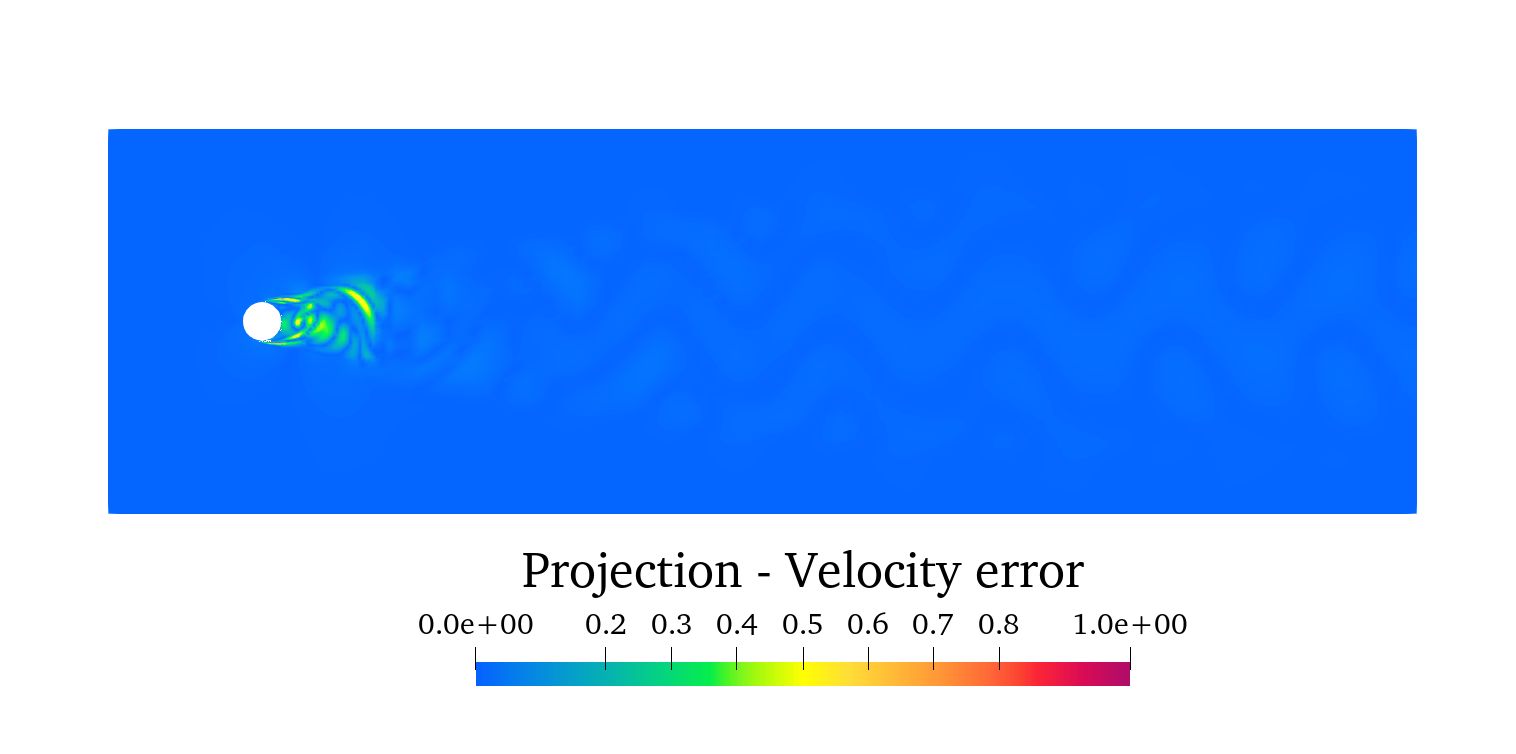}}
    \caption{Graphical absolute error for the velocity magnitude at $\nu=1.15e-4 \, \frac{m^2}{s}$ at the final time instance ($t=8$), of EV-ROM, DD-EV-ROM, DD-EV-ROM$^{\star}$, and for the projected field. The modal regime is $(N_u, N_p, N_{\nu_t})=(10, 24, 30)$.}
    \label{fig:cyl-graphical-u-err}
\end{figure}

\begin{figure}[htpb!]
    \centering
    \subfloat[]{
    \includegraphics[width=0.5\linewidth]{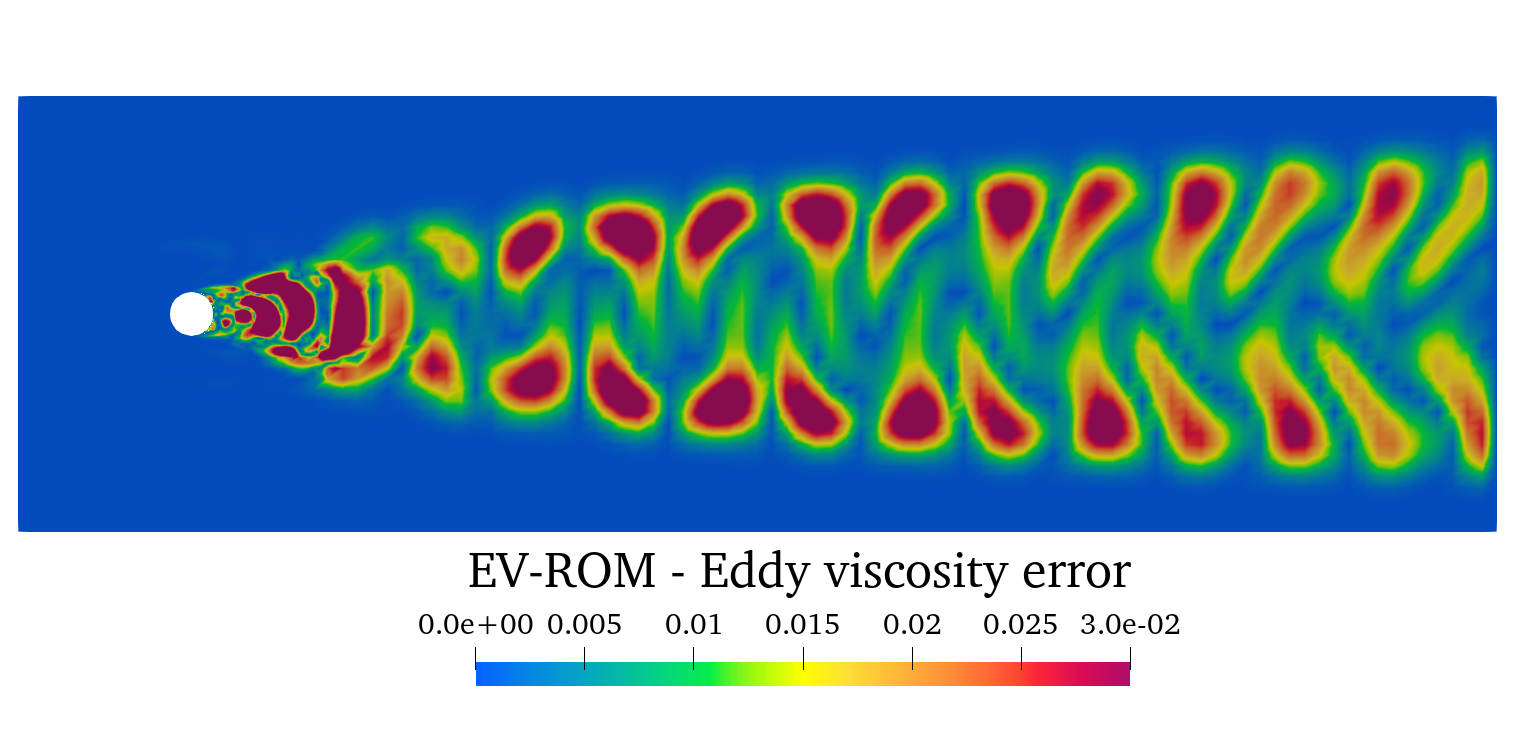}}
    \subfloat[]{
    \includegraphics[width=0.5\linewidth]{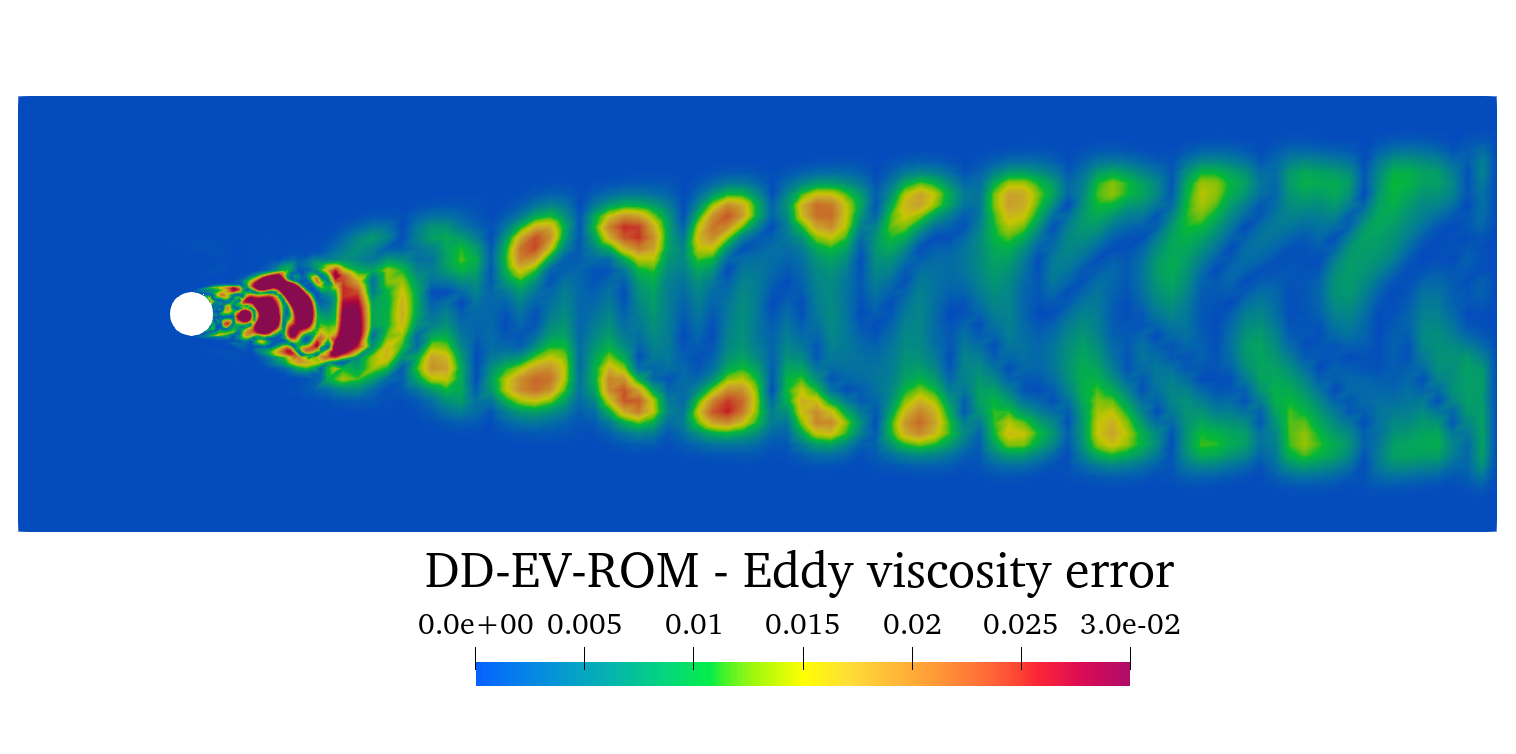}}\\
    \subfloat[]{
    \includegraphics[width=0.5\linewidth]{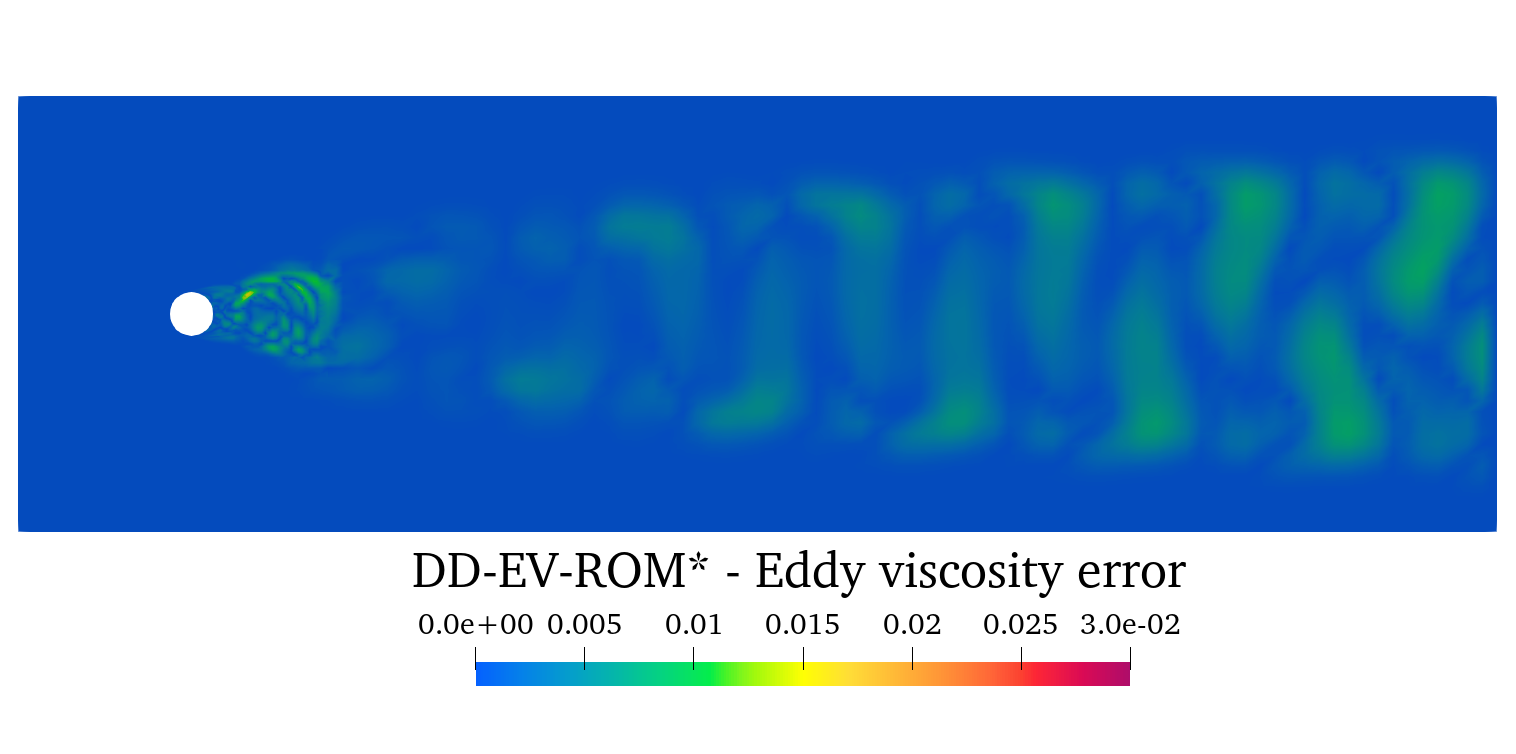}}
    \subfloat[]{
    \includegraphics[width=0.5\linewidth]{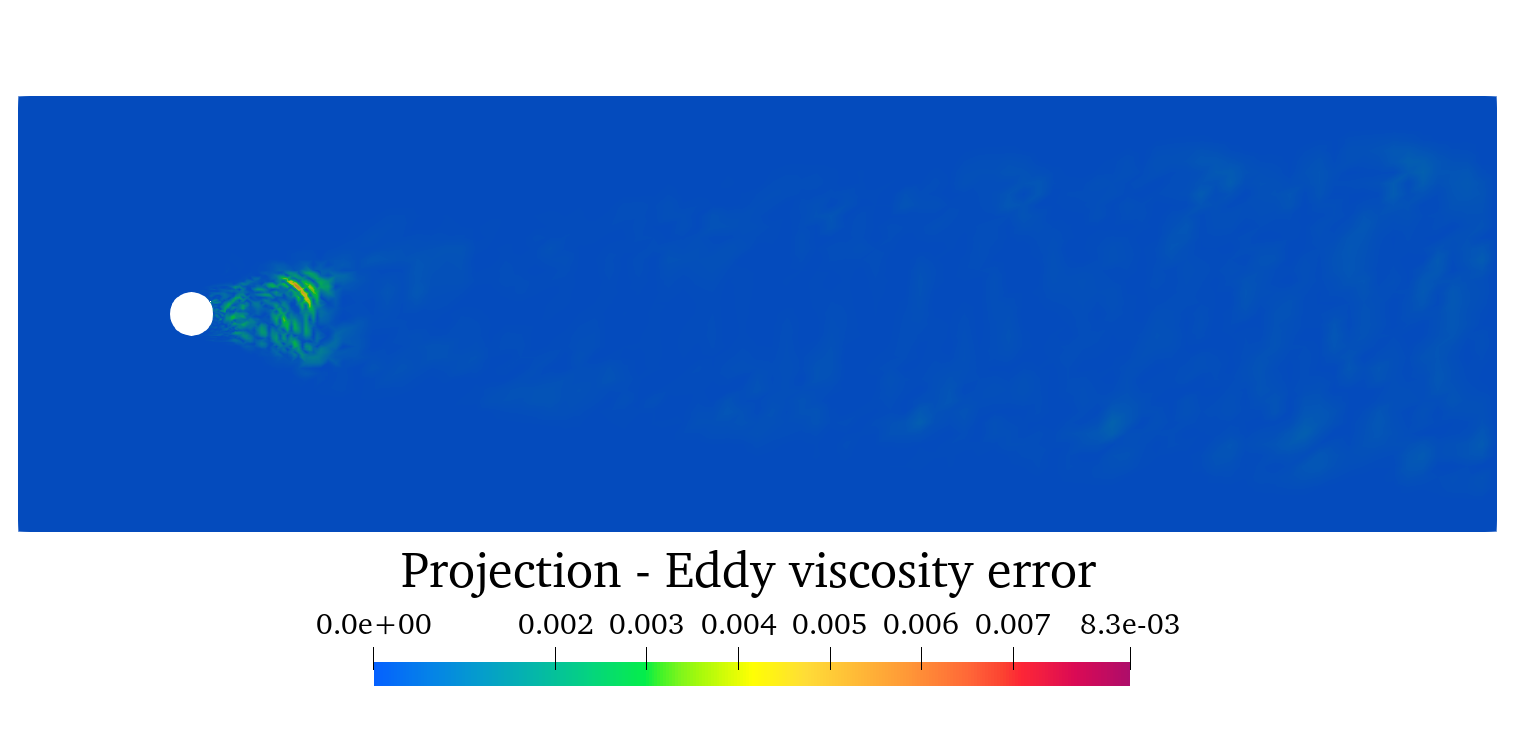}}
    \caption{Graphical eddy viscosity absolute error at $\nu=1.15e-4 \, \frac{m^2}{s}$ at the final time instance ($t=8$), of EV-ROM, DD-EV-ROM, DD-EV-ROM$^{\star}$, and for the projected field. The modal regime is $(N_u, N_p, N_{\nu_t})=(10, 24, 30)$.}
    \label{fig:cyl-graphical-nut-err}
\end{figure}

The graphical error fields of Figures \ref{fig:cyl-graphical-u-err} and \ref{fig:cyl-graphical-nut-err} highlith that the DD-EV-ROM$\star$ method significantly improved both the standard DD-EV-ROM and the baseline EV-ROM. It indeed leads to results close to the projected fields.

\subsection{Test case (\textbf{b})}
\label{app:case-b}

\begin{figure}[htpb!]
    \centering
    \includegraphics[width=\linewidth]{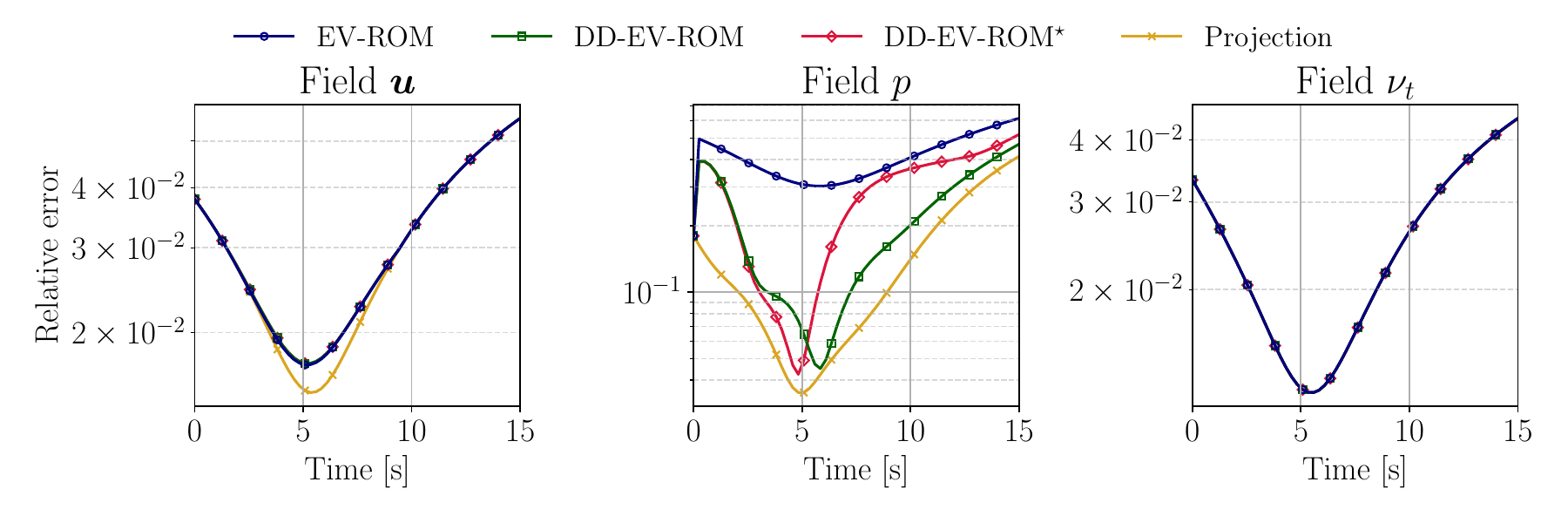}
    \caption{Time trends of the relative $L^2$ errors for the three fields of interest $\bm{u}$, $p$, and $\nu_t$, for one test viscosity $\nu=7e-6 \, \frac{m^2}{s}$. The modal regime is $(N_u, N_p, N_{\nu_t})=(1, 3, 1)$.}
    \label{fig:cavity-errs-1}
\end{figure}

\begin{figure}[htpb!]
    \centering
    \includegraphics[width=\linewidth]{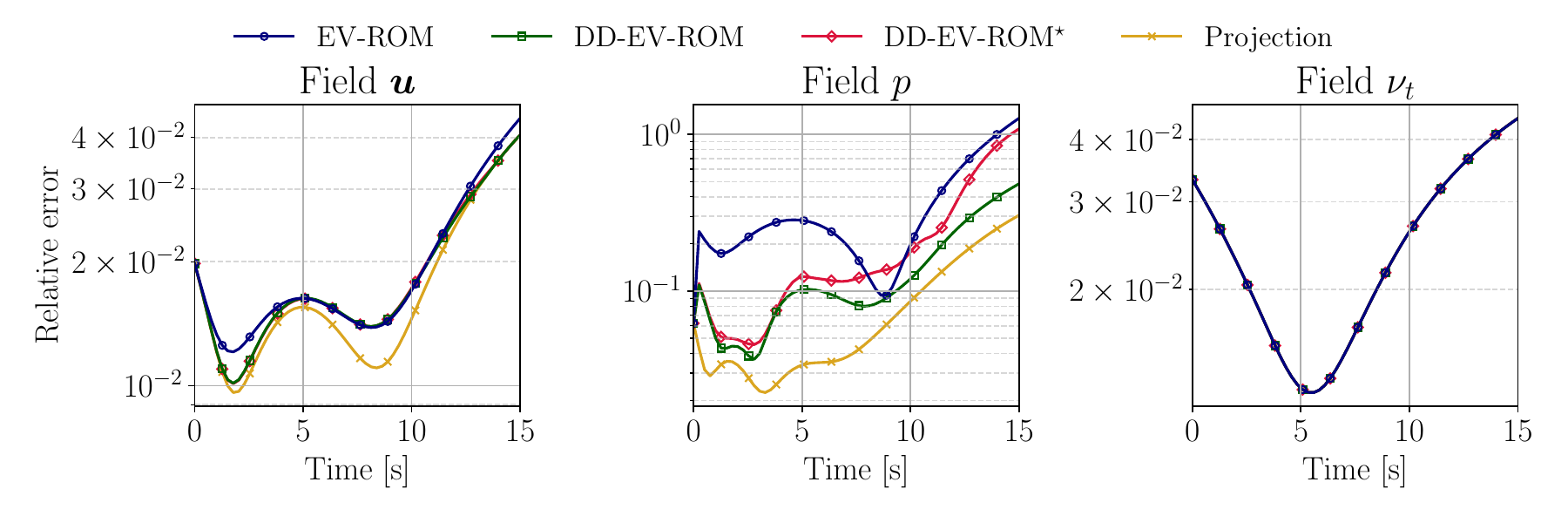}
    \caption{Time trends of the relative $L^2$ errors for the three fields of interest $\bm{u}$, $p$, and $\nu_t$, for one test viscosity $\nu=7e-6 \, \frac{m^2}{s}$. The modal regime is $(N_u, N_p, N_{\nu_t})=(2, 5, 1)$.}
    \label{fig:cavity-errs-2}
\end{figure}

\begin{figure}[htpb!]
    \centering{\includegraphics[width=\textwidth]{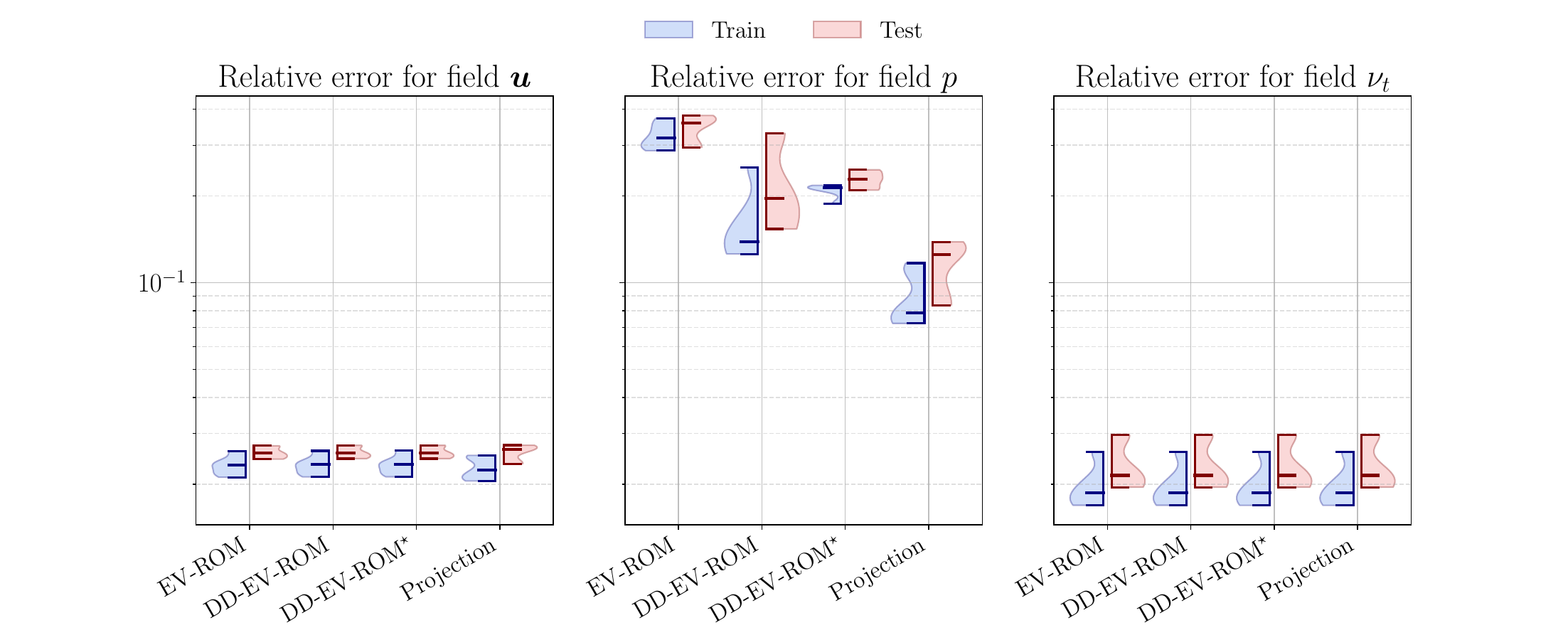}}
    \caption{Statistical performance for train and test viscosities in the time average of the relative $L^2$ error with respect to the high-fidelity solutions. The Figure represents the median and the error bounds for the three fields of interest in the case $(N_u, N_p, N_{\nu_t})=(1, 3, 1)$.}
    \label{fig:cavity-violin-1}
\end{figure}

\begin{figure}[htpb!]
    \centering{\includegraphics[width=\textwidth]{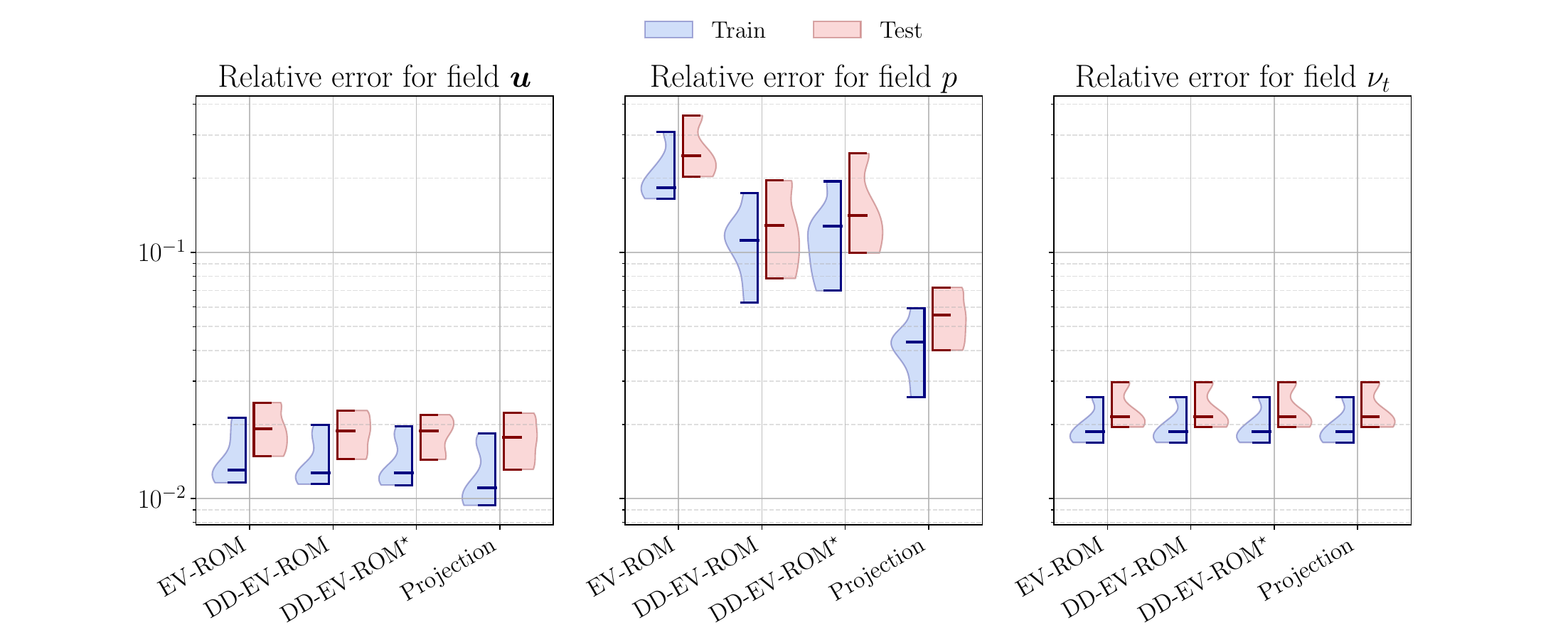}}
    \caption{Statistical performance for train and test viscosities in the time average of the relative $L^2$ error with respect to the high-fidelity solutions. The Figure represents the median and the error bounds for the three fields of interest in the case $(N_u, N_p, N_{\nu_t})=(2, 5, 1)$.}
    \label{fig:cavity-violin-2}
\end{figure}

For the second test case, we display the variation of the relative errors in time for one test viscosity value for two modal regimes (Figures \ref{fig:cavity-errs-1} and \ref{fig:cavity-errs-2}), and the statistical performance of the methods in train and test configurations in the same regimes (\ref{fig:cavity-violin-1} and \ref{fig:cavity-violin-2}).

The above-mentioned plots show that the EV-ROM, DD-EV-ROM, DD-EV-ROM$^{\star}$ predictions are almost exactly overlapped with the reconstruction error. As already mentioned in Subsection \ref{subsec:test-case-b}, the EV-ROM cannot be improved here since the projection error itself is the best accuracy we can reach.

Moreover, the Figures correspond to the regimes where the two methods have similar performances (see the heatmap in Figure \ref{fig:cavity-heatmaps}).
In such cases, DD-EV-ROM performs slightly better than the DD-EV-ROM$^{\star}$, but with comparable results.
This is also reflected in the time trends of Figures \ref{fig:cavity-errs-1} and \ref{fig:cavity-errs-2}, which show the ROM results in a time extrapolation setting, in the time window $[0, 15]$ seconds.
The time-dependent error trends confirm that the ROM velocity and eddy viscosity fields are close to the projection. However, for the pressure we can observe that the projection itself is characterized by large error, especially in the time interval $[10, 15]$ s, which is outside the time window considered to collect the snapshots.

\subsection{Test case (\textbf{c})}
\label{app:case-c}
\begin{figure}[htpb!]
    \centering{\includegraphics[width=\textwidth]{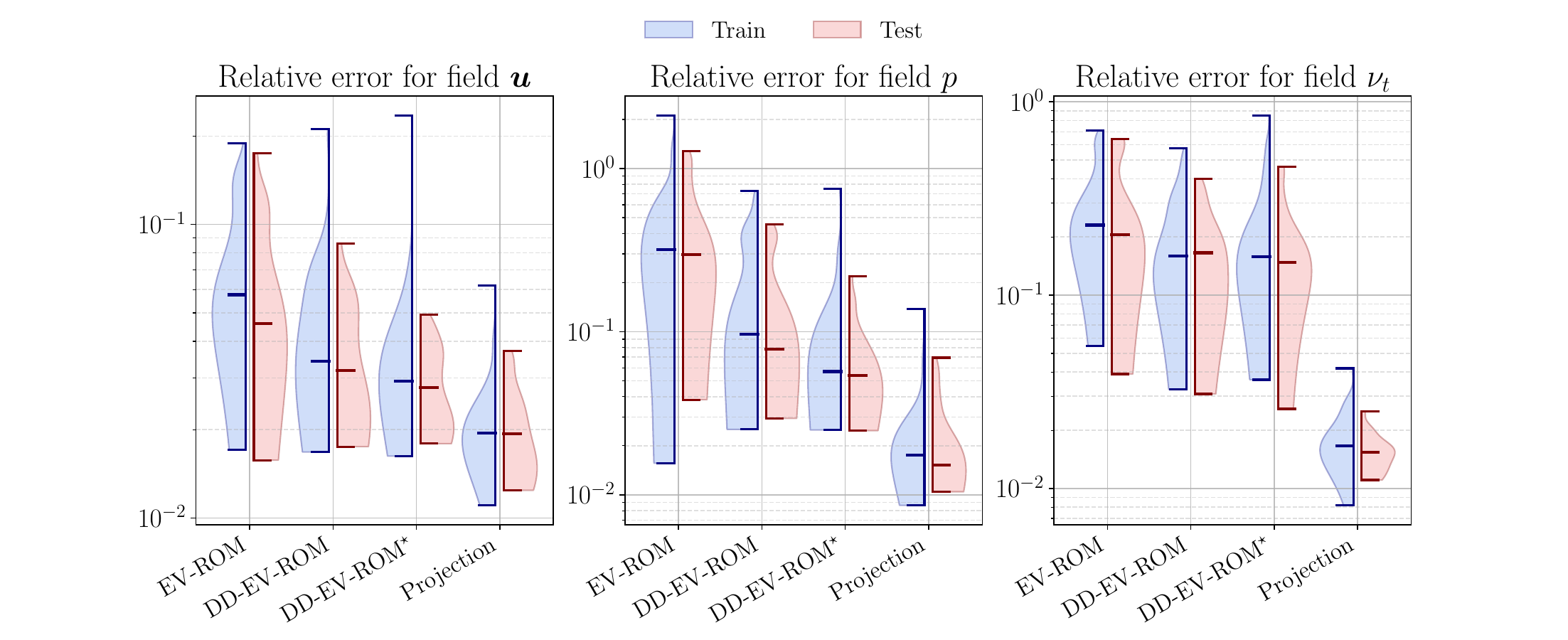}}
    \caption{Statistical performance for train and test parameters of the relative $L^2$ error with respect to the high-fidelity solutions. The Figure represents the median and the error bounds for the three fields of interest in the case $(N_u, N_p, N_{\nu_t})=(4, 5, 20)$.}
    \label{fig:backstep-violin-2}
\end{figure}

\begin{figure}[htpb!]
    \centering{\includegraphics[width=\textwidth]{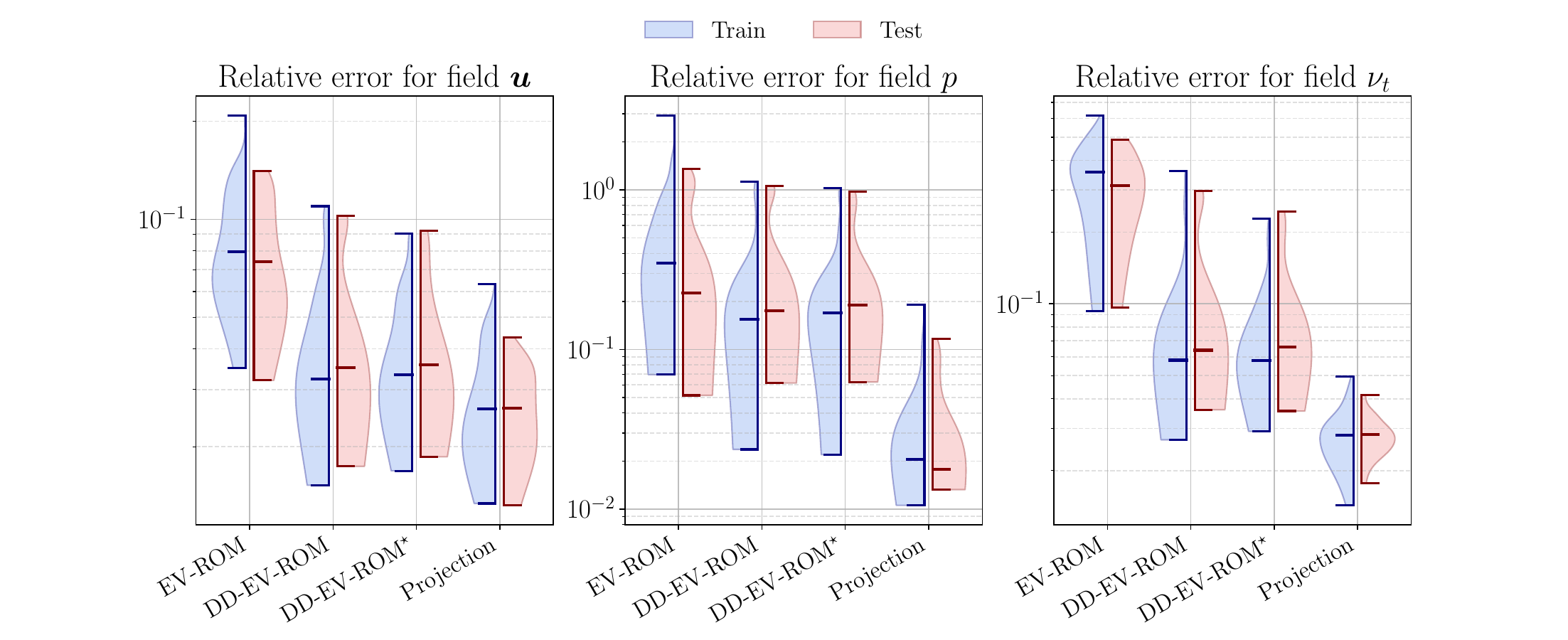}}
    \caption{Statistical performance for train and test parameters of the relative $L^2$ error with respect to the high-fidelity solutions. The Figure represents the median and the error bounds for the three fields of interest in the case $(N_u, N_p, N_{\nu_t})=(3, 4, 12)$.}
    \label{fig:backstep-violin-1}
\end{figure}

Since this test case is steady, we only perform the statistical train/test analysis for two modal regimes (Figures \ref{fig:backstep-violin-1} and \ref{fig:backstep-violin-2}).
As already noticed in Section \ref{subsec:test-case-c}, the DD-EV-ROM$^{\star}$ improves the standard DD-EV-ROM approach, especially in test configurations in the regime $(N_u, N_p, N_{\nu_t})=(4, 5, 20)$ (Figure \ref{fig:backstep-violin-2}).

\begin{figure}[htpb!]
    \centering
    \subfloat[]{
    \includegraphics[width=0.5\linewidth, trim={3cm 5cm 3cm 5cm}, clip]{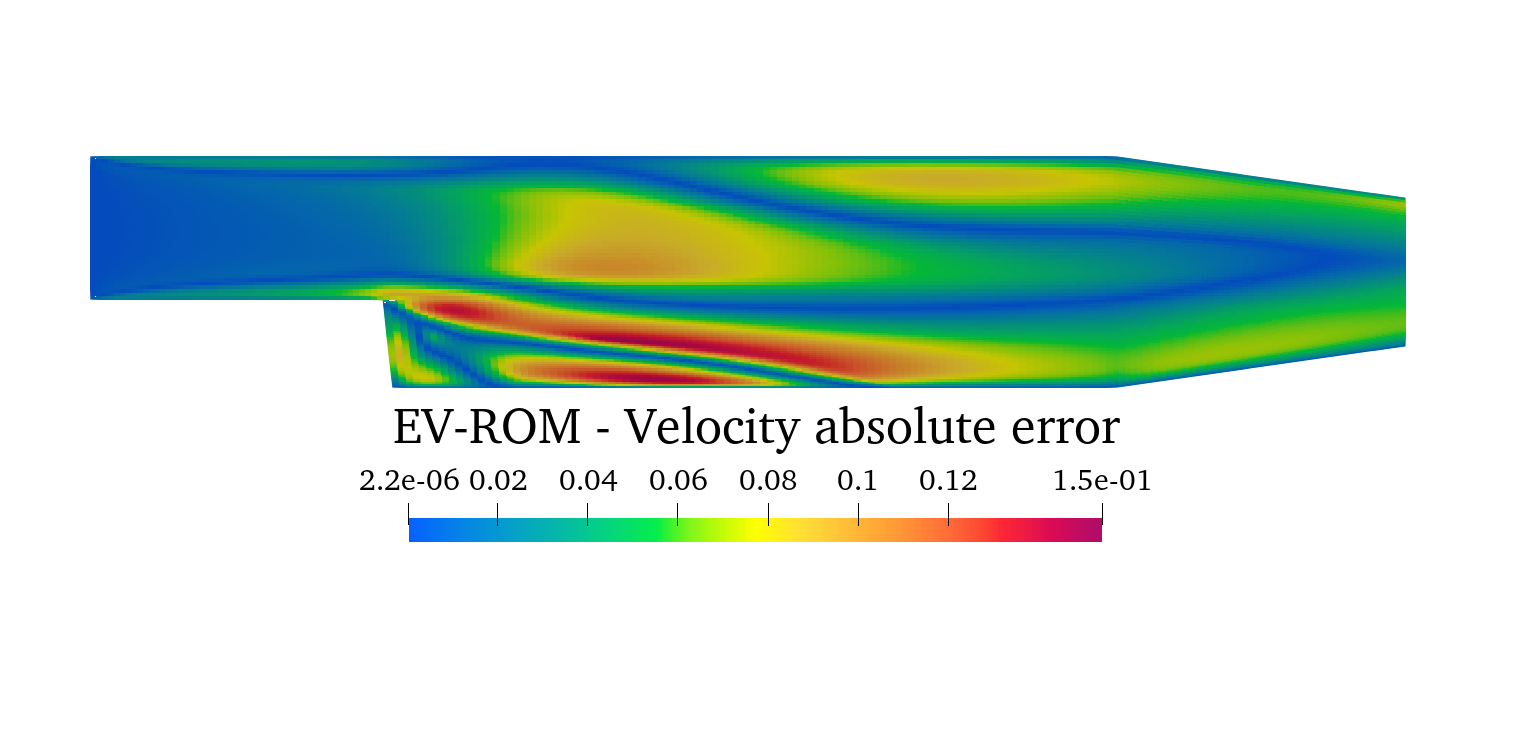}}
    \subfloat[]{
    \includegraphics[width=0.5\linewidth, trim={3cm 5cm 3cm 5cm}, clip]{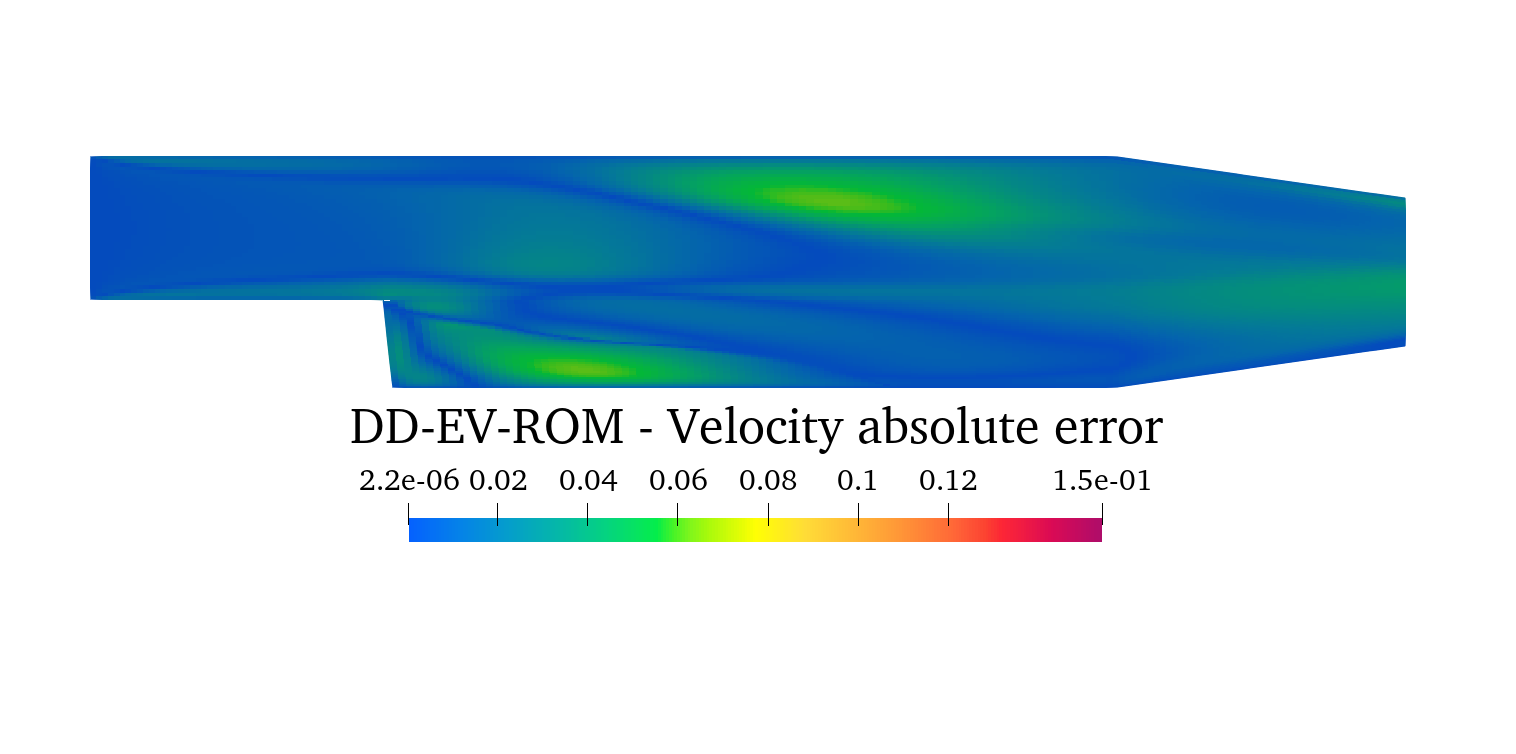}}\\
    \subfloat[]{
    \includegraphics[width=0.5\linewidth, trim={3cm 5cm 3cm 5cm}, clip]{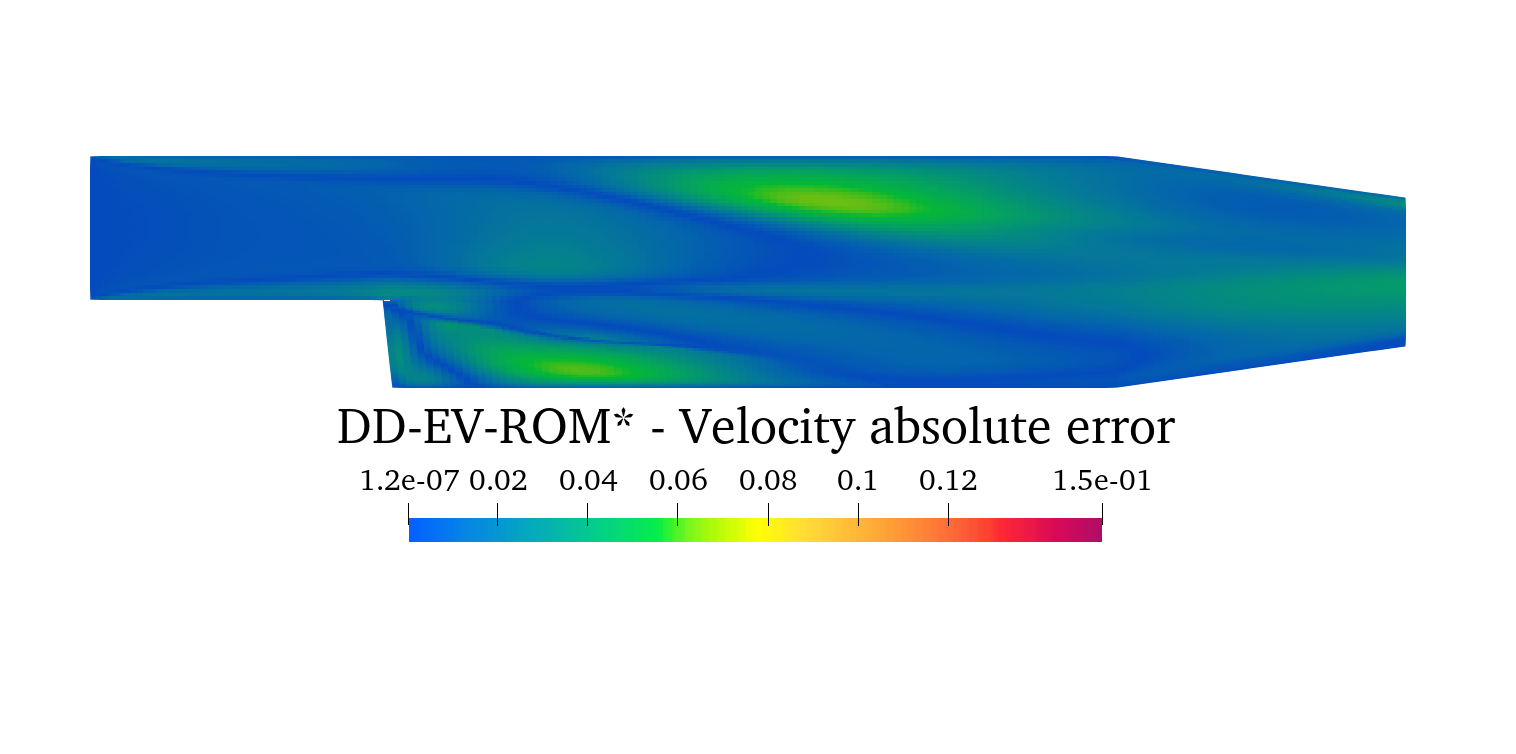}}
    \subfloat[]{
    \includegraphics[width=0.5\linewidth, trim={3cm 5cm 3cm 5cm}, clip]{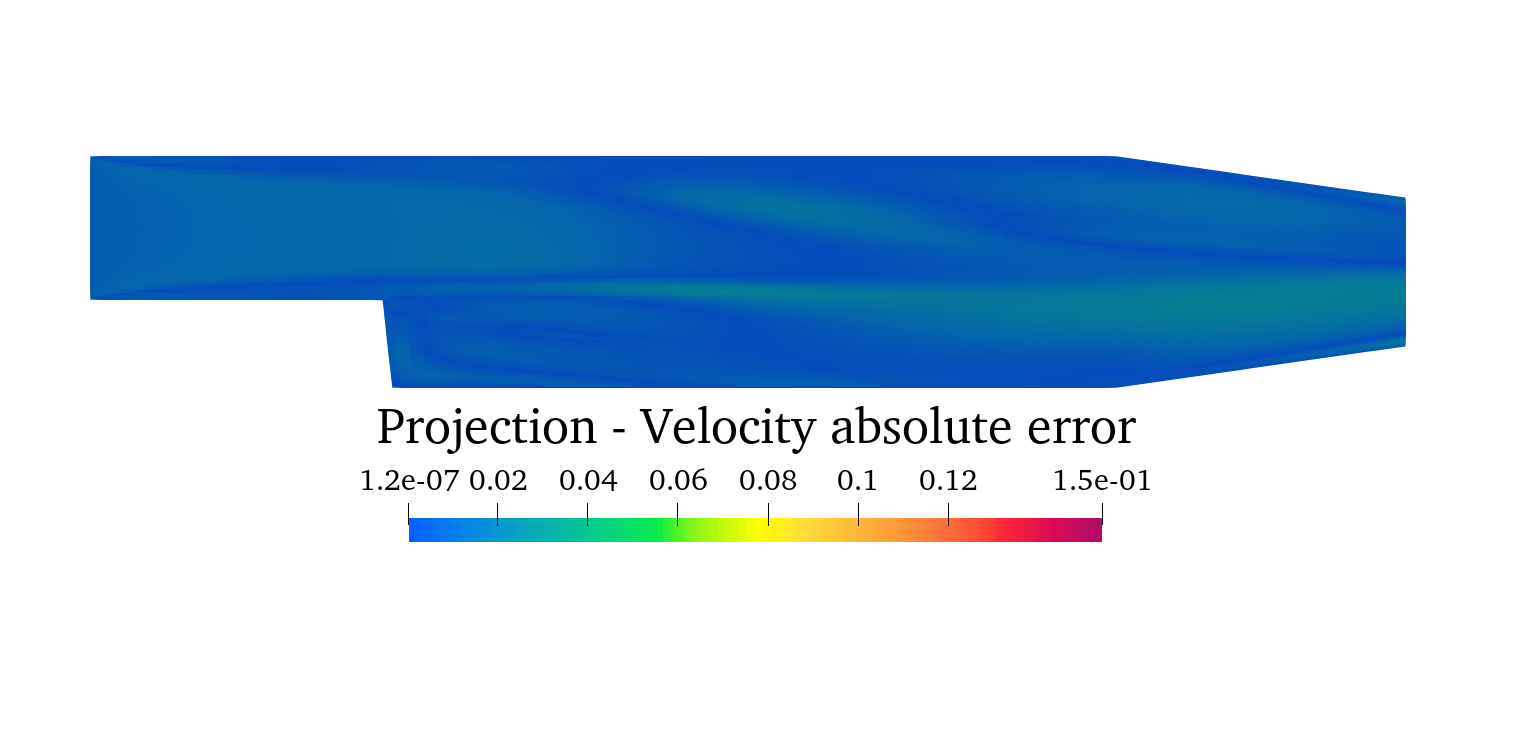}}
    \caption{Graphical absolute error for velocity magnitude field at parameters $(\alpha, h_1, h_2)=(6.3^{\circ}, 0.99 \text{ m}, 1.58\text{ m})$ of EV-ROM, DD-EV-ROM, DD-EV-ROM$^{\star}$, and for the projected field. The modal regime is $(N_u, N_p, N_{\nu_t})=(3, 4, 12)$.}
    \label{fig:backstep-graphical-u-err}
\end{figure}

\begin{figure}[htpb!]
    \centering
    \subfloat[]{
    \includegraphics[width=0.5\linewidth, trim={3cm 5cm 3cm 5cm}, clip]{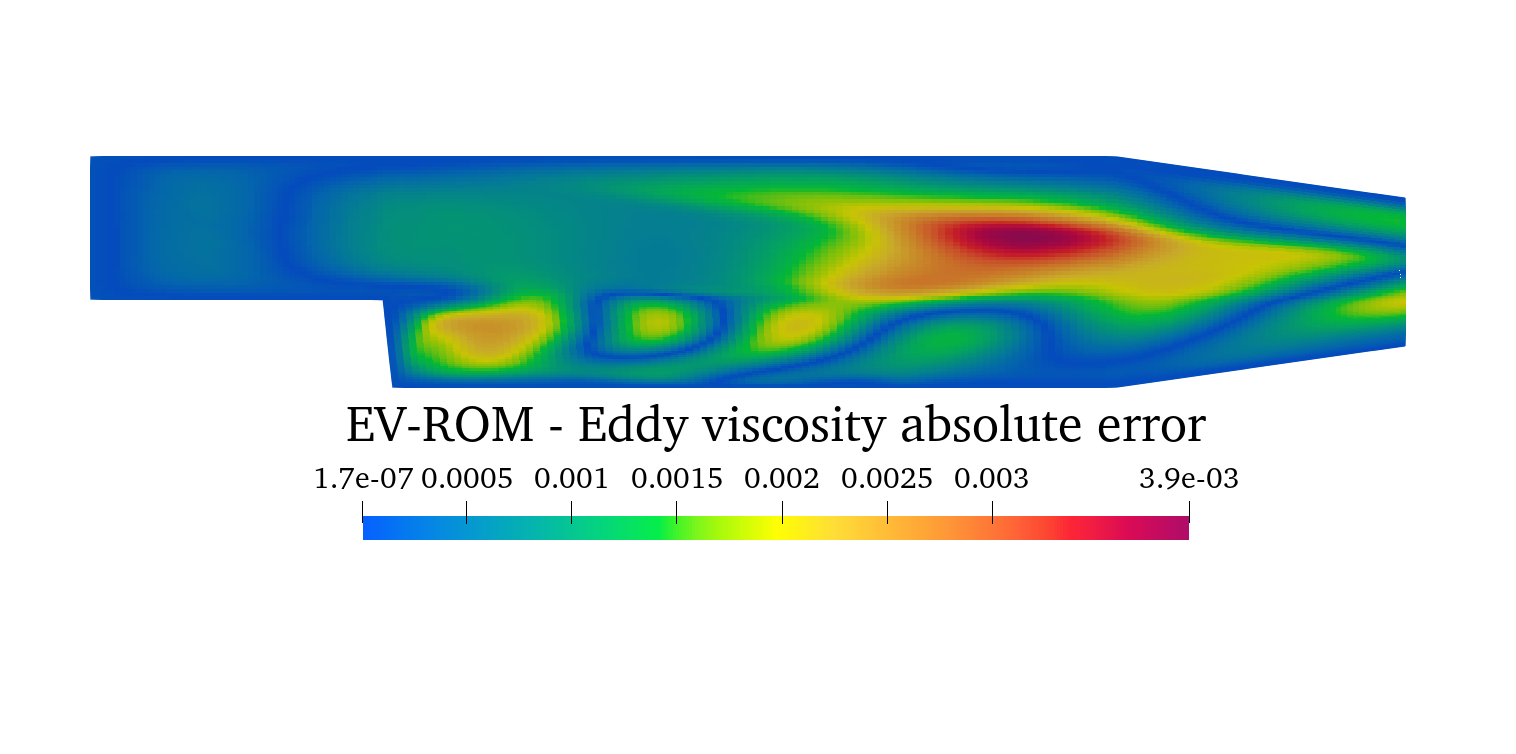}}
    \subfloat[]{
    \includegraphics[width=0.5\linewidth, trim={3cm 5cm 3cm 5cm}, clip]{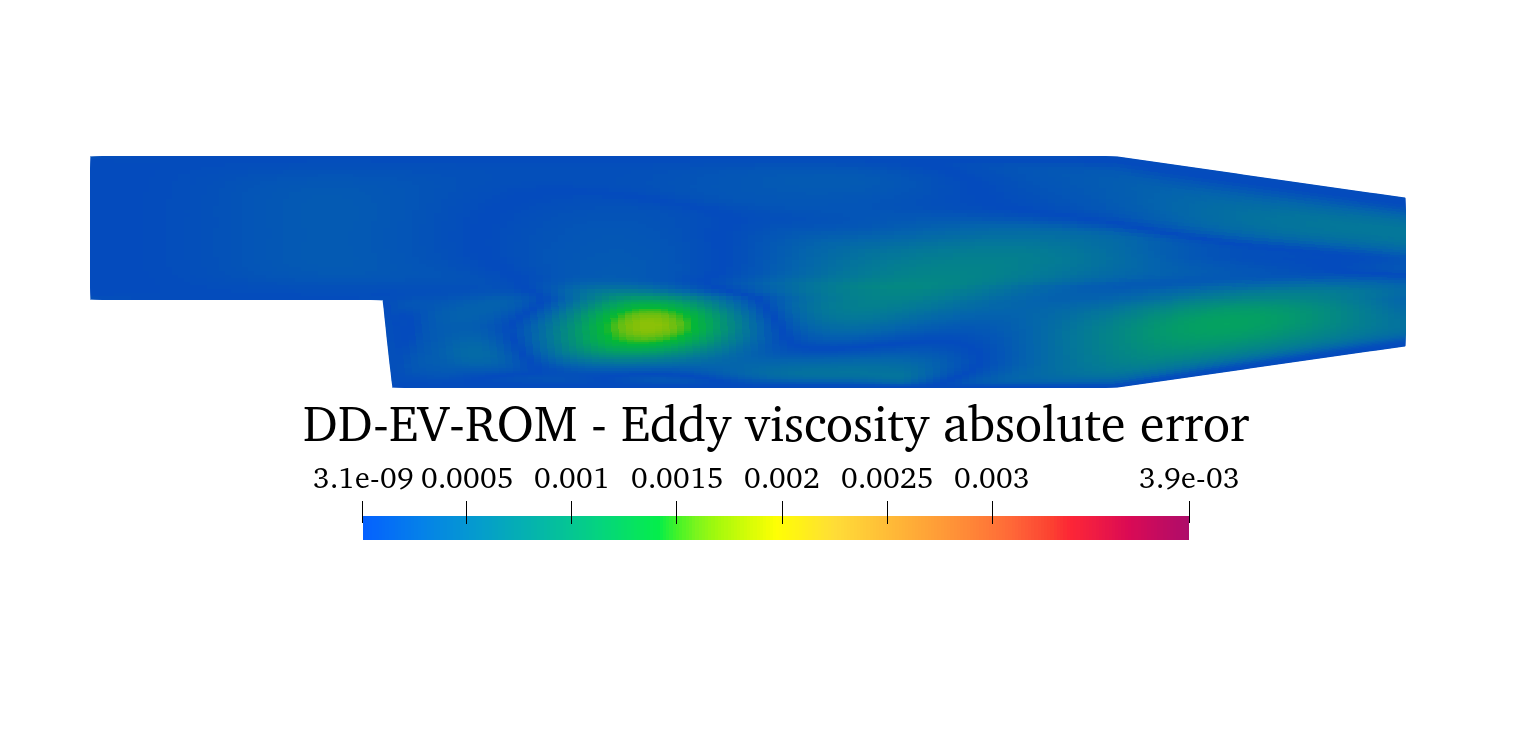}}\\
    \subfloat[]{
    \includegraphics[width=0.5\linewidth, trim={3cm 5cm 3cm 5cm}, clip]{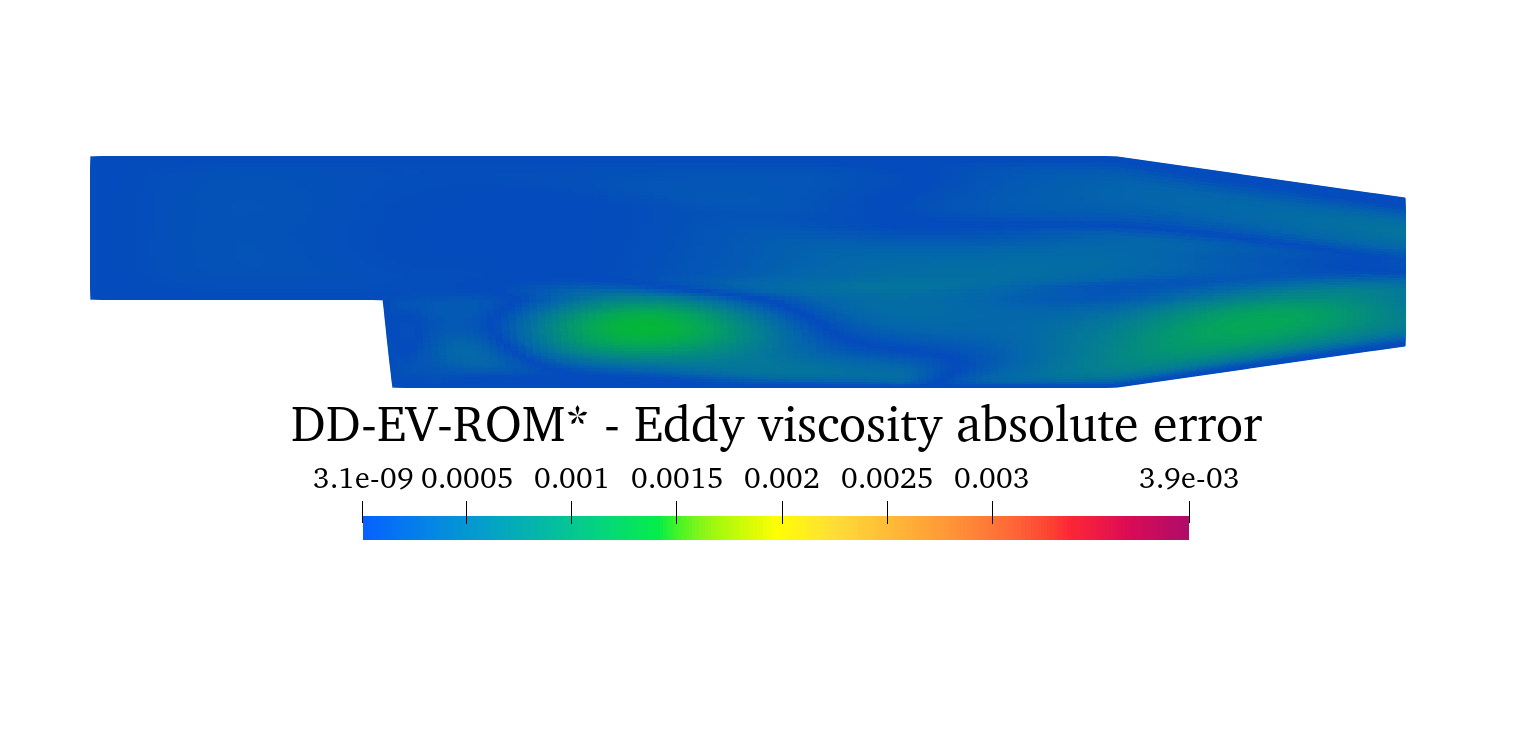}}
    \subfloat[]{
    \includegraphics[width=0.5\linewidth, trim={3cm 5cm 3cm 5cm}, clip]{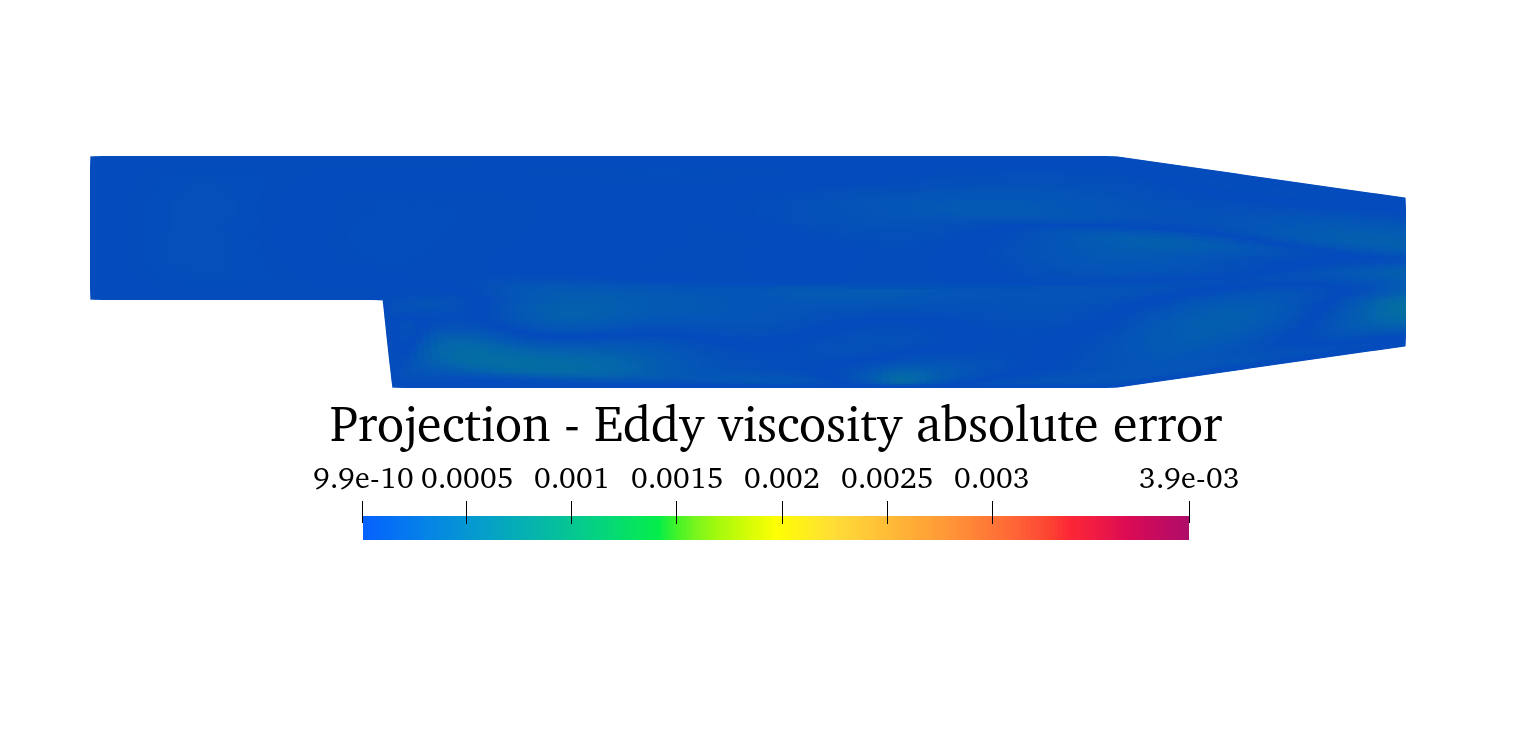}}
    \caption{Graphical absolute error for eddy viscosity field at the test parameters $(\alpha, h_1, h_2)=(6.3^{\circ}, 0.99 \text{ m}, 1.58\text{ m})$ of EV-ROM, DD-EV-ROM, DD-EV-ROM$^{\star}$, and for the projected field. The modal regime is $(N_u, N_p, N_{\nu_t})=(3, 4, 12)$.}
    \label{fig:backstep-graphical-nut-err}
\end{figure}

The overall improvement is confirmed by the graphical error fields in Figures \ref{fig:backstep-graphical-u-err} and \ref{fig:backstep-graphical-nut-err}, which focus on the combination $(N_u, N_p, N_{\nu_t})=(3, 4, 12)$.

\printnomenclature
\label{appendix-nomenclature}

\end{document}